\documentclass[reqno,10pt]{amsart} 
\usepackage[colorlinks = true, linkcolor=blue,
            urlcolor=red,
            citecolor=olive]{hyperref}
\usepackage[utf8]{inputenc}
\usepackage[T1]{fontenc}
\usepackage{amsmath}

\usepackage{mathtools}      
\usepackage{mathabx}        

\usepackage{tikz}
\usepackage{tikz-3dplot}

\usepackage{dynkin-diagrams}
\usepackage{array, multirow, makecell}

\usepackage{amssymb}
\usepackage{amsfonts}
\usepackage{graphicx}

\usepackage{multicol}

\usepackage[colorinlistoftodos]{todonotes}

\tikzset{surface/.style={draw=blue!70!black, fill=blue!40!white, fill opacity=.6}}
\tikzset{circ/.style={shape=circle, inner sep=1.5pt, draw, node contents=}}

\hypersetup{bookmarksdepth  =  3} 

\theoremstyle{plain}
\newtheorem{thm}{Theorem}[section]
\newtheorem{lem}[thm]{Lemma}
\newtheorem{prop}[thm]{Proposition}
\newtheorem{cor}[thm]{Corollary}

\theoremstyle{definition}
\newtheorem{defn}[thm]{Definition}
\newtheorem{que}[thm]{Question}

\newtheorem{exmp}[thm]{Example}

\theoremstyle{remark}
\newtheorem{rem}[thm]{Remark}
\newtheorem{rems}[thm]{Remarks}

\newcommand{\type}{\tau}
\newcommand{\typeS}{{\tau_0}}
\newcommand{\R}{\mathbb{R}}
\newcommand{\C}{\mathbb{C}}
\newcommand{\N}{\mathbb{N}}
\newcommand{\Sph}{\mathbb{S}}
\newcommand{\K}{\mathbb{K}}

\newcommand{\Id}{\text{Id}}
\newcommand{\Tangle}{\angle_\text{Tits}}
\newcommand{\transpose}[1]{\ensuremath{#1^{\scriptscriptstyle T}}}
\DeclareMathOperator{\vis}{vis}
\newcommand{\partialvis}{\partial_{\vis}}

\newcommand{\da}{\mathbf{d}_{\mathfrak{a}}}

\DeclareMathOperator{\Cartan}{\boldsymbol{\mu}}
\DeclareMathOperator{\End}{End}
\DeclareMathOperator{\tr}{Tr}
\DeclareMathOperator{\dev}{dev}
\DeclareMathOperator{\hol}{hol}
\DeclareMathOperator{\Tol}{Tol}
\DeclareMathOperator{\Hom}{Hom}
\DeclareMathOperator{\PSL}{PSL}
\DeclareMathOperator{\SL}{SL}
\DeclareMathOperator{\GL}{GL}
\DeclareMathOperator{\PGL}{PGL}
\DeclareMathOperator{\PSO}{PSO}
\DeclareMathOperator{\PSU}{PSU}
\DeclareMathOperator{\SO}{SO}
\DeclareMathOperator{\Sp}{Sp}
\DeclareMathOperator{\PSp}{PSp}

\DeclareMathOperator{\Diag}{Diag}
\DeclareMathOperator{\Hess}{Hess}
\newcommand{\II}{\mathrm{I\!I}}
\newcommand{\Ein}{\mathbf{Ein}}
\newcommand{\Pho}{\mathbf{Pho}}
\DeclareMathOperator{\ad}{ad}
\DeclareMathOperator{\Ad}{Ad}
\DeclareMathOperator{\I}{I}
\DeclareMathOperator{\Ker}{Ker}
\DeclareMathOperator{\rank}{rank}

\newcommand{\X}{\mathbb{X}}

\newcommand{\ie}{namely }
\title{Nearly geodesic immersions and domains of discontinuity}

\copyrightinfo{2022}{Colin Davalo}

\author{Colin Davalo}
\address{Mathematisches Institut, Ruprecht-Karls Universität Heidelberg, Im Neuenheimer Feld 205, 69120 Heidelberg, Germany}
\email{cdavalo@mathi.uni-heidelberg.de}

\date{\today}
\begin{document}
\begin{abstract}
We study nearly geodesic immersions in higher rank symmetric spaces of non-compact type, which we define as immersions that satisfy a bound on their fundamental form, generalizing the notion of immersions in hyperbolic space with principal curvature in $(-1,1)$. This notion depends on the choice of a flag manifold embedded in the visual boundary, and immersions satisfying this bound admit a natural domain in this flag manifold that comes with a fibration. As an application we give an explicit fibration of some domains of discontinuity for some Anosov representations. Our method can be applied in particular to some $\Theta$-positive representations for each notion of $\Theta$-positivity.


\end{abstract}

\maketitle

\tableofcontents

\section{Introduction.}

Let $\Gamma_g$ be the fundamental group of a closed orientable surface $S_g$ of genus $g\geq 2$. An interesting open subset of the space of representations $\rho:\Gamma_g\to \PSL(2,\C)$ is the space of \emph{quasi Fuchsian} representations, which are representations whose limit set is a Jordan curve, or equivalently representations that are quasi-isometric embeddings.

\medskip

An interesting open subset of this space is the set of \emph{nearly Fuchsian} representations, which is the set of representations that admit a $\rho$-equivariant immersion from the universal cover of $S_g$ to the hyperbolic space $\mathbb{H}^3$ whose image has principal curvature in $(-1,1)$. Epstein studied in \cite{Epstein} immersions satisfying this bound. He proved in particular that these immersions must be proper embeddings and that their Gauss map fibers a domain in $\partial \mathbb{H}^3$.
 Nearly Fuchsian representations must therefore be discrete, and one can even conclude that they are quasi-isometric embeddings. A nearly Fuchsian representation is called \emph{almost Fuchsian} if the immersion $u$ can be chosen to be a minimal immersion, these representations have been studied among others by Uhlenbeck \cite{Uhlenbeck}.
 
 \medskip
 
A generalization of quasi Fuchsian representations for higher rank semi-simple Lie groups are \emph{Anosov representations}, introduced by Labourie \cite{Labourie} for representations of surface groups and by Guichard and Wienhard \cite{GWDod} for representations of any Gromov hyperbolic group. Given a subset $\Theta$ of the set of simple restricted roots of the Lie group $G$, a $\Theta$-Anosov representation is a prepresentation that satisfies some strong contraction property in the associated flag manifold, and it is in particular discrete. The set of $\Theta$-Anosov representations $\rho:\Gamma\to G$ is moreover open.

\medskip

In this paper we consider a generalization of Epstein's bound on the second fundamental form of an immersion in a symmetric space of non-compact type $\X$ associated to a semi-simple Lie group $G$. Given a unit vector $\type$ in the model Weyl chamber, or equivalently a flag manifold $\mathcal{F}_\type$ embedded in the visual boundary of $\X$, we define the notion of \emph{$\type$-nearly geodesic} and  \emph{uniformly $\type$-nearly geodesic} immersions using Busemann functions.
We prove that a \emph{uniformly $\type$-nearly geodesic} immersion is an embedding, a quasi-isometric embedding, and admits a domain in the corresponding flag manifold $\mathcal{F}_\type$, which admits a fibration whose fibers are the base of some \emph{pencils of tangent vectors} in $\X$, which is a notion generalizing pencils of quadrics that we introduce.

\medskip

Given the fundamental group $\Gamma$ of a compact manifold $N$ we also study \emph{$\type$-nearly Fuchsian} representations, namely representations $\rho:\Gamma\to G$ that admit an equivariant and $\type$-nearly geodesic immersion $u:\widetilde{N}\to \X$. We prove that they form an open subset of the space Anosov representations.

\medskip

As an application we study some domains of discontinuity associated to Anosov representations in flags manifolds, constructed by Guichard and Wienhard \cite{GWDod} and Kapovich, Leeb and Porti \cite{KLP-DoD}. Dumas and Sanders conjectured in \cite{DumasSanders} that, in the context of complex Lie groups, the quotient of the domain of discontinuity associated to a Tits-Bruhat ideal fibers equivariantly over the universal cover of the surface $S_g$ for every representation of a surface group $\Gamma_g$ that admit a totally geodesic equivariant immersion. 
We show that for a nearly Fuchsian representation $\rho$ of the fundamental group of $N$, the domain associated with the immersion $u$ coincides with a domain of discontinuity for $\rho$, and fibers equivariantly over the universal cover of $N$. This technique can be applied to representations of surface groups, and in particular to some $\Theta$-positive representations for every notion of $\Theta$-positivity. One can also apply these techniques to fundamental groups of hyperbolic manifolds of higher dimension: we consider two examples in the end of Section \ref{subsec:HigherDimention}.

\medskip

Independently and with different techniques, Alessandrini, Maloni, Tholozan and Wienhard proved that any domain of discontinuity associated to any representation that factor through a cocompact representation into a rank $1$ simple Lie group admits an equivariant fibration \cite{AMTW}. However the fibers are not easy to describe in general.

\subsection{Nearly geodesic immersions.}

We introduce and study a generalization of the condition of having principal curvature in $(-1,1)$ in the setting of higher rank semi-simple Lie groups of non-compact type $G$. More specifically we choose a unit vector $\tau$ in the associated model Weyl chamber $\mathfrak{a}^+$. This choice is equivalent to the choice of a flag manifold $\mathcal{F}_\type$ embedded in the visual boundary $\partialvis \X$ of the symmetric space $\X$ associated to $G$. Let $\iota:\Sph\mathfrak{a}^+\to \mathfrak{a}^+$ be the antipodal involution, to a point $a\in \mathcal{F}_\type\cup \mathcal{F}_{\iota(\type)}\subset \partialvis \X$ and a base-point $o\in \X$ one can associate a \emph{Busemann function} $b_{a,o}$ which can be interpreted as the distance to $a$ in $\X$, relative to $o$. 

\begin{defn}[Definition \ref{Defn:NearlyGeodesicSurface}]
\label{Defn:NearlyGeodesicSurfaceIntro}
An immersion $u: M\to \X$ is called \emph{$\type$-nearly geodesic} if for all $a\in \mathcal{F}_\type\cup \mathcal{F}_{\iota(\type)}$  the function $b_{a,o}\circ u: M\to \R$ has positive Hessian in any critical direction $\mathrm{v}\in TM$, for the metric on $M$ induced by the immersion.
\end{defn}

This property is satisfied for totally geodesic immersions whose tangent vectors are $\type$-regular, namely whose Cartan projection is not orthogonal to $w\cdot \type$ for any $w$ in the Weyl group (see Definition \ref{defn:pRegularVector}). It is equivalent to an open bound on the second fundamental form that depends on the Cartan projection of the image of the differential of the immersion (see Proposition \ref{prop:nearlyfuchsianLocalProp}). When $G=\PSL(2,\C)$ a $\type$-nearly geodesic immersion for the only $\type\in\Sph\mathfrak{a}^+$ is exactly an immersion with sectional curvature in $(-1,1)$ in $\X=\mathbb{H}^3$ (see Proposition \ref{prop:nearlyfuchsianH3}).

\medskip

If the immersion is complete and \emph{uniformly} $\type$-nearly geodesic, namely if the Hessian of Busemann functions in critical directions are uniformly bounded from below, we show moreover that it is a $\tau$-regular embedding, a quasi isometric embedding (see Proposition \ref{prop:QuasiGeod}) and that the nearest point projection $\pi^\type_u\circ u:\X\to u(M)$ is well defined for some explicit appropriate Finsler pseudo distance on $\X$ depending on $\type$ (see Proposition \ref{prop:ProjectionFinsler}).

\medskip

If a $\type$-nearly geodesic immersion is equivariant with respect to a representation of a group acting cocompactly and properly on $M$, it provides information about the representation.

\begin{thm}[Theorem \ref{thm:NearlyfuchsianAnosov}]
\label{thm:NearlyfuchsianAnosovIntro}
A $\type$-nearly Fuchsian representation $\rho:\Gamma\to G$ is $\Theta$-Anosov for some non-empty set of roots $\Theta$ that depends on $\type$ and $\rho$.
\end{thm}

This theorem can be stated in a simpler form if we consider a vector $\tau\in \Sph\mathfrak{a}^+$ colinear to a coroot. Given a simple Lie group, the set $\Delta$ of simple roots can be partitioned into one or two \emph{Weyl orbits of simple roots}, see Figure \ref{fig:TableWeylOrbit}. To a Weyl orbit of simple roots $\Theta$ one can associate a unique normalized coroot $\type_\Theta\in \Sph\mathfrak{a}^+$ in the model Weyl chamber. The previous theorem can be restated in this case.

\begin{thm}[Theorem \ref{thm:NearlyfuchsianAnosovNormal}]
\label{thm:NearlyfuchsianAnosovIntro2}
A $\type_\Theta$-nearly Fuchsian representation for some Weyl orbit of simple roots $\Theta\subset \Delta$ is $\Theta$-Anosov.
\end{thm}

When $\type=\type_\Theta$ we also prove a sufficient condition for a surface to be $\type$-nearly geodesic.

\begin{thm}[Theorem \ref{thm:Sufficient condition}]
Let $\Theta$ be a Weyl orbit of simple roots. Let $u:S\to \X$ be an immersion that satisfies for all $\mathrm{v}\in TS$ and $\alpha\in \Theta$:
\begin{equation}\label{eq:AlmostfuchsianIntro}\lVert\II_u(\mathrm{v},\mathrm{v}) \rVert<c_\Theta \alpha\left(\Cartan\left(\mathrm{d}u(\mathrm{v})\right)\right)^2.
\end{equation}
Then $u$ is a $\type_\Theta$-nearly geodesic immersion.
\end{thm}

Here $\Cartan:T\X\to \mathfrak{a}^+$ denotes the \emph{Cartan projection}. The constant $c_\Theta$ depends on the scaling of the metric chosen on $\X$, and on $\Theta$. 


\subsection{An associated domain in a flag manifold.}
In Section \ref{sec:Pencils} we introduce and study \emph{pencils of tangent vectors}, or $d$-pencils, which are vector subspaces $\mathcal{P}\subset T_x\X$ of dimension $d$ for some $x\in \X$. When $G=\PSL(n,\R)$ these can be thought of as \emph{pencils of quadrics} with zero trace with respect to some scalar product (see Proposition \ref{prop:PencilQuadrics}). To a pencil we associate a subset of the flag manifold $\mathcal{F}_\type$ that we call its \emph{base}, which is a smooth submanifold if the pencil is $\type$-regular, \ie if all non-zero vectors $\mathrm{v}\in \mathcal{P}$ are $\tau$-regular, see Lemma \ref{lem:SmoothSymmetricPencils}. When $G=\PSL(n,\R)$ and $\mathcal{F}_\type \simeq \mathbb{RP}^{n-1}$ the base of the pencil corresponds to the intersection of all the quadric hypersurfaces defined by the elements of the corresponding pencil of quadrics.

\medskip

To a complete and uniformly $\type$-nearly geodesic immersion $u:M\to \X$, we associate an open domain $\Omega^\type_u\subset \mathcal{F}_\type$ consisting of points $a\in \mathcal{F}_\type$ for which $b_{a,o}\circ u$ is proper and bounded from below for one and hence any base-point $o\in \X$.
We define a projection $\pi_u:\Omega^\type_u\to M$ that associates to $a\in \Omega^\type_u$ the point in $M$ at which $b_{a,o}\circ u$ is minimal. This point will be unique because $b_{a,o}\circ u$ is convex and strictly convex in critical directions.

\begin{thm}[Theorem \ref{thm:FoliationDoD}]
\label{thm:FoliationDoDIntro}

Let $u:M\to \X$ be a complete and uniformly $\type$-nearly geodesic immersion. The map $\pi_u:\Omega^\type_u\to M$ is a fibration. The fiber $\pi_u^{-1}(x)$ at a point $x\in M$ is the base $\mathcal{B}_\type(\mathcal{P}_x)$ of the $\type$-regular pencil of tangent vectors $\mathcal{P}_x=\mathrm{d}u(T_xM)$.

\end{thm}

A \emph{domain of discontinuity} $\Omega$ for a representation $\rho:\Gamma\to G$ in a $G$-homogeneous space is an open, $\rho(\Gamma)$-invariant set on which $\rho(\Gamma)$ acts properly discontinuously\footnote{Some authors call these domains of \emph{proper} discontinuity.}. A \emph{cocompact} domain of discontinuity is a domain on which the action of $\rho(\Gamma)$ is cocompact.

\medskip

If the representation $\rho:\Gamma\to G$ admits an equivariant $\type$-nearly geodesic immersion $u:\widetilde{N}\to \X$, the domain $\Omega_\rho^\type:=\Omega_u^\type$ does not depend on $u$, and is a cocompact domain of discontinuity for $\rho$. The fibration $\pi_u$ from Theorem \ref{thm:FoliationDoDIntro} is $\rho$-equivariant so the quotient of $\Omega_\rho^\type$ fibers over the surface.
We show in Theorem \ref{thm:ComparaisonTotallyGeodesic} that the domain $\Omega^\type_\rho$ coincides with some domains of discontinuity associated to \emph{Tits-Bruhat ideals} for Anosov representations, constructed by Kapovich, Leeb and Porti \cite{KLP-DoD}. The Tits Bruhat ideals needed to describe the domains $\Omega^\type_\rho$ can always be constructed as a \emph{metric thickening}.

\medskip

We prove in Section \ref{subsec:Ehresman Principle} that the diffeomorphism type of the domains of discontinuity obtained by Tits-Bruhat ideals is invariant under deformation in the space of Anosov representations. We therefore understand the topology of some domains of discontinuity for all representations in some connected component of $\Theta$-Anosov representation.

\subsection{Applications.}

%

A union of connected components of representations of a surface group into a higher rank simple Lie groups $G$ consisting only of discrete representations is called a \emph{higher rank Teichmüller space} \cite{WInvitation}. The Hitchin component is a higher rank Teichmüller space for any split real simple Lie group. When $G$ is of Hermitian type, the space of \emph{maximal} representations, \ie representations with maximal Toledo invariant, is also a higher rank Teichmüller spaces. More generally, Guichard and Wienhard defined a notion of
$\Theta$-positive representations for some simple Lie groups and some subsets of simple roots $\Theta$, of which Hitchin and maximal representations are a
particular instance. Beyer and Pozzetti proved that the spaces of $\Theta$-positive representations
in $\SO(p,q)$ for $3\leq p < q$ are higher rank Teichmüller spaces \cite{BeyerPozzetti}.
Bradlow, Collier, Garcia-Prada,  Gothen, and Oliveira \cite{magical} and Guichard, Labourie and Wienhard \cite{positivityLabourie} proved that for each notion of $\Theta$-positivity there exist higher rank Teichmüller spaces consisting of $\Theta$-positive representations. 

\medskip

The classical Teichmüller space can be defined a connected component of the space of conjugacy classes of discrete and faithful representations of a surface group into $\PSL(2,\R)$. It can also be interpreted as the space of marked hyperbolic structures on the corresponding surface, a particular instance of $(G,X)$-structures in the sense of Klein. 
To a manifold $M$ equipped with a $(G,X)$-structure one can associate its \emph{holonomy representation} $\rho:\pi_1(M)\to G$, defined up to conjugation by elements of $G$.

\medskip

Hitchin underlined similarities between the Teichmüller space and the Hitchin component, and asked about the geometric significance of elements in the Hitchin component. In the spirit of this question it is natural to ask the following:

\begin{que}
Can $\Theta$-postitive representations $\rho:\Gamma_g\to G$ be characterized as the holonomies on some $(G,X)$-structures on a compact manifold $M$ ?
\end{que}


%

In higher rank, Guichard and Wienhard \cite{GWDod} constructed geometric sturctures associated to Anosov representation by constructing cocompact domains of disconinuity in the corresponding flag manifolds. 
%
%
%
%
%
Such domains have been constructed in more generality by Kapovich, Leeb and Porti \cite{KLP-DoD}. Since $\Theta$-positive representations are $\Theta$-Anosov \cite{positivityLabourie}, the construction of Kapovich, Leeb and Porti provides geometric structures associated to any $\Theta$-positive representation.
 
\medskip
 
The topology on the quotient of these domains, which are the manifolds on which we put a geometric structure, is not completely clear. Alessandrini, Li and the author \cite{ADL} studied the topology of this manifold for Hitchin representations in $\PSL(2n,\R)$. In this paper we study more generally the domains of discontinuity from \cite{KLP-DoD} that coincide with the domains $\Omega_\rho^\type$.

\medskip

Suppose that $G$ is center-free. Let $\mathfrak{h}$ be a $\mathfrak{sl}_2$ Lie subalgebra in the Lie algebra $\mathfrak{g}$ of $G$, \ie a Lie subalgebra isomorphic to $\mathfrak{sl}_2(\R)$. One can associate to $\mathfrak{h}$ a representation $\iota_\mathfrak{h}:\SL(2,\R)\to G$ and a $\iota_\mathfrak{h}$-equivariant and totally geodesic embedding $u_\mathfrak{h}:\mathbb{H}^2\to \X$. Choose $\type\in \Sph\mathfrak{a}^+$ such that $u_\mathfrak{h}$ is $\type$-regular so that $u_\mathfrak{h}$ is $\type$-nearly geodesic.
 A representation preserving and acting properly and cocompactly on $u_\mathfrak{h}\left(\mathbb{H}^2\right)$ will be called \emph{$\mathfrak{h}$-generalized Fuchsian}.

\begin{thm}[Theorem \ref{thm:HighTeich}]
\label{thm:IntroHighTeich}
Let $\rho:\Gamma_g\to G$ be a $\mathfrak{h}$-generalized Fuchsian representation. Suppose that the associated domain of discontinuity $\Omega_\rho^\type$ is non-empty.
Let $\mathcal{C}$ be the connected component of the space of $\Theta$-Anosov representations that contains $\rho$, for some non-empty $\Theta$ depending on $\mathfrak{h}$. Every representations in $\mathcal{C}$ is the restricted holonomy of a $\left(G,\mathcal{F}_{\type}\right)$-structure \emph{on a fiber bundle} over $S_g$.
\end{thm}

 The holonomy of the structure along the fibers is be trivial, so one can consider the restricted holonomy even if the fiber has non-trivial fundamental group. The fiber bundle is induced as the reduction of a $\SO(2,\R)\times \mathcal{C}_K(\mathfrak{h})$-principal bundle over $S_g$ via the action of $\SO(2,\R)\times \mathcal{C}_K(\mathfrak{h})$ on a codimention two submanifold of the flag manifold $\mathcal{F}_\type$ associated to $\mathfrak{h}$. This codimention $2$ submanifold, which is the fiber of the fiber bundle, can be described as union of connected components of the $\type$-base of a pencil of tangent vectors in $\X$ associated to $\mathfrak{h}$.

\medskip

This theorem applies to connected components of $\Theta$-positive representations that contain a $\mathfrak{h}_\Theta$-generalized Fuchsian representation, where $\mathfrak{h}_\Theta$ is a $\Theta$-principal $\mathfrak{sl}_2$ Lie subalgebra. These components are known to form higher rank Teichmüller spaces of $\Theta$-positive representations \cite{positivityLabourie}. Suppose that $G$ admits a notion of $\Theta$-positivity and is not isomorphic to $\PSL(2,\R)$.

\begin{cor}[Corollary \ref{cor:ThetaStruct}]
\label{cor:CayleyIntro}
Every $\Theta$-positive representation $\rho:\Gamma_g\to G$ in a component containing a $\mathfrak{h}_\Theta$-generalized Fuchsian representation is the restricted holonomy of a $(G, \mathcal{F}_{\type_{\Theta'}})$-structure on a fiber bundle over $S_g$, for every Weyl orbit of simple roots $\Theta'\subset \Theta$.
\end{cor}

In this theorem  $\mathcal{F}_{\type_\Theta}$ is a flag manifold associated to a set of simple roots $\Theta(\type_\Theta)\subset \Delta$. In Figure \ref{fig:TableWeylOrbit} we give the the list of all Weyl orbit of simple roots, together with the corresponding sets of simple roots $\Theta(\type_\Theta)$. Corollary \ref{cor:CayleyIntro} describes exactly one geometric structure  for each notion of $\Theta$-positivity, except for Hitchin representations in non-simply laced split Lie groups for which it describes two geometric structures.

\medskip
Any maximal representation $\rho:\Gamma_g\to\Sp(2n,\R)$ for $n\geq 3$ is the restricted holonomy of a projective structure \emph{on a fiber bundle over $S_g$} with fiber a Stiefel manifold. In the case $n=2$, such structures are already well described by Collier, Tholozan, Toulisse in \cite{Collier_2019}, and it is not clear if our method can apply to the \emph{Gothen components}. In Appendix \ref{appendix:CCT} we compare their fibration of the corresponding domain of discontinuity with the one that we construct for nearly Fuchsian representations.

\medskip

Every Hitchin representation $\rho:\Gamma_g\to \PSL(n,\R)$ is the restricted holonomy of a structure modeled on $\mathcal{F}_{1,n-1}$ \emph{on a fiber bundle over $S_g$}, where $\mathcal{F}_{1,n-1}$ is the space of partial flags in $\R^n$ consisting of the form $(\ell,H)$ where $\ell\subset H\subset \R^n$, $\dim(\ell)=1$ and $\dim(H)=n-1$.

\medskip

Similarly every $\Theta$-positive representation $\rho:\Gamma_g\to\SO(p,q)$ for $3\leq p$, $p+1<q$ is the restricted holonomy of a structure modeled on $\mathcal{F}_2$ \emph{on a fiber bundle over $S_g$} where $\mathcal{F}_2$ is the space of isotropic planes in $\R^{p,q}$. If $q=p+1$, there are some exceptional components that are conjectured to consist only of Zariski-dense representations, so it is not clear if our method can be applied to these.

\subsection*{Aknowledgments}

I would like to thank Max Riestenberg and Gabriele Viaggi for very precious discussions on the geometry of symmetric symmetric spaces.  I would like to thank Thi Dang and Johannes Horn for raising interesting questions about nearly geodesic immersions. Finally I would like to thank Beatrice Pozzetti for her comments and support throughout this work and the redaction of the paper. The author was funded through the DFG Emmy Noether project 427903332 of B. Pozzetti. 

\section{Symmetric spaces of non-compact type.}
\label{sec:symspaces}

In this section we recall the general theory of symmetric spaces of non compact type and fix some notations. References for the results mentioned can be found in \cite{Helgason} and \cite{Eberlein}. We then illustrate some of these notions for some families of Lie groups. Finally we introduce the notion of Weyl orbit of simple roots.

\subsection{Symmetric space associated to a semi-simple Lie group.}

Let $G$ be a connected, semi-simple Lie group with finite center and no compact factors, \ie of non-compact type.

\medskip

Let $\mathfrak{g}$ be the Lie algebra of $G$, and let $B$ be the \emph{Killing fom} on $\mathfrak{g}$. Since $\mathfrak{g}$ is semi-simple it admits a \emph{Cartan involution} \ie an involutive automorphism $\theta:\mathfrak{g}\to \mathfrak{g}$ such that $(\mathrm{v},\mathrm{w})\mapsto-B(\mathrm{v} ,\theta (\mathrm{w}))$ is a scalar product on $\mathfrak{g}$. Any two Cartan involutions are conjugated by $\Ad_g$ for some $g\in G$.

\medskip

Let $\X$ be the space of Cartan involutions of $\mathfrak{g}$. For any $x\in X$ we will write the corresponding Cartan involution $\theta_x:\mathfrak{g}\to \mathfrak{g}$. This involution determines a $B$-orthogonal decomposition $\mathfrak{g}=\mathfrak{t}_x\oplus \mathfrak{p}_x$, where $\mathfrak{t}_x$ is the $+1$ eigenspace of $\theta$, and by $\mathfrak{p}_x$ the $-1$ eigenspace. 

\medskip

For $x\in X$, define $K_x$ to be the group of elements $k\in G$ such that $\Ad_k$ commutes with $\theta_x$. This subgroup is a maximal compact subgroup of $G$. Given any $x\in \X$, one can identify $\X$ with the homogeneous space $G/K_x$. The Lie subalgebra $\mathfrak{t}_x$ is the Lie algebra of the compact $K_x$, and thus the space $\mathfrak{p}_x$ is naturally identified with $T_x\X$.

\medskip

Let $\langle\cdot,\cdot\rangle_x$ be the scalar product defined for $\mathrm{v},\mathrm{w}\in \mathfrak{g}$ as :
\begin{equation}
\label{eq:metric}
\langle\mathrm{v},\mathrm{w} \rangle_x=B\left(\mathrm{v},\theta_x(\mathrm{w})\right).
\end{equation}

This scalar product restricted to $\mathfrak{p}_x\simeq T_x\X$ defines a Riemannian metric $g_\X$ on $\X$. We will denote by $d_\X$ the induced Riemannian distance on $\X$. With this metric the space $\X$ is a symmetric space in the sense that for all $x\in \X$ there is an isometry $\sigma_x$ of $X$ such that $\mathrm{d}_x\sigma=-\Id$.

\medskip

The symmetric space $\X$ is of \emph{non-compact type}. It is simply connected and has non-positive sectional curvature. In particular it is a \emph{Hadamard manifold}.

\begin{rem}
We only consider symmetric spaces $\X$ associated to semi-simple Lie groups $G$, having their Riemannian metric defined via the Killing form. 
\end{rem}

\subsection{Reduced root systems.}
\label{subsec:RootSystems}
Fix a base point $o\in \X$. Let $\mathfrak{a}$ be a choice of a \emph{maximal abelian subalgebra} of $\mathfrak{p}_o$. These maximal abelian subalgebras are all conjugated by elements of $K_x$. The dimension $\rank(\X)$ of $\mathfrak{a}$ will be called the \emph{rank} of $\X$.

\begin{rem}
In general $\rank(\X)\leq \rank(G)$, where $\rank(G)$ is the dimension of any Cartan subalgebra in $\mathfrak{g}$.
\end{rem}

Let $\alpha\in \mathfrak{a}^*$ be a linear form. Let $\mathfrak{g}_\alpha$ be the set of elements $\mathrm{v}\in\mathfrak{g}$ such that for all $\type\in \mathfrak{a}$:
$$\ad_\type(\mathrm{v})=\alpha(\type) \mathrm{v}.$$

The \emph{reduced root system} $\Sigma$ is the set of linear forms $\alpha\in \mathfrak{a}^*$ such that $\mathfrak{g}_\alpha\neq \lbrace0\rbrace$.  An element $\type\in\mathfrak{a}$ is \emph{regular} if for all $\alpha \in \Sigma\setminus \lbrace0\rbrace$, $\alpha(\type)\neq 0$. 

\medskip

Let us choose a regular element $\typeS\in \mathfrak{a}$. Let $\Sigma^+$ be the associated set \emph{positive roots}, \ie the set of $\alpha \in \Sigma$ such that $\alpha(\typeS)>0$. There exists a unique set $\Delta$ of linearly independent roots in $\Sigma^+$ such that any root if $\Sigma^+$ can be written as a linear combination of roots in $\Delta$. The roots in $\Delta$ are called \emph{simple roots}. 

\medskip

Let the \emph{Weyl group} $W$ be quotient of the subgroup of elements in $K_x$ whose adjoint action stabilizes $\mathfrak{a}$ by the subgroup of elements who fix $\mathfrak{a}$ point-wise.

For any root $\alpha \in \Sigma\setminus\lbrace0\rbrace$ there is an element $\sigma_\alpha\in W$ such that its action on $\mathfrak{a}$ is the orthogonal symmetry with respect to $\Ker(\alpha)$ in $\mathfrak{a}$. The Weyl group acts linearly on $\mathfrak{a}$, and it is generated by the elements $(\sigma_\alpha)_{\alpha\in \Delta}$.
The \emph{model Weyl chamber} $\mathfrak{a}^+$ is the cone $\lbrace\type\in \mathfrak{a}|\forall \alpha \in \Delta,\, \alpha(\type)\geq 0\rbrace$. For any $\typeS\in \mathfrak{a}$ there is a unique $\type\in \mathfrak{a}^+$ such that for some $w\in W$, $w\cdot \typeS=\type$.  An element $\type\in \mathfrak{a}^+$ is $\Theta$-regular for $\Theta\subset \Delta$ if for all $\alpha \in \Theta$, $\alpha(\type)\neq 0$.

\medskip

We denote by $\Sph\mathfrak{a}$ and $\Sph\mathfrak{a}^+$ reprectively the unit sphere in $\mathfrak{a}$ and the unit sphere intersected with the model Weyl chamber $\mathfrak{a}^+$, for the metric \eqref{eq:metric}.

\medskip 

Let $w_0\in W$ be the only element such that $w_0\cdot \mathfrak{a}^+=-\mathfrak{a}^+$. This element is called \emph{the longest element} of the Weyl group. Let $\iota:\mathfrak{a}^+\to\mathfrak{a}^+$ be the involution such that for $\type\in \Sph\mathfrak{a}^+$, $\iota(\type)=-w_0\cdot \type$. An element $\type\in \Sph\mathfrak{a}^+$ is called \emph{symmetric} if $\iota(\type)=\type$.

\medskip

Given a set of simple roots $\Theta\subset \Delta$, we say that the \emph{model $\Theta$-facet} is the set of elements $\type\in \Sph\mathfrak{a}^+$ such that for all $\alpha\in \Delta\setminus \Theta$, $\alpha(\type)=0$. The \emph{open} model $\Theta$-facet is the set of elements $\type\in \Sph\mathfrak{a}^+$ such that for all $\alpha\in \Delta\setminus \Theta$, $\alpha(\type)=0$, and for all $\alpha \in \Theta$, $\alpha(\type)>0$. For an element $\type\in \Sph\mathfrak{a}^+$ we will write $\Theta(\type)$ the unique set of simple roots such that $\alpha$ lies in the open model $\Theta(\type)$-facet.

\subsection{Maximal Flats, visual boundary and parabolic subgroups.}

A \emph{flat} in $\X$ is a complete totally geodesic subspace of $\X$ on which the sectional curvature completely vanishes. A flat $F$ is maximal if $\dim(F)=\rank(\X)$. Flats passing through a point $x\in \X$ are in one to one correspondence with abelian subalgebras of $\mathfrak{p}_x$, and maximal flats correspond to maximal subalgebras. As a consequence the action of $G$ on the space of maximal flats is transitive. The maximal flat corresponding to $\mathfrak{a}$ will be called the \emph{model flat}.
Moreover for any $x\in X$ and $\mathrm{v}\in T_x\X$, there is a maximal flat $F$ such that $x\in F$ and $\mathrm{v}\in T_xF$.

\medskip

We say that two geodesic rays parametrized with unit length $\eta_1,\eta_2:\R_{\geq 0}\to \X$ are \emph{asymptotic} if there exist a positive constant $C$ such that for all $t>0$, $d_\X(\gamma_1(t),\gamma_2(t))\leq C$. This defines an equivalence relation on the space of rays.

\medskip

The \emph{visual boundary} $\partialvis\X$ of the symmetric space $\X$ is the space of classes of asymptotic geodesic rays parametrized with unit speed. The group $G$ acts by isometries on $\X$, hence it acts on $\partialvis \X$.

A unit vector $\mathrm{v}\in T_x\X$ at a point $x\in \X$ \emph{points towards} $a\in \partialvis \X$ if the geodesic ray $\gamma$ such that $\gamma(0)=x$, $\gamma'(0)=\mathrm{v}$ is in the class corresponding to $a$. Since $\X$ is a Hadamard manifold, for any $x\in \X$ and $a\in \partialvis \X$ there is a unique unit vector that points towards $a$. We will denote this vector by $\mathrm{v}_{a,x}$ throughout the paper.
There exist a unique topology on $\partialvis \X$ such that for any $x\in X$ the map $\phi_x:a\mapsto \mathrm{v}_{a,x}$ is an homeomorphism between $\partialvis \X$ and $T^1_x\X$.

\medskip

The visual boundary of the model flat can be identified with $\Sph\mathfrak{a}$, and is included in the visual boundary of $\X$. The $G$-orbit of a point $\type\in \Sph\mathfrak{a}^+$ will be denoted by $\mathcal{F}_\type$.
 A \emph{$\Theta$-facet} is a subset of $\partialvis \X$ that is the image of the model $\Theta$-facet by the action of an element of $G$. We define similarly the notion of \emph{open $\Theta$-facet}. The stabilizer of the open model $\Theta$-facet will be denoted by $P_\Theta$, and we will denote by $\mathcal{F}_\Theta$ the associated \emph{flag manifolds}, \ie the quotient $G/P_\Theta$. Since $P_\Theta$ is also the stabilizer of any point in the open $\Theta$-facet, there exist a natural $G$-equivariant diffeomorphism between the $G$-homogeneous spaces $\mathcal{F}_\Theta$ and $\mathcal{F}_\type$ for any $\type$ in the open model $\Theta$-facet.
 
\begin{exmp}
Let us consider the case when $G=\PSL(n,\R)$ to illustrate these notions. A parabolic subgroup is this case is the stabilizer of a partial flag $f$. Any point in $\partialvis \X$ belongs to a unique open facet, which corresponds to a partial flag. The type of the partial flag, \ie the dimensions of the subspaces that form the flag, determine a set of roots $\Theta_f$. The points in $\partialvis \X$ are in $1$ to $1$ correspondence with partial flags decorated with a point in the open $\Theta_f$-model facet. This decoration can be interpreted as a collection of weights associated to the subspaces of the partial flag.
\end{exmp}

Given any two points $a,a'\in \partialvis \X$, one can define their \emph{Tits angle} $\Tangle(a,a')$ as the minimum of $\angle(\mathrm{v}_{a,x},\mathrm{v}_{a',x})$ for $x\in \X$. This minimum is obtained when $x\in \X$ lies in a common flat with $a$ and $a'$.

\subsection{Cartan and Iwasawa decomposition.}
\label{subsec:CartanDec}
The \emph{Cartan projection} $\Cartan:T \X\to \mathfrak{a}^+$ is the function that maps any vector $\mathrm{w}\in T_x\X$ to the unique element $\Cartan(\mathrm{w})\in\mathfrak{a}^+$ of the model Weyl Chamber such that for some $g\in G$, $g\cdot \mathrm{w}=\Cartan(\mathrm{w})$. 

\medskip

The \emph{generalised distance} $\da(x,y)$ based at $x\in \X$ of a point $y\in \X$ is the Cartan projection $\Cartan(\mathrm{w})$ of the unique vector $\mathrm{w}\in T_x\X\simeq\mathfrak{p}_x$ such that $\exp(\mathrm{w})\cdot x=y$. This generalized distance is $1$-Lipshitz in the following sense:

\begin{lem}
\label{lem:GeneralizedDistanceLip}

Let $x,y,z\in\X$:
$$|\da(x,z)-\da(x,y)|\leq d_\X(y,z).$$

\end{lem}

A proof of this lemma can be found for instance in {\cite{Riestenberg}, Corollary 3.8}. Here $|\cdot|$ means the norm induced by the metric \eqref{eq:metric}.

\medskip

We say that a vector $\mathrm{v}\in T\X$ is \emph{$\Theta$-regular} for a set of simple root $\Theta$ if it's Cartan projection is $\Theta$-regular, \ie it avoids the walls of the Weyl chamber associated with elements of $\Theta$. We will later introduce a similar notion of a $\type$-regular vector in Definition \ref{defn:pRegularVector}.

\medskip

Let $T_{a,x}:P_a\to G$ be the map that associates to $g\in G_a$ the limit : $$\lim_{t\to +\infty} \exp(-t\mathrm{v}_a)g\exp(t\mathrm{v}_a)$$
 
  It is a well defined continuous morphism.
 Let $N_{a,x}$ be the kernel of $T_{a,x}$, and $\mathfrak{n}_{a,x}$ its Lie algebra. The generalised Iwasawa decomposition is useful to compute Busemann functions.

\begin{thm}[Generalized Iwasawa decomposition]
Let $x\in \X$ and $a\in \partial_{\vis}\X$. then the following map:
\begin{align*}
 N_{a,x}\times \exp(\mathfrak{a}_{a,x})\times K_x\longrightarrow G \\
 (n,\exp(\mathrm{v}),k)\longmapsto n\exp(\mathrm{v})k
\end{align*}
 
is a diffeomorphism. In particular for every $x\in \X$ and $a\in \partial_{\vis}\X$ there is a splitting :
 $$\mathfrak{g}=\mathfrak{n}_{a,x}\oplus \mathfrak{a}_{a,x}\oplus \mathfrak{k}_{x}.$$

The sum $\mathfrak{n}_{a,x}\oplus \mathfrak{a}_{a,x}$ is orthogonal with respect to $\langle\cdot,\cdot\rangle_x$.
\end{thm}

In this Theorem, $\mathfrak{a}_{a,x}\subset \mathfrak{p}_x$ is the centralizer of $\mathrm{v}_{a,x}$.

\subsection{Examples.}
\label{subsec:ExamplesRoots}
In this subsection we consider the case when the semi-simple Lie group $G$ is equal to $\PSL(n,\R)$, $\PSL(n,\C)$, $\PSp(2n,\R)$ or $\PSO(p,q)$. The notations introduced here will be used in the examples throughout the paper.

\medskip

\paragraph*{Let $G=\PSL(n,\R)$} and let us fix a volume form on $\R^n$. Let $\mathcal{S}_n$ for $n\geq 2$ be the space of all scalar products on $\R^n$ having volume one.  The group $\PSL(n,\R)$ acts  transitively on $\mathcal{S}_n$ by changing the basis, i.e for $g\in \PSL(n,\R)$, $q\in X$ and $v,w\in \R^n$: $g\cdot q(v,w)=q(g^{-1}(v),g^{-1}(w))$. For any $q\in \mathcal{S}_n$ the space $\mathcal{S}_n$ can be identified with the quotient $\PSL(n,\mathbb{R})/\PSO(q)\cong \PSL(n,\mathbb{R})/\PSO(n,\mathbb{R})$.

\smallskip

Let $\theta_q$ at a point $q\in X$  be the involutive automorphism of $\mathfrak{sl}(n,\R)$ defined by $u\mapsto -\transpose{u}$ where $\transpose{u}$ is the transpose of $u$ with respect to the scalar product $q$. This is a Cartan involution. The space $\mathcal{S}_n$ is the symmetric space of non-compact type associated to $G=\PSL(n,\R)$. 

The space $\mathfrak{p}_q$ is the space of symmetric endomorphsisms with respect to $q$, and $\mathfrak{t}_q$ is the space of antisymmetric endomorphsisms with respect to $q$. The scalar product $\langle \cdot,\cdot \rangle_q$ at a point $q\in \mathcal{S}_n$ is equal to $\langle u,v \rangle_q=2n\tr(\transpose{u}v)$ for $u,v\in \mathfrak{sl}(n,\R)$.

\smallskip

We choose the standard scalar product $q\in \mathcal{S}_n$ on $\R^n$ to be our base point of $\mathcal{S}_n$. A maximal abelian subalgebra $\mathfrak{a}\subset \mathfrak{p}_q\subset \mathfrak{sl}(n,\R)$ is equal to the algrbra of diagonal matrices:
$$\mathfrak{a}=\left\lbrace\Diag(\sigma_1,\cdots,\sigma_n)|\sigma_1,\cdots,\sigma_n\in \R^n, \sum_{i=1}^n\sigma_i=0\right\rbrace.$$

A choice of simple root is $\Delta=\lbrace\alpha_1,\cdots,\alpha_{n-1}\rbrace$ where for any $1\leq i\leq n-1$ and any $\type=\Diag(\sigma_1,\cdots,\sigma_n)\in \mathfrak{a}$, $\alpha_i(\type)=\sigma_i-\sigma_{i+1}$. 

\medskip

The Weyl chamber associated to this choice is :
$$\mathfrak{a}^+=\left\lbrace\Diag(\sigma_1,\cdots,\sigma_n)|\sigma_1\geq \cdots\geq \sigma_n\in \R^n, \sum_{i=1}^n\sigma_i=0\right\rbrace.$$

The Weyl group $W$ is isomorphic to $\mathfrak{S}_n$. It acts on $\mathfrak{a}$ by permuting the entries.

\medskip

\begin{figure}[h]
\begin{center}
\begin{tikzpicture}[scale=2]
\draw[red, thin, dashed] (0,0) circle(1);
\fill[blue, opacity=.3] (0,0) -- (0:1.6) -- (60:1.6);
\draw (0,0) -- (0:1.6);
\draw (0,0) -- (60:1.6);
\draw (0,0) -- (120:1.6);
\draw (0,0) -- (180:1.6);
\draw (0,0) -- (240:1.6);
\draw (0,0) -- (300:1.6);
\draw[->,very thick, red] (0,0) -- (30:1);
\node[above right,red] () at (28:.95) {$\type_\Delta$};
\draw[->,very thick, red] (0,0) -- (60:1);
\node[right,red] () at (60:1.06) {$\type_1$};
\node[above right,blue] () at (30:1.35) {$\mathfrak{a}^+$};
\node[below,red] () at (270:1) {$\Sph\mathfrak{a}$};
\end{tikzpicture}
\end{center}

\caption{The model restricted Cartan algebra $\mathfrak{a}$ and its Weyl chamber $\mathfrak{a}^+$ for $\PSL(3,\R)$.}
\label{fig:WeylSL3}
\end{figure}
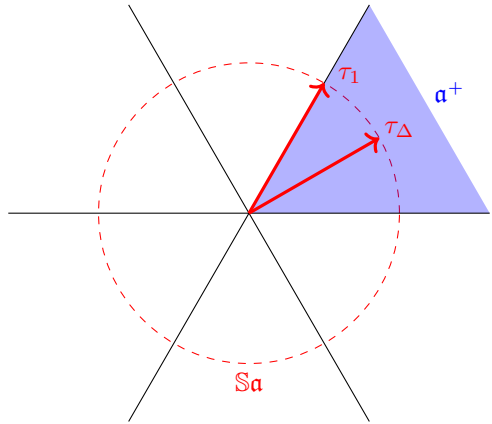

\begin{figure}[h]
\begin{center}
%
%
%
%
%
%
\begin{tikzpicture}[scale=3]

\fill[blue, opacity=.3] (-1,0) -- (1,0) -- (0,1) -- cycle ;
\draw[thick] (-1,0) -- (1,0) -- (0,1) -- cycle ;
\fill[ red] (0,0) circle (.02);
\node[above,red] () at (0,0) {$\type_{\Delta}$};

\node[left] () at (-.5,0.5) {$Ker(\alpha_1)$};

\node[right] () at (.5,0.5) {$Ker(\alpha_3)$};

\node[below] () at (0,0) {$Ker(\alpha_2)$};
\node[blue] () at (0,0.5) {$\Sph\mathfrak{a}^+$};

\end{tikzpicture}
\end{center}

\caption{The projectivization of the Weyl chamber $\Sph\mathfrak{a}^+$ for $G=\PSL(4,\R)$ in an affine chart.}
\label{fig:WeylSL4}
\end{figure}
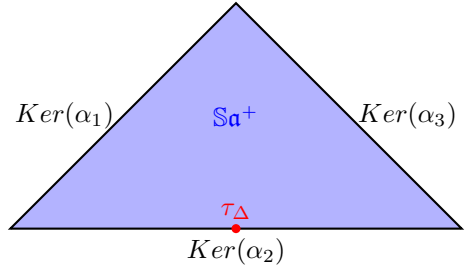

\paragraph*{Let $G=\PSL(n,\C)$}The space $\mathcal{H}_n$ of all definite positive Hermitian bilinear forms of $\C^n$ having volume one can identified with $\PSL(n,\mathbb{C})/\PSU(n,\mathbb{C})$. It can be given in a similar way a Riemannian metric that makes it a symmetric space of non-compact type associated to $G=\PSL(n,\mathbb{C})$.

The subalgebra $\mathfrak{a}\subset \mathfrak{sl}(n,\R)\subset \mathfrak{sl}(n,\C)$ defined previously is still a maximal abelian subalgebra of $\mathfrak{p}_q$. One has $\rank(\mathcal{S}_n)=\rank(\mathcal{H}_n)=\rank(\PSL(n,\R))=n-1$, but $\rank(\PSL(n,\C))=2n-2$.
\smallskip

\bigskip

\paragraph*{Let $G=\PSp(2n,\R)$} Let $\omega$ be a symplectic form on $\R^{2n}$. Let $\mathcal{X}_n$ be the space of endomorphisms $J$ on $\R^{2n}$ such that $J^2=-\Id$ and $(\mathrm{v},\mathrm{w})\mapsto \omega(\mathrm{v} ,J (\mathrm{w}))$ is a scalar product on $\R^{2n}$. The semi-simple Lie group $\PSp(2n,\R)$ acts on $\mathcal{X}_n$ by conjugation. The space $\mathcal{X}_n$ can be identified with $\PSp(2n,\R)/\PSU(n,\R)$. This is one of the models for the \emph{Siegel space}, see for instance \cite{BPSiegel}. 

\smallskip

For $J\in \mathcal{X}_n$, let us write  $\theta_J=\Ad_J$. This is Cartan involution of $\mathfrak{sp}(2n,\R)$. The Siegel space is the symmetric space of non-compact type associated with $G=\PSp(2n,\R)$.

\smallskip

 Let $\omega$ and $J\in \mathcal{X}_n$ be such that for $\mathrm{x}=(x_1,\cdots,x_{2n})$ and $\mathrm{y}=(y_1,\cdots,y_{2n})$ :
\[\omega(\mathrm{x},\mathrm{y})=\sum^{n}_{i=1} x_iy_{2n-i}-\sum^{2n}_{i=n+1} x_iy_{2n-i},\]
\[ J(x)=(-x_{2n},\cdots,-x_{n+1}, x_n, \cdots ,x_1).\]

A maximal abelian subalgebra $\mathfrak{a}\subset \mathfrak{p}_q\subset \mathfrak{sp}(2n,\R)$ is:
$$\mathfrak{a}=\left\lbrace\Diag(\sigma_1,\cdots,\sigma_n, -\sigma_n, \cdots, -\sigma_1)|\sigma_1,\cdots,\sigma_n\in \R^n\right\rbrace.$$

A choice of simple roots is $\Delta=\lbrace\alpha_1,\cdots,\alpha_{n}\rbrace$ where $\alpha_i(\type)=\sigma_i-\sigma_{i+1}$ for $1\leq i\leq n-1$ and $\alpha_n(\type)=2\sigma_n$.

The Weyl chamber associated to this choice is :
$$\mathfrak{a}^+=\left\lbrace\Diag(\sigma_1,\cdots,\sigma_n, -\sigma_n, \cdots, -\sigma_1)|\sigma_1\geq \cdots\geq \sigma_n\geq 0\in \R^n\right\rbrace.$$

The Weyl group $W$ is isomorphic to the subgroup of elements in $\mathfrak{S}_{2n}$ that commutes with the involution $\iota:i\mapsto 2n+1-i$. It acts on $\mathfrak{a}$ by permuting the entries.

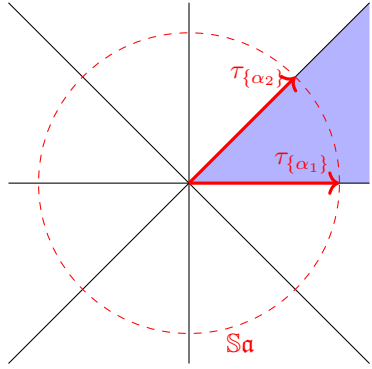
\begin{figure}[h]
\begin{center}
\begin{tikzpicture}[scale=2]
\fill[blue, opacity=.3] (0,0) -- (0:1.2) -- (1.2,1.2);

\draw (0,0) -- (0:1.2);
\draw (0,0) -- (90:1.2);
\draw (0,0) -- (180:1.2);
\draw (0,0) -- (270:1.2);
\draw (0,0) -- (1.2,1.2);
\draw (0,0) -- (-1.2,1.2);
\draw (0,0) -- (1.2,-1.2);
\draw (0,0) -- (-1.2,-1.2);

\draw[red, thin, dashed] (0,0) circle(1);
\draw[->,very thick, red] (0,0) -- (0:1);
\node[above left,red] () at (0:1) {$\type_{\lbrace\alpha_1\rbrace}$};
\draw[->,very thick, red] (0,0) -- (0.707,0.707);
\node[left,red] () at (.707,0.707) {$\type_{\lbrace\alpha_2\rbrace}$};
\node[below,red] () at (290:1) {$\Sph\mathfrak{a}$};
\end{tikzpicture}
\end{center}

\caption{The model restricted Cartan algebra $\mathfrak{a}$ and its Weyl chamber $\mathfrak{a}^+$ for $\PSp(4,\R)$.}
\label{fig:WeylSp4}
\end{figure}

\paragraph*{Let $G=\SO(p,q)$} with $p< q$. Let $\R^{p,q}$ be the vector space $\R^{p+q}$ equipped with a symmetric bilinear form $\langle\cdot,\cdot \rangle$ of signature $(p,q)$ defined in the standard basis by:
$$\langle\mathrm{x},\mathrm{y}\rangle=\sum_{i=i}^p \left(x_iy_{p+q-i}+x_{p+q-i}y_i\right)-\sum_{i=1}^{p-q}x_{p+i}y_{p+i}.$$

 A model for the associated symmetric space $\X$ is the space of spacelike subspaces $U\subset \R^{p,q}$, \ie subspaces on which $\langle\cdot,\cdot \rangle$ is definite positive.

A maximal abelian subalgebra $\mathfrak{a}\subset \mathfrak{p}_U\subset \mathfrak{so}(p,q)$ is the algebra :
$$\mathfrak{a}=\left\lbrace\Diag(\sigma_1,\cdots,\sigma_p,0,\cdots,0, -\sigma_p, \cdots, -\sigma_1)|\sigma_1,\cdots,\sigma_p\right\rbrace.$$

A choice of simple roots is $\Delta=\lbrace\alpha_1,\cdots,\alpha_{p}\rbrace$ where $\alpha_i(\type)=\sigma_i-\sigma_{i+1}$ for $1\leq i\leq p-1$, and $\alpha_p(\type)=\sigma_p$.

The Weyl chamber associated to this choice is :
$$\mathfrak{a}^+=\left\lbrace\Diag(\sigma_1,\cdots,\sigma_p,0,\cdots,0, -\sigma_p, \cdots, -\sigma_1)|\sigma_1\geq \cdots\geq \sigma_p\geq 0\right\rbrace.$$

The Weyl group $W$ is isomorphic to the subgroup of elements in $\mathfrak{S}_{2p}$ that commutes with the involution $\iota:i\mapsto 2p+1-i$. It acts on $\mathfrak{a}$ by permuting the first and last $p$ entries.

\subsection{Weyl orbits of simple roots.}
\label{sec:RootSpecial}

In this subsection we introduce \emph{Weyl orbits of simple roots}, which are special sets of simple roots. To a Weyl orbits of simple roots $\Theta$ one can associate a unit vector in the Weyl chamber $\type_\Theta\in \Sph\mathfrak{a}^+$ which is colinear to a coroot. 

\medskip

We consider the restricted root system $\Sigma$ associated with the semi-simple Lie group $G$, with a choice of a set of positive roots $\Sigma^+$ and of simple roots $\Delta$. 
Two simple roots $\alpha$ and $\beta$ are \emph{conjugates} if there is an element $w$ in the Weyl group $W$ such that $\alpha=w\cdot \beta=\beta\circ w^{-1}$.
 
\begin{defn}
A set of simple roots $\Theta\subset \Delta$ is called a \emph{Weyl orbit of simple roots} if it is an equivalence class for the conjugation relation on the set of simple roots $\Delta$.
\end{defn}

\begin{prop}
Let $\Theta$ be a Weyl orbit of simple roots. There exists a unique unit vector $\type_\Theta\in \Sph\mathfrak{a}^+$ such that for any $\alpha \in \Theta$ there is some $w\in W$ such that $w\cdot \type_\Theta$ is orthogonal to $\ker(\alpha)$. 
The vector $\type_\Theta\in \Sph\mathfrak{a}^+$ will be called the \emph{normalized coroot} associated to $\Theta$.
\end{prop}

The normalized coroot associated to $\Theta$ is colinear to a coroot which is itself conjugate via the Weyl group to the corroot associated to any $\alpha\in \Theta$.

\medskip

\begin{proof}
Let $\alpha\in \Theta$. Let $\type_0\in \Sph\mathfrak{a}$ be a unit vector orthogonal to $\ker(\alpha)$. Since every orbit for the action of the Weyl group on $\Sph\mathfrak{a}^+$ intersects exactly once the model Weyl chamber, there exist a unique vector $\type_\Theta\in \Sph\mathfrak{a}^+$ such that $\type_\Theta=w\cdot \type_0$ for some $w\in W$. 

This definition does not depend of the choice of $\alpha\in \Theta$, because if $\beta\in \Theta$ then for some $w_0\in W$, $w_0\cdot \beta=\alpha$ and hence any vector $\type'_0$ orthogonal  to $\ker(\beta)$ can be written $\type'_0=w_0\cdot \type_0$ or $\type'_0=(w_0\sigma_\beta)\cdot \type_0$, and hence $W\cdot \type_0=W\cdot \type'_0$. Therefore $W\cdot \type_0\cap \Sph\mathfrak{a}^+=W\cdot \type'_0\cap \Sph\mathfrak{a}^+=\lbrace \type_\Theta\rbrace$.
\end{proof} 

Note that in particular $\type_\Theta$ is symmetric.

\begin{rem}
If $\Theta$ is a Weyl orbit of simple roots, $\mathcal{F}_{\type_\Theta}$ is not in general the same flag manifold as $\mathcal{F}_\Theta=G/P_\Theta$.
\end{rem}

The \emph{Dynkin diagram} associated to the restricted root system $\Sigma$ is the graph with vertex set $\Delta$ such that for all $\alpha,\beta\in \Delta$ distinct roots there is a link between $\alpha$ and $\beta$ of multiplicity depending of the order $k$ of $\sigma_\alpha \sigma_\beta$ or $\sigma_\beta \sigma_\alpha$, where $\sigma_\alpha,\sigma_\beta\in W$ are the symmetries associated with the roots $\alpha, \beta$. If $k=2$ we consider that there is no link, there is a simple link if $k=3$, a double link if $k=4$ and a triple link if $k=6$. These are the only cases that occur for spherical Dynkin diagrams. If two roots have different norms, we orient the edge towards the root with largest norm. 

\begin{prop}
Consider the Dynkin diagram associated with the reduced root system $\Sigma$, and remove all the double or triple edges. A set $\Theta\subset \Delta$ is a Weyl orbit of simple roots if and only if it is a connected component of this graph.
\end{prop}

\begin{proof}

Let $\alpha,\beta$ be two simple roots such that $\sigma_\alpha=w\sigma_\beta w^{-1}$ for some $w\in W$. Then the simple root $\alpha$ and $w\cdot \beta$ are proportional, hence $\alpha=w\cdot \beta$ or $\alpha=-w\cdot \beta$. In any case $\alpha$ and $\beta$ are conjugated. Reciprocally if $\alpha$ and $\beta$ are conjugates, then $\sigma_\alpha$ and $\sigma_\beta$ are conjugated.

\medskip

The system $(W, (\sigma_\alpha)_{\alpha\in \Delta})$ is a Coxeter system. Generators of a Coxeter system are conjugated to one another if and only if there is a path of single edges between the corresponding vertices in the Dynkin diagram (\cite{Coxeter} Proposition 2.1). Hence Weyl orbits of simple roots correspond exactly to connected components for the modified Dynkin diagram.\end{proof}

We can now describe the Weyl orbits of simple roots in $\Delta$ for the restricted system of roots associated to a simple group $G$. For this we use the classification of the Dynkin diagrams that occur as reduced root system for a symmetric space $\X$ associated to $G$.

\begin{cor}
\label{cor:RootsWeylOrbit}
If the restricted root system $\Sigma$ is of type $A_n,D_n$ for $n\geq 2$ or $E_6,E_7,E_8$, then the only Weyl orbit of simple roots in $\Delta$ is $\Delta$.

If the root system $\Sigma$ is of type $B_n,C_n$ for $n\geq 2$ or $F_4,G_2$, then $\Delta$ can be partitioned into its only two Weyl orbits of simple roots .
\end{cor}

\begin{exmp}
We keep notations from Section \ref{subsec:ExamplesRoots}. If $G=\PSL(n,\R)$, with previous notations $\Delta$ is the only Weyl orbit of simple roots  and:
$$\type_\Delta=\frac{1}{2\sqrt{n}}\Diag(1,0,\cdots,0,-1). $$

The flag manifold $\mathcal{F}_{p_\Delta}$ can be identified with :
$$\mathcal{F}_{1,n-1}=\lbrace(\ell,H)|\ell\subset H\subset \R^n, \dim(\ell)=1,\dim(H)=n-1\rbrace.$$

\medskip

If $G=\Sp(2n,\R)$,  with previous notations $\Theta=\lbrace\alpha_n\rbrace$ and $\Theta'=\lbrace\alpha_1,\alpha_2,\cdots,\alpha_{n-1}\rbrace$ are the two Weyl orbits of simple roots of $\Delta$. One has :
$$\type_\Theta=\frac{1}{2\sqrt{n}}\Diag(1,0,\cdots,0,-1). $$
$$\type_{\Theta'}=\frac{1}{2\sqrt{n}}\Diag(1,1,0,\cdots,0,-1,-1). $$

The flag manifold $\mathcal{F}_{\type_\Theta}$ can be identified with $\mathbb{RP}^{2n-1}$ and $\mathcal{F}_{\type_{\Theta'}}$ can be identified with the Grassmannian of planes $P$ in $\R^{2n}$ that are isotropic for $\omega$, \ie such that $\omega_{|P}=0$.
\end{exmp}
\begin{figure}
\setlength\extrarowheight{4pt}
\begin{tabular}{|c|c|c|c|}
    \hline
     $\Delta$  & $\Theta$ & $\type_\Theta$ & $\Theta(\type_\Theta)$\\ 
    \hline
     $A_n$ & 
     \dynkin A{**.**} & $\frac{1}{\sqrt{2}}\left(e_1-e_{n+1}\right)$ & \dynkin A{*o.o*}\\ 
    \hline
    $B_2$ & 
    \dynkin B{*o}& $e_1$ & \dynkin B{o*}\\
    \cline{2-4}&\dynkin B{o*}&  $\frac{1}{\sqrt{2}}\left(e_1+e_{2}\right)$ & \dynkin B{*o}\\
     \hline
    $B_n$ & 
    \dynkin B{oo.oo*}& $e_1$ & \dynkin B{*o.ooo}\\
    \cline{2-4}& \dynkin B{**.**o}&  $\frac{1}{\sqrt{2}}\left(e_1+e_{2}\right)$ & \dynkin B{o*.ooo}\\
     \hline
    $C_n$ & 
    \dynkin C{oo.oo*}& $e_1$ & \dynkin C{*o.ooo}\\
    \cline{2-4}& \dynkin C{**.**o}&  $\frac{1}{\sqrt{2}}\left(e_1+e_{2}\right)$ & \dynkin C{o*.ooo}\\
     \hline
    $D_n$ & 
    \dynkin D{**.****} & $\frac{1}{\sqrt{2}}\left(e_1-e_{n+1}\right)$ & \dynkin D{*o.oo**} \\
    \hline
    $E_6$ & 
    \dynkin E6 & $\sqrt{2}e$ & \dynkin E{o*oooo} \\
    \hline
    $E_7$ & 
    \dynkin E7 & $\frac{1}{\sqrt{2}}\left(e_8-e_{7}\right)$ & \dynkin E{*oooooo} \\
    \hline
    $E_8$ & 
    \dynkin E8 & $\frac{1}{\sqrt{2}}\left(e_1-e_{9}\right)$ & \dynkin E{ooooooo*} \\
     \hline
    $F_4$ & 
    \dynkin F{**oo}& $\frac{1}{\sqrt{2}}\left(e_1+e_{2}\right)$ & \dynkin F{*ooo}\\
    \cline{2-4}& \dynkin F{oo**}&  $e_1$ & \dynkin F{ooo*}\\
      \hline
    $G_2$ & 
    \dynkin G{*o}& $\frac{1}{\sqrt{2}}\left(e_1-e_{3}\right)$ & \dynkin G{*o}\\
    \cline{2-4}& \dynkin G{o*}&  $\frac{1}{\sqrt{6}}\left(2e_1-e_{2}-e_{3}\right)$ & \dynkin G{o*}\\
     \hline
\end{tabular}
\caption{Weyl orbits of simple roots, and their associated normalized coroots.}
\label{fig:TableWeylOrbit}
\end{figure}

In general, the Weyl orbits of simple roots for any root system are summarized in Figure \ref{fig:TableWeylOrbit}. The table also includes an illustration of the set of roots $\Theta(\type_\Theta)$ such that $\mathcal{F}_{\Theta(\type_\Theta)}\simeq  \mathcal{F}_{\type_\Theta}$. The sets of roots are illustrated in the diagram as the set of filled vertices. Using notations from \cite[Table 1, page 293]{LieGroupTable}, the basis $(e_i)$ is an orthonormal basis such that $e_i^\vee=\epsilon_i$ for $B_n$, $C_n$, $D_n$, $F_4$ and $e_i^\vee-\frac{1}{n+1}\sum_{k=1}^{n+1}e_k^\vee=\epsilon_i$ for $A_n,E_7,E_8$ and $G_2$. For $E_6$, we write $e=\epsilon^\vee$. 

The table can be checked as follows: for each Weyl orbit of simple root one can check that the vector $\type_\Theta$ is orthogonal to the kernel of a root conjugate to a root in $\Theta$, and lies in the model Weyl chamber. Then one an check that the simple roots that do not vanish on $\type_\Theta$ are the one in $\Theta(\type_\Theta)$, as depicted in Figure \ref{fig:TableWeylOrbit}.


\section{Representations of hyperbolic groups.}
\label{sec:Representations}
\subsection{Gromov hyperbolic groups.}

Let $\Gamma$ be a finitely generated group. Let $F$ be any finite generating system for $\Gamma_g$ that is symmetric, \ie such that  $s^{-1}\in F$ for all $s\in F$. We can define the norm of an element $\gamma\in \Gamma_g$ as:
$$|\gamma|_F=\min\left\lbrace n|n=s_1s_2\cdots s_n, \, s_i\in S\right\rbrace.$$

This norm defines the word distance on $\Gamma_g$ by taking $d_F(\gamma_1,\gamma_2)=|\gamma_1^{-1}\gamma_2|_F$ for $\gamma_1,\gamma_2\in \Gamma_g$.

\medskip

A map $f:Y\to X$ between two metric spaces $X,Y$ is called a \emph{quasi-isometric embedding} if there exist $C,D$ such that for all $x_1,x_2\in X$:
$$\frac{1}{C}d_Y(x_1,x_2)-D\leq d_X(f(x_1),f(x_2))\leq Cd_Y(x_1,x_2)+D .$$

By extension, we say that a representation $\rho$ is a \emph{quasi-isometric embedding} if some and hence any $\rho$-equivariant map $u_0:\Gamma\to \X$ is a quasi-isometry, where $\Gamma$ acts on itself by left multiplication. 

 This notion does not depend on the choice of $F$: indeed if $F'$ is an other finite generating system, the identity map $(\Gamma, d_F)\to (\Gamma, d_{F'})$ is a quasi-isometric embedding.
 
 \medskip
 
The group $\Gamma$ is called \emph{hyperbolic} if as a metric space it is hyperbolic in the sense of Gromov. We denote by $\partial\Gamma$ the Gromov boundary of an hyperbolic group $\Gamma$, that we equip with the usual topology \cite{Gromov}.

\medskip

Given a discrete representation, we will need to consider the limit cone of the Cartan projections of elements of the group.

\begin{defn}
\label{defn:CartanCone}
The \emph{limit cone} of a discrete representation $\rho:\Gamma\to G$ is the closed subset $$\mathcal{C}_\rho=\bigcap_{n\in \mathbb{N}}\overline{\lbrace\left[\da(o,\rho(\gamma)\cdot o)\right],|\gamma|_w\geq n\rbrace}=\bigcap_{n\in \mathbb{N}}\overline{\left\lbrace\frac{\da(o,\rho(\gamma)\cdot o)}{d_\X(o,\rho(\gamma)\cdot o)},|\gamma|_w\geq n\right\rbrace} \subset \Sph\mathfrak{a}^+.$$
\end{defn}

Recall that the generalized distance $\da$ was defined in Section \ref{sec:symspaces}. This definition does not depend on the choice of the base point $o\in \X$.

\begin{lem}
\label{lem:CartanConeConnected}
The limit cone $\mathcal{C}_\rho$ of a discrete representation $\rho:\Gamma\to G$ of a  non-elementary hyperbolic group is connected.\end{lem}

A hyperbolic group is non-elementary if it is neither finite or virtually cyclic.

\begin{rem}If $\rho$ is a Zariski dense representation, Benoist proved that $\mathcal{C}_\rho$ is convex \cite{BenoistLin1}. We give a proof of the connectedness of $\mathcal{C}_\rho$ to avoid considerations of Zariski density.
\end{rem}

\begin{proof}
We define the distance $d$ on $\Sph\mathfrak{a}$, as $d([\mathrm{x}],[\mathrm{y}])=|\frac{\mathrm{x}}{|\mathrm{x}|}-\frac{\mathrm{y}}{|\mathrm{y}|}|$ for $\mathrm{x},\mathrm{y}\in \mathfrak{a}$.

\medskip

Assume that there exists a partition $A\cup B$ of $\mathcal{C}_\rho$ into two open and closed sets. These sets are  compact and hence are at uniform distance $\epsilon>0$.
Let $F$ be a finite symmetric generating set for $\Gamma$. Let $M$ be the maximum for $\gamma\in F$ of $d_\X(o,\rho(\gamma)\cdot o)$. Let $E\subset\Gamma$ be the subset of elements $\gamma$ such that either $d_\X(o,\rho(\gamma)\cdot o)\leq \frac{12M}{\epsilon}$ or:
 $$d(\left[\da(o,\rho(\gamma)\cdot o)\right],\mathcal{C}_\rho)\geq\frac{\epsilon}{3}.$$
 
 Since $\rho$ is discrete and by the definition of $\mathcal{C}_\rho$, $E$ is finite. It is included in the ball of radius $R$ centered at the identity of $\Gamma$ for $|\cdot|_F$ and some $R>0$.

\medskip

We say that two elements $\gamma',\gamma''$ in $\Gamma\setminus E$ are \emph{$E$-connected} if one can construct a sequence $(\gamma_n)$ of elements of $\Gamma\setminus E$ such that for all $1\leq i< N$,$\gamma_{i+1}=a_i\gamma_{i}b_i$ for some $a_i,b_i\in F$, with $\gamma'=\gamma_0$ and $\gamma''=\gamma_N$. 

\medskip

Using only right translations, \ie $a_i=\Id$, one can see that any point in $ \partial\Gamma$ has a neighborhood whose intersection with $\Gamma$ is $E$-connected. Indeed given $D>0$ for any $\gamma\in \Gamma$ with $|\gamma|_F> R+D$ the set of elements $\gamma'$, such that the geodesic from $\Id$ to $\gamma'$ passes at distance at most $D$ of $\gamma$, is $E$-connected. Any point in the visual boundary admits neighborhood whose intersection with $\Gamma$ has this form.

\medskip

Using only left translations , \ie $b_i=\Id$, one can see that the intersection of a neighborhood of $\partial \Gamma$ and $\Gamma$ is $E$-connected. Indeed, the action of $\Gamma$ on $\partial \Gamma$ is minimal because $\Gamma$ is non-elementary. Therefore there exist $b_1\cdots b_n=\gamma\in \Gamma$ with $b_i\in F$ such that $\gamma\cdot x$ lies in the interior of $V$. Hence any $\gamma'\in V\cap \Gamma$ close enough to $\gamma\cdot x$ is $E$-connected to $U\cap \Gamma$. 
Therefore given $x,y\in \partial \Gamma$ with respective neighborhoods $U,V$ whose intersection with $\Gamma$ is $E$-connected, the intersection $\left(U\cup V\right)\cap \Gamma$ is also $E$-connected.
Since $\partial \Gamma$ is compact one can find a cofinite $E$-connected subset of $\Gamma$.

\medskip

Since $A,B\neq \emptyset$, the definition of the limit cone allows us to pick two large enough elements $\gamma',\gamma''\in\Gamma\setminus E$ such that :
 $$d(\left[\da(o,\rho(\gamma')\cdot o)\right],A)\leq\frac{\epsilon}{3},$$
 $$d(\left[\da(o,\rho(\gamma'')\cdot o)\right],B)\leq\frac{\epsilon}{3}.$$
 
 One can assume that $\gamma',\gamma''$ are $E$-connected by taking them large enough. Therefore there exist $\gamma\in \Gamma\setminus E$ and $a,b\in F$ such that:
  $$d(\left[\da(o,\rho(\gamma)\cdot o)\right],A)\leq\frac{\epsilon}{3},$$
  $$d(\left[\da(o,\rho(a\gamma b)\cdot o)\right],B)\leq\frac{\epsilon}{3}.$$

 Using Lemma \ref{lem:GeneralizedDistanceLip} one gets that :
$$|\da(o,\rho(\gamma)\cdot o)-\da(o,\rho(a\gamma b)\cdot o)|\leq |\da(o,\rho(\gamma)\cdot o)-\da(o,\rho(a\gamma )\cdot o)|+d_\X( o,\rho(b)\cdot o).$$ 
$$|\da(\rho(\gamma)^{-1}\cdot o, o)-\da(\rho(\gamma)^{-1}\rho(a)^{-1}\cdot o, o)|\leq d_\X(\rho(a)^{-1}\cdot o,o).$$

And therefore:

$$\left|\frac{\da(o,\rho(\gamma)\cdot o)}{d_\X(o,\rho(\gamma)\cdot o)}-\frac{\da(o,\rho(a\gamma b)\cdot o)}{d_\X(o,\rho(a\gamma b)\cdot o)}\right|< \frac{\epsilon}{3}.$$

This contradicts the fact that $A$ and $B$ are at distance $\epsilon$, so $\mathcal{C}_\rho$ is connected.

\end{proof}

The assumption that $\Gamma$ is non-elementary is necessary, see the end of Section \ref{sec:AnosovNearlyFuch}.

\subsection{Anosov representations.}

The Anosov properties are more restrictive for a representation than the property of being a quasi-isometric embedding. These notions are interesting in high rank because the Anosov properties hold for an open set of representations, whereas the property of being a quasi-isometric embedding is not necessarily open in $\Hom(\Gamma_g,G)$ when the rank of $\X$ is at least $2$.

\begin{defn}[{\cite{BPS}, Section 4}]
\label{defn:Anosov}
Let $\Theta\subset \Delta$ be a non-empty set of simple roots. A representation $\rho:\Gamma \to G$ is \emph{$\Theta$-Anosov} if for every root $\alpha\in \Theta$ there exists some constants $b,c>0$ such that for every $\gamma\in \Gamma$:
$$\alpha\left(\da(o,\rho(\gamma)\cdot o)\right)\geq b|\gamma|_F-c.$$

This definition does not depend on the choice of the generating set $F$ and the base-point $o$.

\medskip

A $\Delta$-Anosov representation in the case when $G$ is a split real simple Lie group is be called a \emph{Borel-Anosov} representation.

\end{defn}

\begin{rem}
A representation is $P$-Anosov for a parabolic subgroup $P$ if it is $\Theta$-Anosov for the corresponding set of simple roots $\Theta\subset \Delta$.
\end{rem}

Let $\alpha\in \Delta$ and $\typeS\in \Sph \mathfrak{a}$ be orthogonal to $\Ker(\alpha)$. The evaluation of $\alpha$ to $\da(x,y)$ satisfies:
$$\alpha(\da(x,y))=\alpha(\type_0)d_\X(x,y)\cos\left(\angle(\da(x,y),\type_0)\right) .$$

Anosov representations are necessarily quasi-isometric embeddings. Reciprocally a quasi-isometric embedding is $\lbrace \alpha\rbrace $-Anosov if and only if the angle $\langle\da(o,\rho(\gamma)\cdot o),\type_0\rangle$ is not too small in absolute value for $\gamma\in \Gamma_g$ large enough. 

\medskip

In particular we have the following characterization of Anosov representations:

\begin{thm}[\cite{KLP}]
\label{thm:CaracAnosocQuasiIsom}

A representation $\rho:\Gamma_g\to G$ is $\Theta$-Anosov for $\Theta\subset\Delta$ if and only if it is a quasi-isometric embedding and if $\Ker(\alpha)\cap \mathcal{C}_\rho=\emptyset$ for all $\alpha\in \Theta$.
\end{thm}

\medskip

 Representations that are $\Theta$-Anosov admit a natural continuous and $\rho$-equivariant map $\xi^\Theta_\rho:\partial\Gamma_g\to \mathcal{F}_\Theta=G/\type_\Theta$, where $\partial\Gamma_g$ is the Gromov boundary of $\Gamma_g$. 
 
 \medskip
 
In the proof of Theorem \ref{thm:ComparaisonTotallyGeodesic} we will use the following results about the boundary maps of Anosov representations.
For two points $o,x\in \X$ let $\ell(o,x)\in \partialvis \X$ be the class of the unique geodesic ray with unit speed starting from $o$ and passing through $x$.

\begin{thm}[{\cite{BPS}, Section 4}]
\label{thm:AnosovLimitMap}
Let $\rho:\Gamma_g\to G$ be a $\Theta$-Anosov representation for a non-empty set $\Theta\subset \Delta$. There exist a unique $\rho$-equivariant continuous and dynamic preserving map $\xi^\Theta_\rho:\partial\Gamma_g\to \mathcal{F}_\Theta$. 
This map is such that for any $o\in \X$ and any sequence $(\gamma_n)_{n\in \N}$ of elements of $\Gamma_g$ converging to $\zeta\in \partial \Gamma_g$, the $\Delta$-facet containing any limit point of the sequence $(\ell(o, \rho(\gamma_n)\cdot o))_{n\in \N}$ also contains the $\Theta$-facet $\xi^\Theta_\rho(\zeta)$.
\end{thm}

For instance, when $G=\PSL(n,\R)$ and if $\Theta=\lbrace \alpha_k\rbrace$, one can associate a partial flag to any point in $\partialvis \X$. If the the representation $\rho$ is $\lbrace \alpha_k\rbrace$-Anosov, the partial flag associated to any limit point of $(\ell(o, \rho(\gamma_n)\cdot o))_{n\in \N}$ contains the same $k$-dimensional plane, that will be denoted by $\xi^k_\rho(\zeta)$.

\medskip

Kapovich, Leeb and Porti also proved a generalization of the Morse lemma. Here is a version of this result. Let us fix any metric on $\Gamma$ quasi-isometric to a word metric.

\begin{thm}[{\cite{KLPMorse}, Theorem 1.3}]
\label{thm:MorseLemma}
Let $\rho:\Gamma\to G$ be a $\Theta$-Anosov representation. Let $o\in \X$ be a base-point. There exist a constant $D>0$ such that for every $\gamma\in \Gamma$, there exist a geodesic ray $\eta:\R_{>0}\to \X$ at distance at most $D$ from $\rho(\gamma)\cdot o$ with $\eta(0)=o$, whose class $[\eta]\in \partialvis\X$ lies in a common $\Delta$-facet with $\xi^\Theta_\rho(\zeta_\gamma)$. Here $\zeta_\gamma\in \partial \Gamma$ is the endpoint of any geodesic
ray in $\Gamma$ starting at the identity and going through $\gamma$.
\end{thm}

\section{Busemann functions on symmetric spaces.}
\label{sec:BusemannFunctions}

Busemann functions are natural functions on Hadamard manifolds associated to points in the visual boundary. These functions will play a key role in the definition of $\type$-nearly geodesic immersions, and in the fibration of domains of discontinuity. In this section we prove the main properties of Busemann functions and compute their Hessian.

\subsection{Main properties of Busemann functions.}

 Busemann functions can be interpreted as the distance of a point $x\in \X$ to a point $a$ in the visual boundary relative to a base-point $o\in \X$.

\begin{defn}
 The \emph{Busemann function} associated to $a\in \partialvis\X$ and based at $o\in \X$ is the map $b_{a,o}:X\to \R$ that associates to $x\in \X$ the limit :
 $$\lim_{t\to +\infty} d_\X(x,\gamma(t))-d_\X(o,\gamma(t)),$$
 
 for any geodesic ray $\gamma:\R^+\to \X$ in the class of $a$.
 
\end{defn}

This definition makes sense because $\X$ is a Hadamard manifold \cite{Eberlein}. The definition implies that for any $x,o,o'\in X$ and $a\in \partialvis \X$, the Busemann cocycle holds:
\begin{equation}
\label{eq:BusemannCocycle}
b_{a,o'}(x)=b_{a,o}(x)+b_{a,o'}(o).
\end{equation}

\medskip

For symmetric spaces, this function can be computed using the generalized Iwasawa decomposition. First we prove that unipotent elements preserve the level lines of Busemann functions.

\begin{lem}

Let $x,o\in \X$ and $a\in \partialvis \X$ be two points. Let $n$ be an element of the unipotent subgroup $N_{a,o}$ of $G$:
$$ b_{a,o}(n\cdot x)=b_{a,o}(x).$$

\end{lem}

\begin{proof}

The Busemann cocycle implies that $b_{a,o}(n\cdot x)-b_{a,o}(x)=b_{a,x}(n\cdot x)$. This is by definition the limit when $t\to \infty$ of the difference:
\begin{flalign*}
d_\X\left(n\cdot x,\exp(t\mathrm{v}_{a,x})\cdot x\right)&-d_\X\left(x,\exp(t\mathrm{v}_{a,x})\cdot x\right)\\
&= d_\X\left(x,n^{-1}\exp(t\mathrm{v}_{a,x})\cdot x\right)-d_\X\left(x,\exp(t\mathrm{v}_{a,x})\cdot x\right)\\
&\leq d_\X\left(n^{-1}\exp(t\mathrm{v}_{a,x})\cdot x,\exp(t\mathrm{v}_{a,x})\cdot x\right).
\end{flalign*}
But since $n\in N_{a,x}$, $\exp(-t\mathrm{v}_{a,x})n\exp(t\mathrm{v}_{a,x})$ converges to the identity when $t\to +\infty$, so this distance converges to $0$. Hence $b_{a,x}(n\cdot x)=0$.

\end{proof}
 
Recall that $\mathrm{v}_{a,o}$ is the unit vector in $T_o\X$ pointing towards $a\in \partialvis\X$. To compute a Busemann function one needs to understand it on maximal flats. Let $x=\exp(\mathrm{w})\cdot o$ for $\mathrm{w}\in \mathfrak{p}_o$. Suppose that $a,o,x$ lie in the same flat subspace, \ie $[\mathrm{w},\mathrm{v}_{a,o}]=0$.
The Busemann function on this Euclidean space is equal to:
$$b_{a,o}(x)=-d_\X(x,o)\cos(\angle_o(a,x))= \langle-\mathrm{v}_{a,x},\mathrm{w} \rangle_x.$$

Using these facts we can write Busemann functions in the symmetric space $\X$ explicitely. Let $o,x\in \X$ be a base point and $a\in \partial_{\vis} \X$. 

\begin{cor}
Let $o,x\in X$ and $a\in \partialvis \X$. The Busemann function can be computed as :
$$ b_{a,o}(x)=\langle-\mathrm{v}_{a,o},\mathrm{w}\rangle_o.$$

Where $\mathrm{w}\in \mathfrak{a}_{a,o}$ is given by the generalized Iwasawa decomposition, \ie  is the unique element such that one can write $x=n\exp(\mathrm{w})k\cdot o$ with $n\in N_{a,o}$ and $k\in K_o$.
\end{cor}

Since $G$ acts by isometries on $\X$, Busemann functions are $G$-equivariant in the following sense.

\begin{cor}
\label{cor:EquivariantBusemann}
Let $o\in \X$ and $a\in \partial_{\vis} \X$. For any $g\in G$, and any $x\in \X$,
$b_{g\cdot a,g\cdot o}(g\cdot x)=b_{a,o}(x)$.
\end{cor}

The gradient of Busemann functions is characterized as follows.

\begin{prop}
The gradient of the Busemann function based at any point $o\in \X$ associated to $a\in \partial_{\vis}\X$ is the vector field  $(-\mathrm{v}_{a,x})_{x\in \X}$ of unit vectors pointing towards $a$.
\end{prop}

\begin{proof}
The differential $\mathrm{d}_xb_{a,o}$ of $b_{a,x}$ at $x$ associates to an element $\mathrm{w}\in \mathfrak{p}_x$ the value $\langle\mathrm{w}',-\mathrm{v}_{a,x}\rangle_x$ where $\mathrm{w}'$ is the projection of $\mathrm{w}$ to $\mathfrak{a}_{a,x}$ with respect to the decomposition $\mathfrak{g}=\mathfrak{n}_{a,x}\oplus \mathfrak{a}_{a,x}\oplus \mathfrak{k}_{x}$. Note that $\mathfrak{n}_{a,x}$ and $\mathfrak{k}$ are orthogonal to $\mathrm{v}_{a,o}\in \mathfrak{a}_{a,x}$ with respect to $\langle \cdot ,\cdot \rangle_x$. Hence $\mathrm{w}=\mathrm{w}'$ so the gradient of $b_{a,o}$ at $x$ is $-\mathrm{v}_{a,x}$.
\end{proof}

Busemann functions vary smoothly when the base flag varies in a flag manifold.

\begin{lem}
\label{lem:SmoothBusemann}
For any $o\in \X$, and $\type\in \Sph \mathfrak{a}^+$. The map $\mathcal{F}_\type \times \X\to \mathbb{R}$, $(a,x)\mapsto b_{a,o}(x)$ is smooth.
\end{lem}

\begin{proof}
Let $P$ be the stabilizer of an element $a\in \mathcal{F}_\type$. By Corollary \ref{cor:EquivariantBusemann}, for $g\in G$ and $x,y\in \X$, $b_{g\cdot a_0,o}(y)=b_{a_0,o}(g^{-1}\cdot y)-b_{a_0,o}(g^{-1}\cdot o)$.

\medskip

 Hence the map $G \times \X\to \mathbb{R}$, $(g,x)\mapsto b_{g\cdot a_0,o}(x)$ is smooth, and defines a smooth map from the quotient $G/P \times X\simeq \mathcal{F}_\type \times \X$.
 
\end{proof}

\begin{exmp}

Let $\X=\mathcal{S}_n$ or $\mathcal{H}_n$ the symmetric space associated with $\PSL_{n}(\K)$ with $\K=\R$ or $\C$, as in Section \ref{subsec:ExamplesRoots}. Let $(e_1,\cdots,e_n)$ be a basis of $\K^n$. The projective space $\mathbb{P}(\K^n)$ can be identified with the $G$-orbit $\mathcal{F}_{\type_1}$ of the point $a\in \partialvis \X$ corresponding to the limit point where $t$ goes to $+\infty$ of the geodesic ray:

$$t\mapsto
\begin{pmatrix}
e^{-t(n-1)} & 0 & \cdots& 0\\
0 & e^t &  & 0\\
\cdots &   & \cdots&  \\
0 & 0 &  &  e^t
\end{pmatrix} $$

The point $a\in \mathcal{F}_\type\simeq\mathbb{RP}^{n-1}$ is identified with the first basis vector since the stabilizer of both points by the respective actions of $\PSL_n(\K)$ on are equal.

\medskip

The Busemann function $b_{[v],q_0}$ where $q_0\in \X$ and $[v]\in \mathbb{P}(\K^n)$ associates to $q\in \X$ the value $\sqrt{\frac{n-1}{n}} \log(q(v,v))$ where $v$ is a representative of $[v]$ such that $q_0(v,v)=1$.
\end{exmp}

The asymptotic behavior of Busemann functions along geodesic rays is determined by the Tits angle between the endpoints.

\begin{lem}
\label{lem:asymptoticsBusemann}
Let $a\in \partial_{\vis}\X$ and $x\in \X$. Let $\eta$ be a geodesic ray converging to $b\in \partial_{\vis}\X$. Then there exists a constant $C>0$ such that for all $t\in \R$ :

$$ \left|b_{a,x_0}(\eta(t))+t\cos(\Tangle(a,b))\right|\leq C$$
\end{lem}

\begin{proof}
There exist some element $g\in G$ such that $g\cdot a$ and $g\cdot b$ belong to $\partial_{\vis}F$ with $F$ the model flat in $\X$. Moreover there exist a geodesic ray $\eta'$ at bounded distance from $g\cdot \eta$ that belongs to the flat $F$. On the flat subspace $F$, the Busemann function can be computed:
$$b_{a,\eta'(0)}(\eta'(t))=-t\cos(\Tangle(g\cdot a,g\cdot b)).$$ 

This proves the lemma.
\end{proof}

\subsection{Computation of the Hessian.}

We compute here the Hessian of Busemann functions in the symmetric space $\X$. This computation will be used in the proof of Theorem \ref{thm:Sufficient condition}.

\begin{lem}
\label{lem:HessianofBusemannFunctions}
Let $a\in \partial_{\vis}\X$, and $x,o\in \X$. The Hessian of the Busemann function $b_{a,o}$ at a point $x\in \X$ is given by the following quadratic form on $T_x\X$:

\begin{equation}
\label{eq:HessianBusemann}
\mathrm{v}\mapsto  \left\langle \sqrt{\ad_{\mathrm{v}_{a,x}}^2}(\mathrm{v}) ,\mathrm{v}\right\rangle_x.
\end{equation}

Here $\sqrt{\ad_{\mathrm{v}_{a,x}}^2}$ is the only root of the endomorphism ${\ad_{\mathrm{v}_{a,x}}\circ \ad_{\mathrm{v}_{a,x}}}_{|\mathfrak{p}_x}:\mathfrak{p}_x\to \mathfrak{p}_x$ that is symmetric and semi-positive for the scalar product $\langle\cdot,\cdot\rangle_x$.

\medskip

This quadratic form is semi-positive, and vanishes exactly on $\mathfrak{z}(\mathrm{v}_{a,x})\cap \mathfrak{p}_x$. For $\mathrm{v}\in (\mathfrak{z}(\mathrm{v}_{a,x})\cap \mathfrak{p}_x)^\perp$, it satisfies: 

\begin{equation}
\label{eq:HessianBusemannInequality}
\Hess_x(\mathrm{v},\mathrm{v})\geq \lVert\mathrm{v}\rVert^2\min_{\alpha \in \Sigma,\alpha(\Cartan(a))\neq 0}|\alpha(\Cartan(\mathrm{v}_{a,x}))|.
\end{equation}
\end{lem}

Recall that $\mathfrak{z}(\mathrm{v})$ for $\mathrm{v}\in \mathfrak{g}$ is the centralizer in $\mathfrak{g}$ of $\mathrm{v}$.

\begin{rem}
\label{rem:secCurvature}
The Hessian of a Busemann function is related to the sectional curvature of the symmetric space. When measured along a tangent plane spanned by two orthogonal unit vectors $\mathrm{v},\mathrm{w}\in T_o\X\simeq \mathfrak{p}_o\subset \mathfrak{g}$ the sectional curvature of $\X$ is equal to:
$$\kappa_{\mathrm{v},\mathrm{w}}=-\langle [\mathrm{v},[\mathrm{v},\mathrm{w}]], \mathrm{w} \rangle_o=-\langle \ad_{\mathrm{v}}^2(\mathrm{w}), \mathrm{w} \rangle_o.$$
\end{rem}

This Lemma implies that Busemann functions are strictly convex except on flats. This is a more general fact about Hadamard manifolds, see \cite{Eberlein}.

\medskip

\begin{proof}

The Busemann function with respect to two different base points differ only by a constant. Hence we can assume here without any loss of generality that $x=o$.

Let $\mathrm{v}\in T_o\X$ be a vector. The generalized Iwasawa decomposition, and the fact that the exponential map is a local diffeomorphism on Lie group implies that there exists a neighborhood $I$ of $0$ in $\mathbb{R}$ such that for all $t\in I$:
\begin{equation}
\label{eq:HessianCalcul}
\exp(t\mathrm{v})=\exp(\mathrm{n}_t)\exp(\mathrm{w}_t)\exp(\mathrm{k}_t)
 \end{equation}
 for $\mathrm{n}_t\in \mathfrak{n}_{a,o}$, $\mathrm{w}_t\in \mathfrak{a}_{a,x}$ and $\mathrm{k}_t\in \mathfrak{t}_{x}$, and so that the map $t\mapsto (\mathrm{n}_t,\mathrm{w}_t,\mathrm{k}_t)$ is smooth. Let us denote by $(\dot{\mathrm{n}},\dot{\mathrm{w}},\dot{\mathrm{k}})$ and $(\ddot{\mathrm{n}},\ddot{\mathrm{w}},\ddot{\mathrm{k}})$ the first and second derivative of this map at $t=0$.

\medskip

The limited development at order $2$ at $t=0$ of \eqref{eq:HessianCalcul} yields:
$$\exp(t\mathrm{v})=\exp\left(\dot{\mathrm{n}}t+\frac{\ddot{\mathrm{n}}}{2}t^2\right)\exp\left(\dot{\mathrm{w}}t+\frac{\ddot{\mathrm{w}}}{2}t^2\right)\exp\left(\dot{\mathrm{k}}t+\frac{\ddot{\mathrm{k}}}{2}t^2\right)+o(t^2).$$

 But the Baker–Campbell–Hausdorff formula \cite{Helgason} implies that the right hand of this equality is equal to:
 $$\exp\left(  \dot{\mathrm{n}}t+\dot{\mathrm{w}}t+\dot{\mathrm{k}}t+  \frac{\ddot{\mathrm{n}}}{2}t^2+\frac{\ddot{\mathrm{w}}}{2}t^2+\frac{\ddot{\mathrm{k}}}{2}t^2+\frac{1}{2}\left([\dot{\mathrm{n}},\dot{\mathrm{w}}]+[\dot{\mathrm{n}},\dot{\mathrm{k}}]+[\dot{\mathrm{w}},\dot{\mathrm{k}}]\right)t^2+o(t^2)\right).$$
 
 Hence we get the following two equalities:
$$\mathrm{v}=\dot{\mathrm{n}}+\dot{\mathrm{w}}+\dot{\mathrm{k}},$$
$$0=\frac{\ddot{\mathrm{n}}}{2}+\frac{\ddot{\mathrm{w}}}{2}+\frac{\ddot{\mathrm{k}}}{2}+\frac{1}{2}\left([\dot{\mathrm{n}},\dot{\mathrm{w}}]+[\dot{\mathrm{n}},\dot{\mathrm{k}}]+[\dot{\mathrm{w}},\dot{\mathrm{k}}]\right).$$

However since $\mathrm{v}, \dot{\mathrm{w}}\in \mathfrak{p}_x$, then $\theta_x(\dot{\mathrm{n}}+\dot{\mathrm{k}})=-\dot{\mathrm{n}}-\dot{\mathrm{k}}$. Hence $\dot{\mathrm{k}}=-\frac{\dot{\mathrm{n}}+\theta_x(\dot{\mathrm{n}})}{2}$. This let us simplify the last part of the previous equation :
$$[\dot{\mathrm{n}},\dot{\mathrm{w}}]+[\dot{\mathrm{n}},\dot{\mathrm{k}}]+[\dot{\mathrm{w}},\dot{\mathrm{k}}]=[\dot{\mathrm{n}},\dot{\mathrm{w}}]-\frac{1}{2}[\dot{\mathrm{n}},\theta_x(\dot{\mathrm{n}})]-\frac{1}{2}[\dot{\mathrm{w}},\dot{\mathrm{n}}+\theta_x(\dot{\mathrm{n}})]$$ 

The metric on $\X$ can be written $\langle \cdot,\cdot\rangle_x=B(\cdot, \theta_x(\cdot))$ on $\mathfrak{p}_x$ with $B$ the Killing form, defined on $\mathfrak{g}$.

Since $\mathrm{v}_{a,x}$ is orthogonal to $\mathfrak{n}_{a,o}$ and $\mathfrak{t}_x$ then $B(\mathrm{v}_{a,x},\ddot{\mathrm{n}})=B(\mathrm{v}_{a,x},\ddot{\mathrm{k}})=0$. Moreover $[\dot{\mathrm{n}},\dot{\mathrm{w}}]\in \mathfrak{n}_{a,o}$ so:
$$B(\mathrm{v}_{a,x},[\dot{\mathrm{w}},\dot{\mathrm{n}}+\theta_x(\dot{\mathrm{n}})])=B(\mathrm{v}_{a,x},[\dot{\mathrm{n}},\dot{\mathrm{w}}])=0.$$

In particular one gets:
$$\Hess_x(b_{a,o})(\mathrm{v},\mathrm{v})=\langle -\mathrm{v}_{a,x},\ddot{\mathrm{w}} \rangle_x=\frac{1}{2}B(-\mathrm{v}_{a,x},[\dot{\mathrm{n}},\theta_x(\dot{\mathrm{n}})])=\frac{1}{2}B([-\mathrm{v}_{a,x},\dot{\mathrm{n}}],\theta_x(\dot{\mathrm{n}})).$$

Let $\Sigma_a\subset \Sigma$ be the set of roots  $\alpha$ such that $\alpha(\Cartan(a))\neq 0$. The Lie algebra decomposes into root spaces:
$$\mathfrak{g}=\mathfrak{z}(\mathrm{v}_{a,x})\oplus\bigoplus_{\alpha\in \Sigma_a}\mathfrak{g}_{a,x}^\alpha,$$

where $\Ad_g(\mathfrak{g}_{a,x}^\alpha)=\mathfrak{g}^\alpha$ the model root spaces with $g\in G$ any element such that $g\cdot \mathrm{v}_{a,x}=\Cartan(\mathrm{v}_{a,x})$.

The restriction of $\ad_{\mathrm{v}_{a,x}}$ on $\mathfrak{g}_{a,x}^\alpha$ is an homothety of ratio $\alpha(\Cartan(\mathrm{v}_{a,x}))$. The vector $\mathrm{v}$ can be decomposed in this direct sum.
$$\mathrm{v}=\mathrm{v}^0+\sum_{ \alpha\in \Sigma_a}\mathrm{v}^\alpha.$$

The endomorphism $\sqrt{\ad_{\mathrm{v}_{a,x}}^2}$ associates to $\mathrm{v}$ the vector :
$$\sum_{ \alpha\in \Sigma_a}|\alpha(\Cartan(\mathrm{v}_{a,x}))|\mathrm{v}^\alpha\in \mathfrak{p}_x.$$

Let $\Sigma^{+}_a$ be the set of roots  $\alpha$ such that $\alpha(\Cartan(\mathrm{v}_{a,x}))> 0$. The vector $\dot{\mathrm{n}}$ can be expressed as :

$$\dot{\mathrm{n}}= 2\sum_{\alpha\in \Sigma_a^{+}}\mathrm{v}^\alpha.$$

Hence, we get as desired:
\begin{equation} \label{eq1}
\begin{split}
	\frac{1}{2} B([-\mathrm{v}_{a,x},\dot{\mathrm{n}}],\theta_x(\dot{\mathrm{n}}))
		&  
	= -2\sum_{\alpha\in \Sigma^+_\alpha}\alpha(\Cartan(\mathrm{v}_{a,x}))B(\mathrm{v}^\alpha, \theta_x(\mathrm{v}^\alpha))
		\\
		& 
	=\sum_{\alpha\in \Sigma_\alpha}|\alpha(\Cartan(\mathrm{v}_{a,x}))|\langle \mathrm{v}^\alpha,\mathrm{v}^\alpha\rangle_x 
	=\left\langle\sqrt{\ad_{\mathrm{v}_{a,x}}^2}(\mathrm{v}),\mathrm{v}\right\rangle_x.
\end{split}
\end{equation}

\medskip

This is equal to zero if and only if $\mathrm{v}=\mathrm{v}^0$.
\end{proof}

\section{Nearly geodesic immersions.}

In this section we introduce a local condition for an immersion into the symmetric space of non-compact type $\X$ that generalizes the notion an immersion with principal curvature in $(-1,1)$ inside $\mathbb{H}^n$. 

\label{sec:nearlygeodesic}
\subsection{Curvature bound and Busemann functions.}

We introduce the key definition of a nearly geodesic immersion, which relies on Busemann functions (see Section \ref{sec:BusemannFunctions}).
Let $M$ be a smooth connected manifold, $u:M\to \X$ be an immersion, $o\in \X$ a base point and let $\type\in \Sph\mathfrak{a}^+$ be a unit vector in the model Weyl chamber.

\begin{defn}

\label{Defn:NearlyGeodesicSurface}
An immersion $u: M\to \X$ is called \emph{$\type$-nearly geodesic} if for all $a\in \mathcal{F}_\type\cup \mathcal{F}_{\iota(\type)}$ and $\mathrm{v}\in TM$ such that $\mathrm{d}(b_{a,o}\circ u)(\mathrm{v})=0$, the function $b_{a,o}\circ u$ has positive Hessian in the direction $\mathrm{v}$.
\end{defn}

The Hessian considered in this definition is computed with the induced metric $u^*g_\X$ on $M$. Recall that $\mathcal{F}_{\iota(\type)}$ is the opposite flag manifold to $\mathcal{F}_{\type}$ for $\type\in \Sph\mathfrak{a}^+$.

\medskip

We will first show that the nearly geodesic condition can be written as a bound on the fundamental form $\II_u$, depending on the Cartan projection of the surface tangent vectors. 

\medskip

Since the Hessian of a Busemann function $b_{a,o}$ on $\X$ does not depend on $o$, we will denote it by $\Hess_{b_a}$. Recall that $\mathrm{v}_{a,o}$ is the unit vector in $T_o\X$ pointing towards $a\in \partialvis \X$. The second fundamental form $\II_u$ for $x\in M$ of the immersion $u$ is the difference $u^*\nabla^\X-\nabla^M$ where $\nabla^\X$ is the Levi-Civita connection on $T\X$ associated to $g_\X$ and $\nabla^M$ is the Levi-Civita connection on $TM\subset u^*T\X$ associated to the metric $u^*g_\X$. The second fundamental form is a symmetric $2$-tensor with values in $u^*\mathrm{N}$ when $\mathrm{N}\subset T\X$ is the normal tangent bundle to $u(M)$.

\begin{prop}
\label{prop:nearlyfuchsianLocalProp}
An immersion $u:M\to \X$ is $\type$-nearly geodesic if and only if for all $y\in M$, for all $a\in \mathcal{F}_\type\cup \mathcal{F}_{\iota(\type)}$ and $\mathrm{v}\in T_yM$ such that $\langle \mathrm{d}u(\mathrm{v}),\mathrm{v}_{a,u(y)}\rangle_{u(y)}=0$:
\begin{equation}
\label{eq:HessianNearlyfuchsian}
\Hess_{b_a}(\mathrm{d}u(\mathrm{v}),\mathrm{d}u(\mathrm{v}))+\langle\II_u(\mathrm{v},\mathrm{v}),\mathrm{v}_{a,u(y)}\rangle_{u(y)}>0. \end{equation}
\end{prop}

We will prove a sufficient condition that has a simpler form in Theorem \ref{thm:Sufficient condition} when $\type=\type_\Theta$ for a Weyl orbit of simple roots $\Theta$.

\medskip

\begin{proof}

Let $y\in M$ and $a\in \mathcal{F}_\type\cup \mathcal{F}_{\iota(\type)}$. The function $b_{a,o}\circ u$ is critical at $y$ in the direction $\mathrm{v}\in T_yM$ if and only if $\langle \mathrm{d}u(\mathrm{v}),\mathrm{v}_{a,u(y)}\rangle_{u(y)}=0$. 

\medskip

Let $\gamma:\R\to M$ be a geodesic for the metric $u^*g_\X$ on $M$ such that $\gamma(0)=y$ and $\gamma'(0)=\mathrm{v}$. The Hessian of $b_{a,o}\circ u$ on $M$ is equal to the derivative at $t=0$ of the differential of the Busemann function, \ie:
$$t\mapsto \langle \mathrm{v}_{a,u(\gamma(t))},\mathrm{d}u(\gamma'(t))\rangle_{u(y)} $$
The first term, $\langle\nabla^\X_{\mathrm{d}u(\mathrm{v})}\mathrm{v}_{a,u(\gamma(t))},\mathrm{d}u(\gamma'(t))\rangle_{u(y)}$, is equal to $\Hess_{b_a}(\mathrm{d}u(\mathrm{v}),\mathrm{d}u(\mathrm{v}))$. The second term can be written:
$$\nabla^\X_{\mathrm{d}u(\mathrm{v})}\mathrm{d}u(\gamma')=u_*\nabla^M_{\mathrm{v}}\mathrm{d}u(\gamma')+\II_u(\mathrm{v},\mathrm{v})$$
But $\gamma$ is a geodesic so $\nabla^M_{\mathrm{v}}\mathrm{d}u(\gamma')=0$, therefore the Hessian of $b_{a,o}\circ u$ on $M$ in the direction $\mathrm{v}$ is equal to :
$$\Hess_{b_a}(\mathrm{d}u(\mathrm{v}),\mathrm{d}u(\mathrm{v}))+\langle\II_u(\mathrm{v},\mathrm{v}),\mathrm{v}_{a,u(y)}\rangle_{u(y)}>0. $$
\end{proof}

A consequence of Proposition \ref{prop:nearlyfuchsianLocalProp} is that the property of being $\type$-nearly geodesic is locally an open property for the $\mathcal{C}^2$-topology, which is the topology associated with the uniform convergence over any compact set of the first two differentials.

\begin{cor}
\label{cor:nearlyfuchsianOpen}

Let $u_0: M\to \X$ be a $\type$-nearly geodesic map for some $\type\in \Sph\mathfrak{a}^+$.  For all compact $K\subset M$, there exists a neighborhood $U$ of $u_0$ for the $\mathcal{C}^2$-topology in the space of $\mathcal{C}^2$ maps from $M$ to $\X$  and a neighborhood $V$ of $\type$ in $\Sph\mathfrak{a}^+$ such that for all $\type'\in V$ and $u\in U$, $u$ satisfies the $\type'$-nearly geodesic immersion condition on $K$.
\end{cor}

Let $G$ be the isometry group of the $n$-dimensional hyperbolic space $\mathbb{H}^n$ for some $n\in \mathbb{N}$ with its usual metric with sectional curvature equal to $-1$. We prove that the notion of $\type$-nearly geodesic immersion generalizes the notion of immersion with principal curvatures in $(-1,1)$ in $\mathbb{H}^n$. Principal curvatures are only defined for hypersurfaces, but the following definition allows to generalize the notion of having bounded principal curvature.

\begin{defn}
An immersion $u:M\to \mathbb{H}^n$ has \emph{principal curvature in $(-1,1)$} if and only if for all $\mathrm{v}\in TM$, $\lVert \mathrm{d}u(\mathrm{v})\rVert^2> \lVert\II_u(\mathrm{v},\mathrm{v})\rVert$.
\end{defn}

Since $\mathbb{H}^n$ is a rank one symmetric space, $\Sph\mathfrak{a}^+$ contains a single element.

\begin{prop}
\label{prop:nearlyfuchsianH3}
An immersion $u:M\to \mathbb{H}^n$ is nearly geodesic for the only element $\type\in \Sph\mathfrak{a}^+$ if and only $u$ has principal curvature in $(-1,1)$.
\end{prop}

\begin{proof}
Let $x\in \X=\mathbb{H}^n$ and $a\in \mathcal{F}_\type=\mathcal{F}_{\iota(\type)}=\mathbb{CP}^1=\partial\mathbb{H}^n$. For any $\mathrm{w}\in T_x\X$, $\Hess_{b_a}(\mathrm{w},\mathrm{w})=\lambda\lVert\mathrm{w}^\perp\rVert^2$ where $\mathrm{w}^\perp$ is the orthogonal projection of $\mathrm{w}$ onto the orthogonal in $T_x\X$ of $\mathrm{v}_{a,x}$ by Proposition \ref{lem:HessianofBusemannFunctions}, with some constant $\lambda$  which is equal to $1$ for the metric of sectional curvature equal to $-1$ on $\mathbb{H}^n$ (see Remark \ref{rem:secCurvature}).

If $u$ is $\type$-nearly geodesic, then it is an immersion and for every $y\in M$ and $\mathrm{v}\in T_yM$ there exist $a\in \partial \mathbb{H}^n$ such that $\mathrm{v}_{a,u(y)}$ is positively collinear with $-\II_u\left(\mathrm{d}u(\mathrm{v}),\mathrm{d}u(\mathrm{v})\right)$. By Proposition \ref{prop:nearlyfuchsianLocalProp}, and since $\mathrm{v}\perp\mathrm{v}_{a,u(y)}$ one has:
$$\lVert \mathrm{d}u(\mathrm{v})\rVert^2-\lVert\II_u(\mathrm{v},\mathrm{v})\rVert>0.$$

Therefore the principal curvature of $u$ is in $(-1,1)$.

\medskip

Conversely if $u$ is an immersion with principal curvatures in $(-1,1)$, let $a\in \partial\mathbb{H}^n$, $y\in M$ and $\mathrm{v}\in T_yM$ be such that $b_{a,o}\circ u$ is critical in the direction $\mathrm{v}$. Hence $\mathrm{v}_{a,u(y)}$ is perpendicular to $\mathrm{d}u(\mathrm{v})$ so $\Hess_{b_a}(\mathrm{d}u(\mathrm{v},\mathrm{d}u(\mathrm{v}))=\lVert \mathrm{d}u(\mathrm{v})\rVert^2$. Therefore the fact that $u$ has principal curvature in $(-1,1)$ implies that the hypothesis of Proposition \ref{prop:nearlyfuchsianLocalProp} hold, so $u$ is $\type$-nearly geodesic.
\end{proof}

In general, the property of being $\type$-nearly geodesic implies that the surface is regular in the following sense.

\begin{defn}
\label{defn:pRegularVector}
A tangent vector $\mathrm{v}\in T\X$ is called \emph{$\type$-regular} if its Cartan projection $\Cartan(\mathrm{v})$ does not belong to $\bigcup_{w\in W}(w\cdot \type)^\perp$.
\end{defn}

We say that an immersion $u:M\to \X$ is $\type$ regular if for all $\mathrm{v}\in TM$, $\mathrm{d}u(\mathrm{v})$ is $\type$-regular.

\medskip

 Being regular, namely having the Cartan projection in the interior of $\mathfrak{a}^+$, and being $\type$-regular is in general unrelated. However when $\type=\type_\Theta$ for a Weyl orbit of simple roots $\Theta$, a $\type$-regular vector $\mathrm{v}\in T\X$ is exactly a $\Theta$-regular vector, namely such that for all $\alpha\in \Theta$, $\alpha\left(\Cartan(\mathrm{v})\right)\neq 0$.

\begin{prop}
\label{prop:pRegularSurface}
Let $\type\in \Sph\mathfrak{a}^+$. If $u$ is a $\type$-nearly geodesic immersion, the tangent vectors $\mathrm{d}u(\mathrm{v})$ for $\mathrm{v}\in TM$ are $\type$-regular.
\end{prop}

\begin{proof}
Let $\mathrm{v}\in T_yM$ for some $y\in M$. Assume that $\mathrm{d}u(\mathrm{v})$ is not $\type$-regular, so its Cartan projection is orthogonal to $w\cdot \type$ for some $w\in W$. Therefore there is a unit vector which lies in a common maximal flat with $\mathrm{d}u(\mathrm{v})$, and whose Cartan projection is equal to $\type$. This vector is equal to $\mathrm{v}_{a,u(y)}$ for some $a\in \mathcal{F}_\type\cup \mathcal{F}_{\iota(\type)}$.

\medskip

Since $\mathrm{v}_{a,u(y)}$ and $\mathrm{d}u(\mathrm{v})$ are in a common flat, $\Hess_{b_a}\left(\mathrm{d}u(\mathrm{v}),\mathrm{d}u(\mathrm{v})\right)=0$. One can assume that $\langle \II_u(\mathrm{v},\mathrm{v}),\mathrm{v}_{a,u(y)}\rangle_{u(y)}\leq 0$ up to exchanging $a$ with its symmetric with respect to $u(y)$ which is still in $\mathcal{F}_\type\cup \mathcal{F}_{\iota(\type)}$. Moreover since $\langle\mathrm{v}_{a,u(y)},\mathrm{d}u(\mathrm{v})\rangle_{u(y)}=0$, this is a contradiction with the criterion from Proposition \ref{prop:nearlyfuchsianLocalProp}, so the immersion $u$ cannot be $\type$-nearly geodesic.\end{proof}

The property of being $\type$-nearly geodesic is not necessarily satisfied for totally geodesic immersions, but it is satisfied for $\type$-regular totally geodesic immersions.

\begin{prop}
\label{prop:TotGeodNearlyGeod}
A totally geodesic immersion is $\type$-nearly geodesic if and only if it is $\type$-regular.
\end{prop}

\begin{proof}

An immersion $u$ is totally geodesic if and only if $\II_u= 0$. If $u$ is a $\type$-nearly geodesic immersion that is totally geodesic, for every $y\in M$, $\mathrm{v}\in T_yM$ and every $a\in \mathcal{F}_\type$, $y\in M$ and $\mathrm{v}\in T_yM$ Proposition \ref{prop:nearlyfuchsianLocalProp} implies that :
$$\Hess_{b_a}(\mathrm{d}u(\mathrm{v}),\mathrm{d}u(\mathrm{v}))>0.$$

The implies that for no $a\in \mathcal{F}_\type$ the vector $\mathrm{v}_{a,u(y)}$ lies in a common flat with  $\mathrm{d}u(\mathrm{v})$ by Lemma \ref{lem:HessianofBusemannFunctions}. Hence the Cartan projection of $\mathrm{d}u(\mathrm{v})$ is not orthogonal to $w\in \type$ for any $w\in W$.

\medskip

Conversely, if the totally geodesic immersion is $\type$-regular, $\Hess_{b_a}(\mathrm{d}u(\mathrm{v}),\mathrm{d}u(\mathrm{v}))$ is never equal to $0$ for any $y\in M$, $\mathrm{v}\in T_yM$ and $a\in \mathcal{F}_\type$ such that $\langle \mathrm{v}_{a,u(y)},\mathrm{v}\rangle_y=0$. Since $\Hess_{b_a}$ is non-negative, Proposition \ref{prop:nearlyfuchsianLocalProp} implies that $u$ is $\type$-nearly geodesic.
\end{proof}

\begin{proof}
Consider $y_0\in M$. The function $y\in M\mapsto \exp\left(\lambda d^\type_\X(u(y),u(y_0))\right)$ is strictly convex for some $\lambda>0$ and admits a minimum at $y=y_0$. The completeness of the metric $u^*(g_\X)$ implies that there is a geodesic joining any two points. Hence the minimum of any strictly convex function is unique, so $u$ is injective: it is an embedding.
\end{proof}

\subsection{Uniformly nearly geodesic immersions.}
If the nearly geodesic condition for an immersion is satisfied uniformly, one can prove that the exponential of some multiple of Busemann functions are strictly convex on the image of the immersion. 
%

\begin{defn}
\label{defn:UniformlyNearlyGeod}
Let $\type\in \Sph\mathfrak{a}^+$. An immersion $u:M\to \X$ is \emph{uniformly} $\type$-nearly geodesic if there exist $\epsilon>0$ such that for all $\mathrm{v}\in TM$ such that $\lVert\mathrm{d}u(\mathrm{v})\rVert=1$ one has for all $a\in \mathcal{F}_\type$ satisfying $\mathrm{v}_{a,o}\perp\mathrm{v}$:
$$\Hess_{b_a}(\mathrm{d}u(\mathrm{v}),\mathrm{d}u(\mathrm{v}))+\langle\II_u(\mathrm{v},\mathrm{v}),\mathrm{v}_{a,o}\rangle_o\geq\epsilon.$$
\end{defn}

\begin{rem}
When $\X=\mathbb{H}^n$ being uniformly nearly geodesic is equivalent to having principal curvature in $(-\lambda,\lambda)$ for some $\lambda<1$.
\end{rem}

Suppose that $M=\widetilde{N}$ is the universal cover of a compact smooth manifold $N$.  Let $\Gamma$ be the fundamental group of $N$. A $\rho$-equivariant immersion $u:M\to \X$ for some representation $\rho:\Gamma\to G$ which is $\type$-nearly geodesic is necessarily uniformly $\type$-nearly geodesic since $T^1N$ is compact.

\medskip

If we consider a uniformly $\type$-nearly geodesic immersion $u$, not only are Busemann functions convex in critical directions, but for some $\lambda>0$, $e^{\lambda b_{a,o}\circ u}$ is strictly convex on $M$.

\begin{lem}
\label{lem:nearlyGeodConvex}
Let $\type\in \Sph\mathfrak{a}^+$. Let $u:M\to \X$ be a uniformly $\type$-nearly geodesic immersion. For some $\lambda>0$, for all $a\in\mathcal{F}_\type$ the function $\exp\left(\lambda b_{a,o}\circ u\right)$ has positive Hessian for the metric $u^*(g_\X)$.
Moreover there exists some $\epsilon>0$ such that for any $a\in \mathcal{F}_\type$ and any geodesic $\eta:\R\to M$ the functions $f_\eta=\exp\left(\lambda b_{a,o}\circ u\circ \eta\right)$ satisfy $f''\geq\epsilon f$.
\end{lem}

Recall that the metric on $M$ that we consider to define geodesics is the induced metric $u^*(g_\X)$.

\medskip

\begin{proof}

Let $o\in \X$ and let $U_\type^\epsilon$ be the compact set of pairs $(\mathrm{v},\mathrm\II)\in T^1\X_o\times T_o\X$ such that for all $a\in \mathcal{F}_\type$ satisfying $\mathrm{v}_{a,o}\perp\mathrm{v}$:
$$\Hess_{b_a}(\mathrm{v},\mathrm{v})+\langle\II,\mathrm{v}_{a,o}\rangle_o\geq\epsilon.$$

Let us consider :
$$C=\inf_{a\in \mathcal{F}_\type,(\mathrm{v},\II)\in U_\type^\epsilon}\frac{\Hess_{b_a}(\mathrm{v},\mathrm{v})+\langle\II,\mathrm{v}_{a,x}\rangle_o}{\langle \mathrm{v}_{a,o},\mathrm{v}\rangle_o^2}.$$

This infimum is the infimum of a continuous function taking values in $\R\cup \lbrace +\infty\rbrace $ on a compact set. Indeed the numerator must be strictly positive whenever the denominator vanishes, and the denominator is always positive. Hence $C\in \R\cup\lbrace +\infty\rbrace $. 

\medskip

Let $\lambda$ be any real number greater than $\max(1-C,0)$. Let $\eta$ be any geodesic in $M$. Let us write $g=b_{a,o}\circ u\circ \eta$. Note that $g''\geq C(g')^2$ by definition of $C$. Therefore:
$$\left(e^{\lambda g}\right)''/e^{\lambda g}=\lambda g''+\lambda^2 \left(g'\right)^2\geq  (C\lambda+\lambda^2)\left(g'\right)^2+(\lambda-C)g''\geq\lambda \left(g'\right)^2\geq 0.$$

Note also that $\left(e^{\lambda g}\right)''/e^{\lambda g}\geq \lambda g''$. Consider the following quantity:

 $$M=\inf_{a\in \mathcal{F}_\type,(\mathrm{v},\II)\in U_\type^\epsilon} \max\left(\Hess_{b_a}(\mathrm{v},\mathrm{v})+\langle\II,\mathrm{v}_{a,x}\rangle_o,\langle \mathrm{v}_{a,o},\mathrm{v}\rangle_o^2\right).$$ 
 
Note that $M\leq \max\left(g'',\left(g'\right)^2\right)$. Since $K\subset U_\type$, this quantity is strictly positive as it is an infimum taken on a compact set of a positive function. Hence the function $f=e^{\lambda g}$ is strictly convex and satisfies $f''>\lambda M f$.
\end{proof}

\subsection{Convexity of a Finsler distance.} 
\label{subsec:QIEmbedding}
When $\X=\mathbb{H}^n$, and given $y\in \mathbb{H}^n$, for any nearly geodesic immersion $u:M\to \mathbb{H}^n$ the function $x\mapsto \exp\left(d_{\mathbb{H}^n}(u(x),y)\right)$ is strictly convex. However for a general symmetric space of higher rank the $\type$-nearly geodesic condition doesn't imply the convexity for the Riemannian metric at critical points.

 This leads us to consider a Finsler pseudo distance $d^\type_\X$ on $\X$ associated to an element $\type\in \Sph\mathfrak{a}^+$. We show in this section that this pseudo distance satisfies a similar convexity property for any $\type$-nearly geodesic immersion. This pseudo-distance is symmetric when $\type$ is symmetric and it is equal to the Riemannian distance when $\rank(\X)=1$. The convexity of this distance allows us to prove the injectivity and properness of complete $\type$-nearly geodesic immersions. This Finsler pseudo distance is studied in \cite[{Section 5}]{KLBFinsler}.

\medskip

Let us define for $\typeS\in \mathfrak{a}$:
$$|\typeS|_\type=\max_{w\in W}\langle w\cdot \type,\typeS\rangle.$$

The map $\typeS\mapsto |\typeS|_\type$ is non-negative, homogeneous and subadditive, thus we call it in general a pseudo-norm.

\medskip

This pseudo-norm is not necessarily symmetric: $|\typeS|_\type=|-\typeS|_{\iota(\type)}$. In particular it is symmetric if and only if $\type$ is symmetric. Figure \ref{fig:Hexagon} illustrates the unit ball of this norm in $\mathfrak{a}$ for two examples of semi-simple Lie groups whose associated symmetric space has rank $2$: on the left $G=\SL(3,\R)$ and $\type=\type_\Delta$, in the middle $G=\SL(3,\R)$ and $\type=\type_1$ (such that $\mathcal{F}_{\tau_1}\simeq \mathbb{RP}^2$) and on the right $G=\Sp(4,\R)$ and $\type=\type_{\lbrace\alpha_n\rbrace}$ with the notations from Section \ref{subsec:ExamplesRoots}.

\medskip
\begin{figure}
\begin{center}
\begin{tikzpicture}[scale=1]
\draw[dashed] (-1,-0.866*2) -- (1,0.866*2);
\draw[dashed] (-1,0.866*2) -- (1,-0.866*2);
\draw[dashed] (-2,0) -- (2,0);
\draw[very thick, red,->] (0,0) -- (0.866,1/2);
\fill (0,0) circle (.05);
\node[above right] () at (0.866,1/2) {$\type_\Delta$};
\fill[thick, blue,opacity=0.3] (0.577,1) -- (2*0.577,0) --  (0.577,-1) --  (-0.577,-1) -- (-2*0.577,0) -- (-0.577,1) -- cycle;
\draw[thick, blue] (0.577,1) -- (2*0.577,0) --  (0.577,-1) --  (-0.577,-1) -- (-2*0.577,0) -- (-0.577,1) -- cycle;
\end{tikzpicture} \hspace{2pt}
\begin{tikzpicture}[scale=1]
\draw[dashed] (-1,-0.866*2) -- (1,0.866*2);
\draw[dashed] (-1,0.866*2) -- (1,-0.866*2);
\draw[dashed] (-2,0) -- (2,0);
\draw[very thick, red,->] (0,0) -- (1/2,0.866);
\fill (0,0) circle (.05);
\node[right] () at (.54,0.87) {$\type_1$};
\fill[thick, blue,opacity=0.3] (-1,0.866*2) -- (-1,-0.866*2) --  (2,0) -- cycle;
\draw[thick, blue] (-1,0.866*2) -- (-1,-0.866*2) --  (2,0) -- cycle;
\end{tikzpicture} \hspace{2pt}
\begin{tikzpicture}[scale=1]
\draw[dashed] (0,-0.866*2) -- (0,0.866*2);
\draw[dashed] (-0.866*2,0) -- (0.866*2,0);
\draw[dashed] (1.414,1.414) -- (-1.414,-1.414);
\draw[dashed] (1.414,-1.414) -- (-1.414,1.414);
\draw[very thick, red,->] (0,0) -- (1,0);
\fill (0,0) circle (.05);
\node[above right] () at (1,0) {$\type_{\lbrace{\alpha_2\rbrace}}$};
\fill[thick, blue,opacity=0.3] (1,1) -- (1,-1) --  (-1,-1) --  (-1,1) -- cycle;
\draw[thick, blue] (1,1) -- (1,-1) --  (-1,-1) --  (-1,1) -- cycle;
\end{tikzpicture}
\end{center}

\caption{The unit ball in $\mathfrak{a}$ of $|\cdot|_{\type_\Delta}$ and $|\cdot|_{\type_1}$ for $G=\SL(3,\R)$ and $|\cdot|_{\type_{\lbrace{\alpha_2\rbrace}}}$  for $G=\Sp(4,\R)$.}
\label{fig:Hexagon}
\end{figure}

This pseudo-norm isn't necessarily positive on non-zero vectors. However if a non-zero vector $\mathrm{v}\in \mathfrak{a}$ has zero norm, the Weyl group does not act irreducibly on $\mathfrak{a}$ since the $W$-orbit of $\type$ is orthogonal to $\mathrm{v}$, which means that the underlying Lie group $G$ is not simple.

\begin{exmp}
Let $n\geq 2$ be an integer and $G=\PSL(2,\R)^n$, $\X=\left(\mathbb{H}^2\right)^n$. A model flat in $\X$ is the product of a geodesic in each of the $n$ copies of $\mathbb{H}^2$. Let $\type_k$ be the tangent vector to the geodesic on the $k$-th copy of $\mathbb{H}^2$, the pseudo-distance on $\X$ defined by the pseudo-norm $|\cdot|_{\type_k}$ is the distance in $\mathbb{H}^2$ of the $k$-th components, which is not a distance on $\X$.
\end{exmp}

However if $G$ is simple and $\type$ is symmetric $|\cdot|_\type$ is a norm. This norm is $W$-invariant and hence it defines a $G$-invariant not necessarily symmetric Finsler metric on $\X$ such that for $\mathrm{v}\in T\X$, $\lVert\mathrm{v}\rVert_\type=|\Cartan(\mathrm{v})|_\type$ , where $\Cartan$ is the Cartan projection (\cite{Finsler} Theorem 6.2.1). 

\medskip

For any semi-simple Lie group $G$, we denote by $d^\type_\X:\X\times \X\to \R_{\geq 0}$ the corresponding pseudo-distance on $\X$, \ie $d_\X(x,y)$ is the infimum for all piece-wise $\mathcal{C}^1$-path $\eta$ from $x$ to $y$ of:
 $$\int \lVert\eta'\rVert_\type=\int |\Cartan(\eta')|_\type.$$
  This distance can be characterized in terms of Busemann functions.

\begin{prop}
\label{prop:defNorm}
Let $x,y\in \
X$ be two points. The pseudo distance $d^\type_\X$ between these two points satisfies:
$$d^\type_\X(x,y)=\max_{a\in \mathcal{F}_\type}b_{a,x}(y).$$
\end{prop}

\begin{proof}
Let $o\in \X$, $\mathrm{v}\in  T_o\X$ and $a\in \mathcal{F}_\type$. As usual $\mathrm{v}_{a,o}$ is the unit vector based at $o$ pointing towards $a$. The maximum for $a\in \mathcal{F}_\type$ of $\langle \mathrm{v},\mathrm{v}_{a,o}\rangle$ is reached when $\mathrm{v}$ and $\mathrm{v}_{a,o}$ are in a common flat (\cite{Eberlein} Proposition 24).

\medskip

If we assume that $\mathrm{v}$ and $\mathrm{v}_{a,o}$ are in a common flat, the maximum is equal to $|\Cartan(\mathrm{v})|_\type$. Given two points $x,y\in \X$, any piece-wise $\mathcal{C}^1$ curve $\eta$ such that $\eta(0)=x, \eta(1)=y$ satisfies  for all $a\in \mathcal{F}_\type$:
$$(b_a\circ \eta)'=\langle\mathrm{v}_{a,\eta(t)},\eta'(t)\rangle_{\eta(t)}\leq \lVert\eta'\rVert_\type.$$

Hence:
 $$b_{a,x}(y)\leq \int \lVert\eta'\rVert_\type.$$
 
Moreover equality is reached for the Riemannian geodesic such that $\eta(0)=x, \eta(1)=y$. Indeed there is a point $a\in \mathcal{F}_\type$ that lies in a common flat with $x$ and $y$ such that $|\eta'(t)|_\type=\langle \eta'(t),\mathrm{v}_{a,\eta(t)} \rangle_{\eta(t)}$ for all $t\in [0,1]$. Hence $b_{a,x}(y)=d^\type_\X(x,y)$. Note that the curve reaching this minimum is not unique in general. 
\end{proof}

This pseudo distance satisfies the desired convexity condition.
\begin{prop}
\label{prop:convexNorm}
Let $u:M\to \X$ be a uniformly $\type$-nearly geodesic immersion. There exist $\lambda>0$ such that for all $x\in \X$  the following function is strictly convex for the metric $u^*g_\X$:
$$f:y\in M\mapsto \exp\left(\lambda d^\type_\X(x,u(y))\right).$$

\end{prop}

A strictly convex continuous function on the Riemannian manifold $\left(M, u^*(g_\X)\right)$ is a function that is strictly convex on any geodesic.

\medskip

\begin{proof}
By the Lemma \ref{lem:nearlyGeodConvex} there exist $\lambda>0$ such that for any $a\in \mathcal{F}_\type$, the function $\exp\left(\lambda b_{a,o}\circ u\right)$ is strictly convex on $M$. One can then write $f$ as :

$$ f(y)=\exp\left(\lambda d^\type_\X(x,u(y))\right)=\sup_{a\in\mathcal{F}_\type} \exp\left(b_{a,x}\circ u (y)\right).$$

Hence $f$ is the supremum of a family of convex functions, so it is convex. Moreover the supremum is taken over a compact family of strictly convex functions, so it is strictly convex.

\end{proof}

A consequence of the convexity of this Finsler distance is that the immersion $u$ is injective, which is an interesting property of $\type$-nearly geodesic surface. We say that $u$ is \emph{complete} if $M$ is complete for the induced metric $u^*(g_\X)$.

\begin{prop}
\label{prop:nearlyGeodEmbedding}
Let $u:M\to \X$ be a complete uniformly $\type$-nearly geodesic immersion. Then $u$ is an embedding.
\end{prop}

\begin{proof}
Consider $y_0\in M$. The function $y\in M\mapsto \exp\left(\lambda d^\type_\X(u(y),u(y_0))\right)$ is strictly convex for some $\lambda>0$ and admits a minimum at $y=y_0$. The completeness of the metric $u^*(g_\X)$ implies that there is a geodesic joining any two points. Hence the minimum of any strictly convex function is unique, so $u$ is injective: it is an embedding.
\end{proof}

Moreover the immersion $u$ cannot be too distorded: the metric induced by $u$ is quasi-isometric to the ambient metric on $\X$. The notion of quasi-isometric embedding was recalled in Section \ref{sec:Representations}.

\begin{prop}
\label{prop:QuasiGeod}
Let  $u:M\to \X$ be a complete uniformly $\type$-nearly geodesic immersion. Then $u$ is a quasi-isometric embedding for the induced metric $u^*g_\X$ on $M$. In particular $u$ is proper.
\end{prop}

\begin{proof}

Let $y_0\in M$ and let $o=u(y_0)$. Let $\epsilon>0$ and $\lambda>0$ be the constants provided by Lemma \ref{lem:nearlyGeodConvex}. Let $\gamma:\R_{\geq 0}\to M$ be a geodesic ray parametrized with unit speed in $M$ for the metric $u^*g_\X$ with $\gamma(0)=y_0$.

\medskip

 Let $a\in \mathcal{F}_\type$ be such that $b_{a,o}\left(u\circ\gamma(1)\right)=d_\X^\type\left(o,u\circ\gamma(1)\right)$. Consider the function:
 $$f:t\in \R_{\geq 0}\mapsto \exp\left(\lambda b_{a,u\circ \gamma(t)}\left(u\circ \gamma(t)\right)\right).$$

It is strictly convex and satisfies $f(1)\geq f(0)$ so $f'(1)\geq 0$. Moreover $f''>\epsilon f$ so $f(t)\geq\cosh\left(\epsilon(t-1)\right)> \frac{e^{\epsilon(t-1)}}{2}$. In particular:
$$d_\X^\type\left(o,u\circ\gamma(t)\right)\geq b_{a,o}\left(u\circ\gamma(t)\right)\geq \frac{\epsilon}{\lambda}(t-1)-\frac{\log(2)}{\lambda}$$
For all $y\in M$ there exist a geodesic ray $\gamma$ passing through $y$. If $d_{u}$ is the Riemannian distance on $M$ induced by $u^*g_\X$:
$$d_\X^\type\left(u(x_0),u(x)\right)\geq \frac{\epsilon}{\lambda}(d_u\left(x,y\right)-1) \frac{\log(2)}{\lambda}.$$

This Finsler metric is equivalent to the Riemannian metric $g_\X$ if $G$ is simple, and in general it is dominated by the Riemannian metric. Moreover $u$ is $1$-Lipshitz with respect to the induced metric, so $u$ is a quasi-isometric embedding.

\end{proof}

Using the convexity of this Finsler pseudo-distance one can define a continuous projection from the whole symmetric $\X$ to $M$.  This projection will not be used in what follows, but the fibration of the domains in $\mathcal{F}_\type$ constructed in Section \ref{sec:Fibration} is an extension of it. 

\begin{prop}
\label{prop:ProjectionFinsler}
Let $u:M\to \X$ be a complete uniformly $\type$-nearly geodesic immersion. For every $x\in \X$, there exist a unique point $\pi^\type_u(x)\in M$ that minimizes:
$$y\in M\mapsto d^\type_\X(x,u(y)).$$

The function $\pi^\type_u:\X\to M$ is continuous, and $\pi^\type_u(u(y))=y$ for $y\in M$.
\end{prop}

\begin{proof}
Let $\lambda,\epsilon$ be the two constants provided by the Lemma \ref{lem:nearlyGeodConvex} and let $x\in\X$. The following function is strictly convex on $M$: 
$$y\mapsto\exp\left(\lambda d^\type_\X(x,u(y))\right).$$

It is moreover proper since $u$ is proper by Proposition \ref{prop:QuasiGeod}. Hence it has a unique minimum, so $\pi^\type_u$ is well defined. 

\medskip

If we consider a sequence $(x_n)\in \X$ of points that converge to $x\in \X$, then the sequence $(\pi^\type_u(x_n))$ is bounded since $\rho$ is discrete. Moreover any of its limit points is a minimum of $d^\type_\X(x,u(y))$, so the sequence converges to $\pi^\type_u(x)$. The function $\pi^\type_u$ is hence continuous.
\end{proof}
\subsection{Anosov property for nearly Fuchsian representations.}
\label{sec:AnosovNearlyFuch}
Let $N$ be a compact manifold with fundamental group $\Gamma$. 
We call a representation $\rho:\Gamma\to G$ that admits a $\type$-nearly geodesic equivariant immersion $u:\widetilde{N}\to \X$ a \emph{$\type$-nearly Fuchsian} representation. 

\begin{prop}
The set of $\type$-nearly Fuchsian representations is open in the space of representations $\rho:\Gamma\to G$, for the compact-open topology.
\end{prop}
\begin{proof}
One can continuously deform any $\rho$-equivariant immersion $u:\widetilde{N}\to \X$ to a $\rho'$-equivariant smooth map $u':\widetilde{N}\to \X$ for $\rho'$ close to $\rho$. Indeed fix a Riemannian metric on $N$ let $\eta:\R^+\to \R^+$ be a smooth function that is positive on $[0,R]$ for $R$ large enough and vanishes on $[R',+\infty)$ for some $R'>R$. One can define $u'(y)$ for $y\in \widetilde{N}$ as the barycenter of the points $x^y_\gamma=\rho'(\gamma^{-1})\cdot u(\gamma\cdot y)$ with weight $\lambda^y_\gamma=\eta(d(y,\gamma\cdot y))$ for $\gamma\in \Gamma$. Concretely this means that we consider the unique local minimum of the convex function :
$$D:x\in \X\mapsto \sum_{\gamma\in \Gamma}\lambda^y_\gamma d(x,x^y_\gamma)^2.$$
Note that $\rho'(\gamma_0)\cdot x_\gamma^y=x^{\gamma_0\cdot y}_{\gamma\gamma_0^{-1}}$ and $\lambda^y_\gamma=\lambda^{\gamma_0\cdot y}_{\gamma\gamma_0^{-1}}$ for $\gamma_0\in \Gamma$. Therefore $u'$ is $\rho'$-equivariant.
Since $\X$ is a Hadamard manifold $D$ is strictly convex so the barycenter map is well-defined and smooth. Therefore for $\rho'$ close enough to $\rho$, $u'$ is an immersion which is close to $u$ for the $\mathcal{C}^2$-topology on any compact fundamental domain of the action of $\Gamma$ on $\widetilde{N}$. In particular $u'$ is a $\type$-nearly geodesic immersion for $\rho$ close enough to $\rho'$.
\end{proof}

The condition that $u$ is $\type$-nearly geodesic is local, but it will imply some coarse property on $u$ and therefore on $\rho$. Recall that the limit cone $\mathcal{C}_\rho$ was defined in Section \ref{sec:Representations} (Definition \ref{defn:CartanCone}). Due to flats Busemann functions are not strictly convex in critical directions on $\X$. However Busemann functions are strictly convex in critical directions on $u(\widetilde{N})$. We deduce that $\type$-nearly geodesic surfaces must coarsely avoid these flats, which in turn can be interpreted as a property of the limit cone $\mathcal{C}_\rho$ (see Definition \ref{defn:CartanCone}.

\begin{prop}
\label{prop:CartanCone}
Let $\rho:\Gamma\to G$ be a $\type$-nearly Fuchsian representation :
\begin{equation}
\label{eq:CartanCone}
\mathcal{C}_\rho\cap \bigcup_{w\in W} (w\cdot \type)^\perp=\emptyset .\end{equation}
\end{prop}

Recall that $W$ is the Weyl group associated to $G$.

\medskip

\begin{proof}
Let $x_0,x\in \widetilde{N}$ and $o=u(x_0)$. Let $w\in W$, there exist two points $a\in \mathcal{F}_\type$  and $a'\in \mathcal{F}_{\iota(\type)}$ that are opposite from $o$, \ie $\mathrm{v}_{a,o}=-\mathrm{v}_{a',o}$, and such that :
$$b_{a,o}(u(x))=\langle\da(o,u(x)),w\cdot \type\rangle,$$
$$b_{a',o}(u(x))=\langle\da(o,u(x)),-w\cdot \type\rangle.$$

This holds for $a,a'$ that lie in a maximal flat containing $o$ and $u(x)$. 

\medskip

Let $\eta$ be a geodesic parametrized with unit length in $\widetilde{N}$ for $u^*(g_\X)$ such that $\gamma(0)=x_0$ and $\gamma(d_u(x_0,x))=x$. Let $\lambda,\epsilon>0$ be the constants given by Lemma \ref{lem:nearlyGeodConvex}. Consider the function:
$$f:t\mapsto \exp\left(\lambda b_{a,o}\circ u\circ \gamma(t)\right).$$

Since $\mathrm{v}_{a,o}=-\mathrm{v}_{a',o}$, up to exchanging $a$ and $a'$ we can assume that $f'(0)\geq 0$. By Lemma \ref{lem:nearlyGeodConvex}, one has $f''>\epsilon f$. Together with the fact that $f(0)=1$, this implies that for all $t\in [0,d_u(x_0,x)]$:
$$f(t)\geq \cosh(\epsilon t).$$

Hence $\langle\da(o,u(x)),w\cdot p\rangle\geq \frac{\epsilon}{\lambda} d_u(x_0,x) -\frac{\log(2)}{\lambda}$. Since $u$ is a quasi-isometric embedding, there exist $c,D>0$ such that for all $x\in \widetilde{N}$ the distance for the induced metric $u^*(g_\X)$ between $x$ and $x_0$ is at least:
 $$c d_\X(o, u(x))-D.$$

 In conclusion:
$$\langle\frac{\da(o,u(x))}{d_\X(o,u(x))},w\cdot p\rangle\geq\frac{c\epsilon}{\lambda}- \frac{\log(2)+\epsilon D}{\lambda d_\X(o,u(x))}.$$

Any element of $\mathcal{C}_\rho$ has therefore a scalar product at least $\frac{c\epsilon}{\lambda}>0$ with $w\cdot \type$, for any $w\in W$. This implies that the limit cone cannot intersect $\bigcup_{w\in W} (w\cdot \type)^\perp$.
\end{proof}

The set $\Sph\mathfrak{a}^+\setminus \bigcup_{w\in W} (w\cdot \type)^\perp$ contains a single connected component if and only if $(w\cdot \type)^\perp$ is always a wall of the Weyl chamber decomposition of $\mathfrak{a}$, \ie when $\type=\type_\Theta$ for a Weyl orbit of simple roots  $\Theta\subset \Delta$ (this notion was defined Section \ref{sec:RootSpecial}). In this case:
$$\Sph\mathfrak{a}^+\setminus \bigcup_{w\in W} (w\cdot \type_\Theta)^\perp=\Sph\mathfrak{a}^+\setminus \bigcup_{\alpha\in \Theta} \Ker(\alpha).$$

Hence we get the following.

\begin{thm}
\label{thm:NearlyfuchsianAnosovNormal}
Let $\Theta\subset \Delta$ be a Weyl orbit of simple roots. A $\type_\Theta$-nearly Fuchsian representation $\rho:\Gamma_g\to G$ is $\Theta$-Anosov.
\end{thm}
 In particular only hyperbolic groups admit $\type_\Theta$-nearly Fuchsian representations, (see \cite[Theorem 3.2]{BPS}).
 
 \medskip

If $\type\in \Sph\mathfrak{a}^+$ does not correspond to a Weyl orbit of simple roots, let us assume that $\Gamma$ is a non-elmentary Gromov hyperbolic group, so that the limit cone of $\Gamma$ is connected, (see Lemma \ref{lem:CartanConeConnected}).

\medskip  

To a $\type$-nearly Fuchsian representation $\rho:\Gamma\to G$ one can associate the connected component $\sigma_\rho^\type$ in which $\mathcal{C}_\rho$ lies inside:
 $$\Sph\mathfrak{a}^+\setminus \bigcup_{w\in W} (w\cdot \type)^\perp.$$
 
To a connected component of this space one can associate a  non-empty set $\Theta(\sigma_\rho^\type)$ of simple roots. Recall that for $\typeS\in \Sph\mathfrak{a}^+$, $\Theta(\typeS)\subset \Delta$ is the set of simple roots $\alpha$ such that $\alpha(\typeS)\neq 0$.

\begin{lem}
\label{lem:ConnectedRoot}
Let $\type\in \Sph\mathfrak{a}^+$ and let $\sigma\subset \Sph\mathfrak{a}^+$ be a connected component of :
\begin{equation}
\label{eq:WeylWalls}
\Sph\mathfrak{a}^+\setminus \bigcup_{w\in W} (w\cdot \type)^\perp.\end{equation}

Let $\Theta(\sigma)\subset \Delta$ be the set of simple roots $\alpha$ such that $\sigma\cap \Ker(\alpha)=\emptyset$. This set is non-empty, and there exist some $\typeS\in \sigma$ such that $\Theta(\typeS)=\Theta(\sigma)$.
\end{lem}

\begin{figure}
\begin{center}
\begin{tikzpicture}[scale=3]

\fill[opacity=0.1] (-1,0) -- (1,0) -- (0,1) -- cycle ;
\fill[red,opacity=0.2] (0,1) -- (-1/2,1/2) -- (1/2,1/2) -- cycle ;
\draw (-1,0) -- (1,0) -- (0,1) -- cycle ;

\draw[very thick, blue] (-1/2,1/2) -- (1/2,1/2) -- (0,0) -- cycle ;
\node[below] () at (0,1/3) {$\type$};
\node () at (0,4/6) {$\sigma$};
\node[above] () at (0,1) {$\typeS$};
\fill[red] (0,1) circle(.03);
\fill[blue] (0,1/3) circle(.03);
\end{tikzpicture}
\end{center}
\caption{Illustration for $G=\PSL(4,\R)$ of a connected component $\sigma$ of $\Sph\mathfrak{a}^+\setminus \bigcup_{w\in W} (w\cdot \type)^\perp$ in an affine chart.}
\label{fig:WeylSL4Cut}
\end{figure}
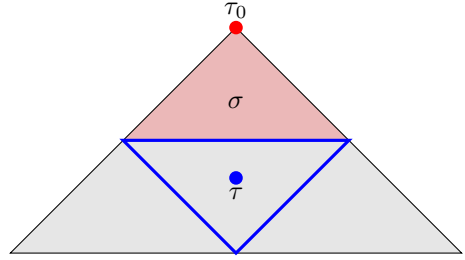

In other words there is $\typeS\in \sigma$ such that for any simple root $\alpha\in \Delta$, $\alpha(\typeS)\neq 0$ if and only if for all $\typeS'\in \sigma$, $\alpha(\typeS')\neq 0$. Figure \ref{fig:WeylSL4Cut} illustrates the lines in $\Sph\mathfrak{a}^+$ corresponding to $\bigcup_{w\in W}(w\cdot \type)^\perp$ for some $\type\in\Sph\mathfrak{a}^+$, as well as some connected component of the complement $\sigma$. In this example $\Theta(\sigma)$ contains only one root.

\medskip

\begin{proof}
Let $W_0\subset W$ be the subgroup of the Weyl group generated by symmetries associated to $\alpha\in \Delta\setminus \Theta(\sigma)$. Let $\hat{\sigma}\subset \Sph\mathfrak{a}$ be the connected component of $\sigma$ in: $$\Sph\mathfrak{a}\setminus \bigcup_{w\in W} (w\cdot \type)^\perp.$$
This connected component $\hat{\sigma}$ is stabilized by $W_0$. Indeed let $\alpha$ be in $\Delta\setminus \Theta(\sigma)$. By definition there is some $\mathrm{v}\in \sigma$ such that $\alpha(\mathrm{v})=0$, and hence that is fixed by the symmetry associated to $\alpha$. Thus the connected component of $\mathrm{v}$ in $\Sph\mathfrak{a}$ is stabilized by this symmetry, and hence by the group $W_0$.

\medskip

 Let $\typeS\in \sigma$ be any element and let $\typeS'\in \Sph\mathfrak{a}^+$ be the element that up to the action of $W$ is positively colinear to:
$$\sum_{w\in W_0} w\cdot \typeS.$$
This sum does not vanish since $\langle \type,w\cdot\typeS\rangle$ has constant sign for $w\in W_0$.
This element $\typeS'$ is $W_0$-invariant, hence for all $\alpha \in \Delta\setminus \Theta(\sigma)$, $\alpha(\typeS')=0$.

\medskip

Moreover, since the Lie group considered $G$ is semi-simple, the action of $W$ has no global fixed point on $\Sph\mathfrak{a}$, and hence $W_0\neq W$, which proves that $ \Theta(\sigma)\neq \emptyset$.
\end{proof}

\begin{thm}
\label{thm:NearlyfuchsianAnosov}
A $\type$-nearly Fuchsian representation $\rho:\Gamma_g\to G$ from a non-elementary hyperbolic group $\Gamma$ is $\Theta(\sigma_\rho^\type)$-Anosov.
\end{thm}

Note that $\Theta(\sigma_\rho^\type)\neq \emptyset$, because of Lemma \ref{lem:ConnectedRoot}.

\medskip

\begin{proof}
We use the characterization of Anosov representations from Theorem \ref{thm:CaracAnosocQuasiIsom}. We already proved that $\type$-nearly Fuchsian representations are quasi-isometric embeddings in Proposition \ref{prop:QuasiGeod}. 

\medskip

Moreover the limit cone $\mathcal{C}_\rho$ lies inside $\sigma^\type_\rho$, which avoids $\Ker(\alpha)$ for $\alpha\in\Theta(\sigma_\rho^\type)$. Hence $\rho$ is $\Theta(\sigma_\rho^\type)$-Anosov.
\end{proof}

The non-elementary assumption is necessary: indeed the following representation $\rho:\mathbb{Z}\to \SL(3,\R)$ is not Anosov for any set of roots:
$$n \mapsto \begin{pmatrix}
4^n &0 &0 \\
 0& 2^{-n}&0 \\
 0&0 &2^{-n}
\end{pmatrix}.$$

However this representation preserves a geodesic which is $\type$-regular for almost every $\type\in \Sph\mathfrak{a}^+$.

\subsection{A sufficient condition for an immersion to be nearly geodesic.}
\label{subsec:Finsler}
Let $\Theta$ be a Weyl orbit of simple roots as in Section \ref{sec:RootSpecial}. Let $\alpha\in \Theta$ be any root. We define the following constant:
\begin{equation}
\label{eq:constantTheta}
c_\Theta= \min_{\beta\in \Sigma, \beta(\type_\Theta)\neq 0}\frac{|\beta(\type_\Theta)|}{\lVert\alpha\rVert^2}.
\end{equation}

Here $\lVert\alpha\rVert$ for $\alpha \in \Theta$ denotes the maximum of $|\alpha(\type)|$ for $\type\in \Sph\mathfrak{a}$ a unit vector. This quantity is the same for any $\alpha\in \Theta$, since $\Theta$ is a Weyl orbit of simple roots.

\medskip

A sufficient condition for the immersion $u$ to be a $\type_\Theta$-nearly geodesic surface is the following.

\begin{thm}
\label{thm:Sufficient condition}

Let $u:S\to \X$ be an immersion that satisfies for all $\mathrm{v}\in TS$ and $\alpha\in \Theta$:
\begin{equation}\label{eq:Almostfuchsian}\lVert\II_u(\mathrm{v},\mathrm{v}) \rVert_{\type_\Theta}<c_\Theta \alpha\left(\Cartan\left(\mathrm{d}u(\mathrm{v})\right)\right)^2.
\end{equation}
Then $u$ is a $\type_\Theta$-nearly geodesic immersion.
\end{thm}

Note that $\lVert\cdot \rVert_{\type_\Theta}<\lVert\cdot\rVert$ so having Inequality \eqref{eq:Almostfuchsian} with the Riemannian metric in the left hand side instead of the Finsler pseudo-distance from Section \ref{subsec:Finsler} is also a sufficient condition.

\medskip

This property is a generalization of the property of having principal curvature in $(-1,1)$, where the norm of the tangent vector is replaced by the evaluation of roots of the Cartan projection.

\medskip

\begin{proof}

Let us show that \eqref{eq:Almostfuchsian} implies the condition of Proposition \ref{prop:nearlyfuchsianLocalProp}: 
$$\Hess_{b_a}(\mathrm{d}u(\mathrm{v}),\mathrm{d}u(\mathrm{v}))+\langle\II_u(\mathrm{v},\mathrm{v}),\mathrm{v}_{a,u(y)}\rangle_{u(y)}>0. $$

\medskip

Let $x=u(y)$. Let us write $\mathrm{d}u(\mathrm{v})=\mathrm{w}_0+\mathrm{w}^\perp$ where $\mathrm{w}_0\in \mathfrak{z}(\mathrm{v}_{a,x})\cap \mathbf{p}_x$ and $\mathrm{w}^\perp\in (\mathfrak{z}(\mathrm{v}_{a,x})\cap \mathbf{p}_x)^\perp$. Because of Lemma \ref{lem:HessianofBusemannFunctions}, one has :

\begin{equation}
\label{eq:HessianPositiveEq1}
\Hess_{b_{a}}\left(\mathrm{d}u(\mathrm{v}),\mathrm{d}u(\mathrm{v})\right)\geq\lVert\mathrm{w}^\perp\rVert^2 \min_{\beta\in \Sigma, \beta(\type_\Theta)\neq 0}|\beta(\type_\Theta)|.
\end{equation}

Lemma \ref{lem:GeneralizedDistanceLip} implies that for any $\alpha\in \Sigma$:
$$\alpha(\Cartan(\mathrm{w}_0))+\lVert\alpha\rVert\times \lVert\mathrm{w}^\perp\rVert\geq \alpha(\Cartan(\mathrm{d}u(\mathrm{v}))).$$ 

\medskip

We assumed that $\langle \mathrm{v}_{a,x},\mathrm{d}u(\mathrm{v})\rangle_x=0$. Since $\mathrm{v}_{a,x}\in \mathfrak{z}(\mathrm{v}_{a,x})\cap \mathbf{p}_x$, one has therefore $\langle \mathrm{v}_{a,x},\mathrm{w}^\perp\rangle_x=0$ and hence $\langle \mathrm{v}_{a,x},\mathrm{w}_0\rangle_x=0$. Moreover $\mathrm{v}_{a,x}$ and $\mathrm{w}_0$ are in a common flat. Since $\Cartan(\mathrm{v}_{a,x})=\type_\Theta$, this implies that $\alpha(\Cartan(\mathrm{w}_0))=0$ for some root $\alpha\in \Theta$. Therefore for this root $\alpha$:
\begin{equation}
\label{eq:HessianPositiveEq2}
\lVert\mathrm{w}^\perp\rVert\geq \lVert\alpha\rVert^{-1}\alpha(\Cartan(\mathrm{d}u(\mathrm{v}))).
\end{equation}

Recall that for any $\mathrm{w}\in T_x\X$ and $a\in \mathcal{F}_{\type_\Theta}$,  $\langle \mathrm{w},\mathrm{v}_{a,x}\rangle_x\leq \lVert\mathrm{w}\rVert_{{\type_\Theta}}$, as a consequence of \cite[Proposition 24]{Eberlein}.

Equation \eqref{eq:HessianPositiveEq1} and \eqref{eq:HessianPositiveEq2} imply together the following inequality, with $c_\Theta$ defined in \eqref{eq:constantTheta}:
$$\Hess_{b_a}\left(\mathrm{d}u(\mathrm{v}),\mathrm{d}u(\mathrm{v})\right)+\langle\II_u(\mathrm{v},\mathrm{v}),\mathrm{v}_{a,u(y)}\rangle_{u(y)}\geq c_\Theta\alpha(\Cartan(\mathrm{d}u(\mathrm{v})))^2-\lVert\II_u(\mathrm{v},\mathrm{v})\rVert_{\type_\Theta}.$$

The rightmost term is strictly positive because of the condition \eqref{eq:Almostfuchsian}. This concludes the proof.
\end{proof}

\begin{exmp}
Let $G=\PSL(n,\R)$. We chose the standard metric on $\X$ that comes from the Killing form. In particular the Euclidean metric on $\mathfrak{a}$ is given by:
$$\langle\Diag(\lambda_1,\cdots,\lambda_n),\Diag(\mu_1,\cdots,\mu_n)\rangle=2n\sum_{i=1}^n \lambda_i\mu_i.$$

In this case $\type_\Delta=\Diag(\frac{1}{2\sqrt{n}},0,\cdots,0,-\frac{1}{2\sqrt{n}})$. 
The minimum non-zero value of $\beta(\type_\Theta)$ for $\beta\in \Sigma$ is reached for the root $\alpha_1:\Diag(\lambda_1,\cdots,\lambda_n)\mapsto \lambda_1-\lambda_2$, and is equal to $\frac{1}{2\sqrt{n}}$ if $n\geq 3$ and $\frac{1}{\sqrt{n}}$ if $n=2$.

\medskip

The norm of any root $\alpha$ is equal to the norm of $\alpha_1$. But $|\alpha_1(\type)|\leq |\lambda_1|+|\lambda_2|\leq {\sqrt{2}}\sqrt{\lambda_1^2+\lambda_2^2}\leq \frac{1}{\sqrt{n}}\lVert \type\rVert$ with equality for some $\type\in \Sph \mathfrak{a}$. Hence $\lVert\alpha_1\rVert=\frac{1}{\sqrt{n}}$, so if $n\geq 3$:
$$c_\Delta=2\sqrt{n}.$$ 

And $c_\Delta=\sqrt{2}$ if $n=2$. Note that if we rescale the metric on $\X=\mathbb{H}^2$ so that the sectional curvature is equal to $-1$, Equation \eqref{eq:Almostfuchsian} is exactly the condition of having principal curvature in $(-1,1)$.
\end{exmp}

\section{Pencils of tangent vectors.}
\label{sec:Pencils}
In this section we recall the classical notion of a pencil of quadrics, then we generalize it to the notion of a pencil of tangent vectors in a symmetric space of non compact type and its base in a flag manifold. Bases of pencils appear as the fibers of the fibration that will be constructed in Section \ref{sec:Fibration}.

\subsection{Pencils of quadrics.}

Some references for the notion of pencil of quadrics can be found in \cite{fevola2020pencils}. Let $V$ be a finite dimensional vector space over $\K=\R$ or $\C$.

\begin{defn}
A \emph{pencil of quadrics}, or more precisely a \emph{$d$-pencil of quadrics}\footnote{In the literature, for instance in \cite{fevola2020pencils}, a \emph{pencil of quadrics} is often just a $2$-pencil of quadrics.} on $V$ is a linear subspace $\mathcal{P}$ of dimension $d$ in the space $\mathcal{S}(V)$ of symmetric bilinear forms on $V$ if $\K=\R$, or in the space $\mathcal{H}(V)$ of Hermitian forms on $V$ if $\K=\C$.

\medskip

The \emph{base} $b(\mathcal{P})$ of a $d$-pencil $\mathcal{P}$ is the set of points $[v]\in\mathbb{P}(V)$ such that for all $q\in \mathcal{P}$, $q(v,v)=0$.

\end{defn}

The following is a criterion for a pencil of quadrics to have a smooth base.

\begin{lem}
\label{lem:SubmertionPencils}
Let $\mathcal{P}$ be a pencil of quadrics such that all non-zero $q\in \mathcal{P}$ are non-degenerate bilinear forms. The map $p:V\to \mathcal{P}^*$, $v\mapsto \left(q\mapsto q(v,v)\right)$ is a submersion at every $v\in V\setminus \lbrace0\rbrace$ such that $[v]\in b(\mathcal{P})$. In particular $b(\mathcal{P})$ is a smooth manifold of codimension $d$.

\end{lem}
\begin{proof}
Let $(q_1,\cdots, q_d)$ be a basis of $\mathcal{P}$. Let us consider some $v\in V\setminus \lbrace0\rbrace$ such that $q_1(v,v)=\cdots=q_d(v,v)=0$. The kernel of the differential of $p$ is the intersection of the orthogonal spaces $[v]^{\perp_{q_i}}$ with respect to $q_i$ of the line generated by $v$ for $1\leq i\leq d$. Since the forms $q_i$ are non-degenerate, these are hyperplanes. 

Suppose that their intersection has not codimension $d$. In particular the linear forms $q_i(v,\cdot)$ for $1\leq i\leq d$ are not linearly independent, so there exist a linear combination of the bilinear forms that is degenerate, but is a non-zero element of $\mathcal{P}$, contadicting our assumption.

\smallskip

Hence the the kernel of $p$ has codimension $d$, so $p$ is a submersion at $v$.

\end{proof}

The base of a pencil of quadric is smooth and has codimension $d$ around each of it's points which are non-singular, meaning that they are not degenerate points for any quadric in the pencil. We generalize this notion of singular points in the next section.

\subsection{Pencils of tangent vectors in symmetric spaces.}

In this subsection we consider pencils of tangent vectors in a symmetric space $\X$ of non compact type, which are related to pencils of quadrics when $G=\PSL(n,\R)$.

\begin{defn}
A \emph{pencil of tangent vectors} at $x\in \X$, or more precisely a \emph{$d$-pencil}, is a vector subspace $\mathcal{P}\subset T_x\X$ of dimension $d$ for some point $x\in \X$.
\end{defn}

 To a pencil one can associate some subsets of any $G$-orbit in the visual boundary. Recall that for $a\in \partialvis \X$ and $x\in \X$, the unit vector $\mathrm{v}_{a,x}\in T_x\X$ is the unit vector pointing towards $a$. Let $\type\in \Sph \mathfrak{a}^+$.

\begin{defn}
The \emph{$\type$-base of the pencil $\mathcal{P}$}, whose base-point is $x\in \X$, is the set $\mathcal{B}_\type(\mathcal{P})$ of elements $a\in \mathcal{F}_\type$ such that $\mathrm{v}_{a,x}$ is orthogonal to $\mathcal{P}$.
\end{defn}

When $G=\PSL(n,\R)$ and $\X=\mathcal{S}_n$, a pencil at $q\in \mathcal{S}_n$ corresponds to a subspace $\mathcal{P}'$ of symmetric bilinear forms on $\R^n$, \ie a pencil of quadrics, that is compatible with $q$ in the sense that the trace of the associated $q$-symmetric matrices vanishes.

\begin{prop}
\label{prop:PencilQuadrics}
Let $\type\in \Sph \mathfrak{a}^+$ be such that $\mathcal{F}_\type\simeq \mathbb{RP}^{n-1}$. The $\type$-base of the pencil $\mathcal{P}$ is identified via this identification with the base of the pencil of quadrics $\mathcal{P}'$. 
\end{prop}

\begin{proof}
The $\type$-base of $\mathcal{P}$ is the space of lines $[C]$ where $C\in \R^n$ is a column vector satisfying $\tr(CC^{\perp_q}M)=0$ for all $M\in \mathcal{P}'$, since $CC^{\perp_q}$ is colinear to $\mathrm{v}_{[C],q}$. Hence the $\type$-base of the pencil is also the set of lines $[C]$ such that $C^{\perp_q}MC=0$, \ie the base of the pencil of quadrics $\mathcal{P}'$.

\end{proof}

We now generalize Lemma \ref{lem:SubmertionPencils} to general pencils of tangent vectors. 

\begin{defn}
A point $a\in \mathcal{F}_\type$ in the base $\mathcal{B}_\type(\mathcal{P})$ of a pencil $\mathcal{P}$ at $x\in X$ is called \emph{singular} if for some $\mathrm{w}\in \mathcal{P}$ one has $[\mathrm{w},\mathrm{v}_{a,x}]=0$.

\medskip

We denote by $\mathcal{B}^*_\type(\mathcal{P})\subset \mathcal{B}_\type(\mathcal{P})$ the set of non-singular points, that we will also call the \emph{regular base}.
\end{defn}

\begin{lem}
\label{lem:SmoothSymmetricPencils}
Let $\mathcal{P}$ be a pencil of tangent vectors at $x$ in $\X$. The function which associates to $a\in \mathcal{F}_\type$ the linear form $\mathrm{v} \mapsto \langle \mathrm{v}_{a,x},\mathrm{v} \rangle_x$ on $\mathcal{P}$ is a submersion at $a\in \mathcal{B}_\type(\mathcal{P})$ if and only if $a\in \mathcal{B}^*_\type(\mathcal{P})$. In particular $\mathcal{B}^*_\type(\mathcal{P})$ is always a smooth codimention $d$ submanifold of $\mathcal{F}_\type$.

\end{lem}

\begin{proof}
Let $\phi: \mathcal{F}_\type\to \mathcal{P}^*$ be the map that associates to $a\in \mathcal{F}_\type$ the linear form $\mathrm{v} \mapsto \langle \mathrm{v}_{a,x},\mathrm{v} \rangle_x$.

Suppose that $a\in \mathcal{B}_\type(\mathcal{P})\setminus \mathcal{B}^*_\type(\mathcal{P})$. Then there exist some $\mathrm{w}\in \mathcal{P}$ such that $[\mathrm{w},\mathrm{v}_{a,x}]=0$. The map $\psi:k\in K_x\mapsto k\cdot a\in \mathcal{F}_\type$ is a submersion, so for every tangent vector in $T_a\mathcal{F}_\type$ the differential of $a\mapsto\mathrm{v}_{a,x}$ in this direction is $\ad_\mathrm{k}(\mathrm{v}_{a,x})$ for some $\mathrm{k}\in \mathfrak{k}_x$. The differential of $a\mapsto \langle \mathrm{v}_{a,x},\mathrm{w}\rangle_x$ in this direction is equal to $\langle \ad_\mathrm{k}(\mathrm{v}_{a,x}),\mathrm{w}\rangle_x=-B(\ad_\mathrm{k}(\mathrm{v}_{a,x}),\mathrm{w})=-B(\mathrm{k},[\mathrm{v}_{a,x},\mathrm{w}])=0$. Hence the image of the differential of $\phi$ is not surjective: it is not a submersion.

\medskip

Suppose that $a\in \mathcal{B}^*_\type(\mathcal{P})$. Let $\mathrm{v}\in \mathcal{P}$ be any non-zero vector and consider $[\mathrm{v}_{a,x},\mathrm{v}]=\mathrm{k}\in \mathfrak{k}_x$. The differential of $a\mapsto \langle \mathrm{v}_{a,x},\mathrm{v}\rangle_x$ in the corresponding tangent direction is equal to $\langle \ad_\mathrm{k}(\mathrm{v}_{a,x}),\mathrm{v}\rangle_x=-B(\ad_\mathrm{k}(\mathrm{v}_{a,x}),\mathrm{v})=-B([\mathrm{v}_{a,x},\mathrm{v}],[\mathrm{v}_{a,x},\mathrm{v}])\neq 0$. Since for all $\mathrm{v}\in \mathcal{P}$ there is a direction in which the differential of $a\mapsto \langle \mathrm{v}_{a,x},\mathrm{v}\rangle_x$ does not vanish, the map $\phi$ is therefore a submersion at $a$.
\end{proof}

A pencil of tangent vectors $\mathcal{P}$ at $x\in \X$ is called \emph{$\type$-regular} if all its non-zero vectors are $\type$-regular as in Definition \ref{defn:pRegularVector}. In particular a $\type$-regular pencil satisfies $\mathcal{B}^*_\type(\mathcal{P})=\mathcal{B}_\type(\mathcal{P})$, so the $\type$-base of a $\type$-regular vector is a smooth codimension $d$ submanifold of $\mathcal{F}_\type$.

\medskip

Because of Lemma \ref{lem:SmoothSymmetricPencils}, the topology of the base of a  regular pencil is not varying if the pencil is deformed continuously.

\begin{cor}
\label{cor:PencilsDiffeo}
Let $\mathcal{P}_0$ and $\mathcal{P}_1$ be two pencils at $x\in \X$ in the same connected component of the space of $\type$-regular pencils at $x$. Then $\mathcal{B}_\type(\mathcal{P}_1)$ and $\mathcal{B}_\type(\mathcal{P}_2)$ are diffeomorphic.
\end{cor}

\begin{proof}
Since the space of regular pencils is open is the Grassmanian of planes in $T_x\X$, there exist a smooth path $(\mathcal{P}_t)_{t\in[0,1]}$ of regular pencils between $\mathcal{P}_0$ and $\mathcal{P}_1$. Because of Lemma \ref{lem:SmoothSymmetricPencils} the set $\lbrace(a,t)|a\in \mathcal{F}_\type, t\in [0,1]\rbrace$ is a submanifold with boundary of $\mathcal{F}_\type\times [0,1]$ that comes with a natural submersion $(a,t)\mapsto t$. Since this manifold is compact all the fibers are diffeomorphic by the Ehresmann fibration theorem.
\end{proof}
\begin{exmp}
\label{exmp:PencilsRegular}
Let $G=\PSL(3,\R)$, $\X=\mathcal{S}_3$. We identify the tangent space $T_{q_0}\mathcal{S}_3$ at the point $q_0$ corresponding to the standard scalar product on $\R^3$ with the space of $3$ by $3$ symmetric matrices with real coefficients  and zero trace. Consider the following two pencils:
$$\mathcal{P}_\text{irr}=\langle\begin{pmatrix}
0& 1 &0 \\
1&0  &1 \\
0&  1& 0
\end{pmatrix}, \begin{pmatrix}
2& 0 &0 \\
0&0  &0 \\
0&  0& -2\end{pmatrix}\rangle,\; \mathcal{P}_\text{red}=\langle\begin{pmatrix}
0& 0 &1 \\
0&0  &0 \\
1&  0& 0
\end{pmatrix}, \begin{pmatrix}
1& 0 & 0\\
0&0  &0 \\
0&  0& -1\end{pmatrix}\rangle.$$ 

Let $\type_1\in \Sph\mathfrak{a}^+$ be such that $\mathcal{F}_\mathrm{p_1}$ is diffeomorphic in a $\PSL(3,\R)$-equivariant way to $\mathbb{RP}^{n-1}$, and let $\type_\Delta$ be the normalized coroot associated to the Weyl orbit of simple roots  $\Delta$. It satisfies $\mathcal{F}_{\type_\Delta}\simeq \mathcal{F}_{1,2}$ the space of complete flags in $\R^3$.

The pencils $\mathcal{P}_\text{irr}$ and $\mathcal{P}_\text{red}$ are not $\type_1$-regular: $\mathcal{B}_{\type_1}(\mathcal{P}_\text{irr})$ is the disjoint union of a point and a line where $\mathcal{B}^*_{\type_1}(\mathcal{P}_\text{irr})$ contains only the point. In this particular case the regular base is a connected component of the base, so it is a smooth compact codimension 2 submanifold. 
The set $\mathcal{B}_{\type_1}(\mathcal{P}_\text{red})$ is a single point that is singular for the pencil. Here we see that a singular point can still be a point around which the base is a smooth codimension 2 submanifold. 

\smallskip

Both pencils are $\type_\Delta$-regular, but their $\type_\Delta$-bases are different.

\smallskip

A flag $(\ell,H)=([\mathrm{x}],[\mathrm{y}]^\perp)$ with $\mathrm{x}=(x_1,x_2,x_3)$ and $\mathrm{y}=(y_1,y_2,y_3)$ non-zero vectors such that $x_1^2+x_2^2+x_3^2=y_1^2+y_2^2+y_3^2$ and $x_1y_1+x_2y_2+x_3y_3=0$ belongs to $\mathcal{B}_{\type_\Delta}(\mathcal{P}_\text{red})$ if and only if:
$$x_1^2-x_3^2=y_1^2-y_3^2,$$
$$2x_1x_3=2y_1y_3.$$
Up to replacing $\mathrm{y}$ by $-\mathrm{y}$, these equations are equivalent to $x_1=y_1$, $x_2=-y_2$, $x_3=y_3$, $x_1^2+x_3^2=x_2^2$.
The corresponding flags $(\ell,H)$ in the affine chart $(x,y) \mapsto [x,1,y]$ of $\mathbb{RP}^2$ are the tangent point and tangent lines to the the circle of radius $1$ centered at the origin.

\smallskip

A flag $(\ell,H)=([\mathrm{x}],[\mathrm{y}]^\perp)$ with $\mathrm{x}=(x_1,x_2,x_3)$ and $\mathrm{y}=(y_1,y_2,y_3)$ non-zero vectors such that $x_1^2+x_2^2+x_3^2=y_1^2+y_2^2+y_3^2$ and $x_1y_1+x_2y_2+x_3y_3=0$ belongs to $\mathcal{B}_{\type_\Delta}(\mathcal{P}_\text{irr})$ if and only if:
$$2x_1^2-2x_3^2=2y_1^2-2y_3^2,$$
$$2x_1x_2+2x_2x_3=2y_1y_2+2y_2y_3.$$

Let $\ell_0=\langle(1,0,1)\rangle$ and $H_0=\langle(1,0,-1), (0,1,0) \rangle$. The corresponding flags $(\ell,H)$ belong to one of the three circles in $\mathcal{F}_{1,2}$ defined by :

\begin{itemize}
\item[-] $\ell=\ell_0$, $H$ any plane through $\ell$,
\item[-] $H=H_0$, $\ell $ any line in $H$,
\item[-] $\ell\subset H_0$, $\ell_0\subset H$.
\end{itemize}

Indeed one can check that these flags satisfy the equations. In order to check that these are the only solutions, one can see that these are the fibers of a fibration over the surface with $3$ connected components, see Section \ref{sec:TeichSL3}.

\smallskip

The $\type_\Delta$-base $\mathcal{B}_{\type_\Delta}(\mathcal{P}_\text{red})$ is a circle whereas $\mathcal{B}_{\type_\Delta}(\mathcal{P}_\text{irr})$ is the union of $3$ circles. Hence Corollary \ref{cor:PencilsDiffeo} implies that they must lie in different connected components of the space of $\type_\Delta$-regular pencils.
\end{exmp}

Since the pencils will be the fibers of the domains of discontinuity that we will construct, proving that the domain is non-empty will be equivalent to having non-empty pencils. We present here a topological argument to prove that some pencils are non-empty.

\begin{prop}
\label{prop:NonEmptyPencils}
Let $\type\in \Sph\mathfrak{a}^+$ and $\mathcal{P}$ be a $\type$-regular pencil of tangent vectors based at $x\in \X$ of dimension $d$. If the $\type$-base of $\mathcal{P}$ is empty, then $\mathcal{F}_\tau$ fibers over the sphere $S^{d-1}$.

 If moreover $d=2$ it implies that the fundamental group of $\mathcal{F}_\type$ is infinite.
\end{prop}

\begin{proof}
To $a\in \mathcal{F}_\type$ we associate $\pi_0(a)\in \mathcal{P}$ the orthogonal projection of $\mathrm{v}_{a,x}\in T_x\X$ onto $\mathcal{P}\subset T_x\X$. Since the $\type$-base of $\mathcal{P}$ is empty, one can define a map $\pi:\mathcal{F}_\type\to \Sph\mathcal{P}$ into the unit sphere of $\mathcal{P}$ where $\pi(a)=\frac{\pi_0(a)}{\lVert\pi_0(a)\rVert}$. This map is a submersion. Indeed let $a\in \mathcal{F}_\type$, let $\mathcal{P}_0\subset \mathcal{P}$ be the orthogonal to $\pi(a)$ in $\mathcal{P}$. Lemma \ref{lem:SmoothSymmetricPencils} applied to $\mathcal{P}_0$ implies that $\pi$ is a submersion at $a$.

\medskip

This submersion is proper since $\mathcal{F}_\type$ is compact, hence it is a fibration.
If $d=2$, this fibration induces a long exact sequence, where $F$ is the fiber.
$$\cdots \to \pi_1\left(\mathcal{F}_\type\right) \to \pi_1\left(\Sph\mathcal{P}\right)\to \pi_0\left(F\right)\to \cdots.$$

Since $F$ is compact,  $\pi_0\left(F\right)$ is finite and $\pi_1\left(\Sph\mathcal{P}\right)\simeq \mathbb{Z}$, so $\pi_1\left(\mathcal{F}_\type\right)$ is infinite.

\end{proof}

We conclude this section by the following remark that regular pencils cannot be tangent to flats.

\begin{prop}
\label{prop:Rank2}
If a $2$-pencil is tangent to a flat, then it is not $\type$-regular for any $\type\in \Sph\mathfrak{a}^+$. 
\end{prop}
\begin{proof}
Up to the action of $G$ one can identify $\mathcal{P}$ with a plane in $\mathfrak{a}$. But  for any $\type\in \Sph\mathfrak{a}^+$, the orthogonal of $\type$ intersects this plane. Hence there is an element of $\mathcal{P}$ whose Cartan projection is orthogonal to $w\cdot \type$ for some $w$ in the Weyl group.
\end{proof}

\section{Fibered domains in flag manifolds.}
\label{sec:Fibration}

In this section we associate an open domain $\Omega^\type_u\subset \mathcal{F}_\type$ to any complete uniformly $\type$-nearly geodesic immersion $u:M\to \X$ with $\type\in \Sph\mathfrak{a}^+$, and show that this domain is a smooth fiber bundle over $M$ where the fibers are $\type$-bases of the pencils that are the tangent planes to the surface. This is the construction is the analog of the \emph{Gauss map} for hypersurfaces in $\mathbb{H}^n$. We also mention what happens with our construction for totally geodesic immersions that are not $\type$-regular.

\medskip

If $M=\tilde{N}$ for some compact manifold $N$ with torsion-free fundamental group $\Gamma$, and if $u$ is equivariant with respect to a representation $\rho$, we show that the domain $\Omega^\type_u$ is a co-compact domain of discontinuity for the action of $\rho$ and its quotient fibers over $N$. This domain always coincides with some domain of discontinuity associated to Tits Bruhat ideals constructed by Kapovich-Leeb-Porti \cite{KLP-DoD}. Finally we prove the invariance of the topology of the quotients of these domains of discontinuity. 

\subsection{A domain associated to a nearly geodesic immersion.}

 Let $\type\in \Sph\mathfrak{a}^+$ be any unit vector and $u:M\to \X$ be a complete uniformly $\type$-nearly geodesic immersion.
 
\medskip

We consider a particular domain of the flag manifold $\mathcal{F}_\type$, defined for any nearly geodesic immersion $u:M\to \X$ using Busemann functions. For this we fix a base-point $o\in \X$, but the definition will dot depend on this choice.

\begin{defn}
\label{defn:DoDNearly}
Let $\Omega^\type_u$ be the set of elements  $a\in\mathcal{F}_\type$ such that the function $b_{a,o}\circ u$ is proper and bounded from below.
\end{defn}

We have additional properties if $u$ is a complete uniformly $\type$-nearly geodesic immersion.

\begin{lem}
\label{lem:FibrationCriticalPoints}
Let $a\in \mathcal{F}_\type$. There exist a critical point $x\in M$ for the function $b_{a,o}\circ u$ if and only if $a\in \Omega^\type_u$. In this case this point is unique, and the Hessian of $b_{a,o}\circ u$ at this point is positive. The domain  $\Omega^\type_u$ is open.
\end{lem}

\begin{proof}
Let $a\in \mathcal{F}_\type$. Suppose that $b_{a,o}\circ u$ is critical at $y\in M$. Since $u$ is $\type$-nearly geodesic the Hessian of $b_{a,o}\circ u$ at $y$ is positive. Moreover, due to Lemma \ref{lem:nearlyGeodConvex} there exist $\lambda>0$ such that $\exp\left(\lambda b_{a,o}\circ u\right)$ has positive Hessian everywhere on $M$. 

A convex function with positive Hessian on a complete connected Riemannian manifold has a unique minimum, and is proper. The function $\exp\left(\lambda b_{a,o}\circ u\right)$, and hence the function $b_{a,o}\circ u$ are hence proper and have a unique minimum.  In particular $a\in \Omega^\type_\rho$.

\medskip

Conversely if $a\in \Omega^\type_\rho$, $b_{a,o}\circ u$ is proper so it admits a global minimum, which is a critical point.

\medskip

If a function has a critical point with positive Hessian, every small deformation of the function for the $\mathcal{C}^0$-topology still admits a local minimum, and hence a critical point. Therefore $\Omega^\type_\rho$ is open.
\end{proof}

We thus can define the projection $\pi_u:\Omega^\type_u\to M$ associated to $u$ as the map that associates to  $a\in \Omega^\type_u$ the unique critical point $\pi_u(a)\in M$ of $b_{a,o}\circ u$. This is an extension at infinity of the nearest point projection from Proposition \ref{prop:ProjectionFinsler}.

\begin{thm}
\label{thm:FoliationDoD}

Let $u:M\to \X$ be a complete and uniformly $\type$-nearly geodesic immersion. The map $\pi_u:\Omega^\type_u\to M$ is a fibration. The fiber $\pi_u^{-1}(x)$ at a point $x\in M$ is the base $\mathcal{B}_\type(\mathcal{P}_x)$ of the $\type$-regular pencil $\mathcal{P}_x=\mathrm{d}u(T_xM)$.

\end{thm}

\begin{figure}
\begin{center}
\begin{tikzpicture}[scale=1.5]
\tdplotsetmaincoords{70}{90}
\tdplotsetrotatedcoords{0}{0}{0};
\draw[dashed, thick,
    tdplot_rotated_coords] (0,2,0) arc (90:270:2);

\fill[fill=blue, opacity=0.3,
    tdplot_rotated_coords
] (0,0,0) circle (2);
\tdplotsetmaincoords{0}{0}
\draw[dashed] (0,1.8,0) -- (0,-1.8,0);

\fill[red,opacity = 0.5] (0,-1.8,0) circle (.04);
\draw (0,-1.8,0) circle (.04);

\draw[fill=orange] (-0.2,0,-0.2) -- (.2,0,-.2) -- (.2,0,.2) -- (-.2,0,.2) -- cycle;
\fill (0,0,0) circle (.05);

\shade[ball color = lightgray,
   opacity = 0.6] (0,0,0) circle (2);
\fill[red] (0,1.8,0) circle (.04);
\draw[thick] (0,1.8,0) circle (.04);
\node[left] () at (0,1.8,0) {$a$};
\node[left] () at (0,-1.8,0) {$a'$};
\node[left] () at (-.6,1) {$\X=\mathbb{H}^3$};
\node[left] () at (2.4,1.8) {$\mathcal{F}_\type=\mathbb{CP}^1$};
\draw[thick,
    tdplot_rotated_coords] (0,2,0) arc (90:-90:2);

\node[right] () at (0.2,0,-0.2) {$\mathcal{P}$};
\node[above] () at (-0.3,0,-0.1) {$\pi_u(a)$};

\end{tikzpicture}
\end{center}
\caption{Fibration of the domain $\Omega^\type_u$ in the rank one case, $G=\PSL(2,\C)$.}
\label{fig:FibrationDoD}
\end{figure}
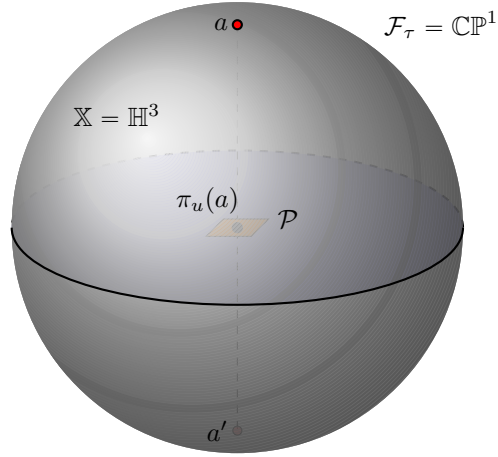

Figure \ref{fig:FibrationDoD} illustrates this construction in the rank one case $G=\PSL(2,\C)$, for a totally geodesic immersion $u$. The  associated symmetric space $\mathbb{H}^3$ is depicted with Poincaré's ball model. Since $\mathbb{H}^3$ has rank $1$, its visual boundary contains a single orbit $\mathcal{F}_\type\simeq\mathbb{CP}^1$.  The image of $u$ is the disk bounded by the equator. The pencil $\mathcal{P}$ is depicted as a parallelogram. its $\type$-base is a fiber of the fibration, and is the co-dimension $2$ submanifold $\mathcal{B}_\type(\mathcal{P})=\lbrace a,a'\rbrace$.

\begin{rem}
Note that if some element $g\in G$ preserves $u(M)$, then the map $\pi_u\circ u$ commutes with the action of $g$. In particular if $M=\tilde{N}$ for a compact manifold $N$ with fundamental group $\Gamma$ and if $u$ is $\rho$-equivariant for some $\rho:\Gamma\to G$, then $\pi_u$ is $\rho$-equivariant, and hence defines a fibration $\overline{\pi}_u:\Omega_u^\type/\rho(\Gamma)\to N$.
\end{rem}

\medskip

The two important steps in the proof of Theorem \ref{thm:FoliationDoD} are to check that the fibers are distinct and far enough from one another using Lemma \ref{lem:FibrationCriticalPoints}, and that these fibers are smooth manifolds using Lemma \ref{lem:SmoothSymmetricPencils}.

\medskip

\begin{proof}[Proof of Theorem \ref{thm:FoliationDoD}]
Consider the set:
$$E=\lbrace(a,x)\in \Omega^\type_u\times M|\mathrm{d}_x(b_{a,o}\circ u)=0\rbrace.$$

Because of Proposition \ref{lem:FibrationCriticalPoints} the Hessian of $b_{a,o}\circ u$ is non-degenerate at critical points, hence $E$ is locally the zero set of a submersion so it is a codimension $2$ submanifold of $\Omega^\type_u\times M$. 


Let $\pi_1:\Omega^\type_u\times M\to \Omega^\type_u$ and $\pi_2:\Omega^\type_u\times M\to M$ be respectively the projections onto the first and second factor.

\medskip

Lemma \ref{lem:FibrationCriticalPoints} implies that $\pi_1$ restricted to $E$ is a bijection. Moreover, again because of the non-degeneracy of the Hessian of $b_{a,o}\circ u$, the tangent space $T_{(a,x)}E$ at $(a,x)\in E$ intersects trivially $T_xM\subset T_{(a,x)}\left(\Omega^\type_u\times M\right)$. Hence $\pi_1$ restricted to $E$ is a local diffeomorphism, and therefore a diffeomorphism.

\medskip

Let $(a,x)\in E$. By definition $\mathrm{d}_a(b_{a,o}\circ u):\mathrm{v}\mapsto \langle \mathrm{d}u(\mathrm{v}),\mathrm{v}_{a,u(x)} \rangle_{u(x)}$ vanishes, so $\mathrm{v}_{a,u(x)}\perp \mathrm{d}u(T_xM)=\mathcal{P}_x$. Hence $a$ belongs to the $\type$-base of $\mathcal{P}_x$. Because of Proposition \ref{prop:pRegularSurface}, this pencil is $\type$-regular and hence its $\type$-base contains no singular points. Lemma \ref{lem:SmoothSymmetricPencils} implies that the tangent space $T_{(a,x)}E$ at $(a,x)\in E$ intersects trivially $T_a\Omega^\type_u\subset T_{(a,x)}\left(\Omega^\type_u\times M\right)$. The map $\pi_2$ restricted to $E$ is therefore a submersion at $(a,x)$.

\medskip

As a conclusion, $\pi_u=\pi_2\circ  \pi_1^{-1}$ is a smooth submersion. The $\type$-base of the pencil $\mathcal{P}_x$ is compact in $\mathcal{F}_\type$, and it is included in $\Omega^\type_u$ because of Lemma \ref{lem:FibrationCriticalPoints}. Hence $\pi_u$ is a proper submersion over a connected manifold: by the Ehresmann fibration theorem it is a fibration.
\end{proof}

\subsection{Totally geodesic immersions that are not nearly geodesic.}

In this subsection, let $u:M\to \X$ be a complete totally geodesic immersion. Let $\type\in \Sph\mathfrak{a}^+$. We don't assume in this subsection that $u$ is $\type$-regular, and hence $\type$-nearly geodesic. 

\medskip

One can still define $\Omega_u^\type$ as the set of $a\in \mathcal{F}_\type$ such that $b_{a,o}\circ u$ is proper and bounded from below, but we can't always expect the domain to have compact fibers in this case. Lemma \ref{lem:FibrationCriticalPoints} can be adapted as follows:

\begin{lem}
\label{lem:TotallyGeodesicDoD}
A point $a\in\mathcal{F}_\type$ belongs to $\Omega^\type_u$ if and only if the function $b_{a,o}\circ u$ admits a critical point $\pi_u(a)\in M$ at which the Hessian is positive. In this case the critical point is unique and is a global minimum of $b_{a,o}\circ u$. The domain $\Omega^\type_u$ is open.
\end{lem}

\begin{proof}
Let $a\in \Omega^\type_\rho$, and let $y\in M$ a point at which $b_{a,o}\circ u$ has a global minimum. The function $b_{a,o}\circ u$ is convex, but not necessarily strictly convex. Assume that the Hessian of $b_{a,o}$ in the direction $\mathrm{d}u(\mathrm{v})$ vanishes for some $\mathrm{v}\in T_yM$.

\medskip

There must exist a flat that contains $a$ and $\mathrm{d}u(\mathrm{v})$ by Lemma \ref{lem:HessianofBusemannFunctions}. Let $\eta$ be the  geodesic ray starting at $\mathrm{d}u(\mathrm{v})$ in $\X$. The function $b_{a,o}$ is linear on $\eta$ since $a$ and $\eta$ belong to a common flat. However the derivative of $b_{a,o}$ along $\eta$ vanishes at $u(y)$, so $b_{a,o}$ is constant along $\eta$. Moreover $u$ is totally geodesic and the whole geodesic ray starting at $\mathrm{d}u(\mathrm{v})$ in $\X$ belongs to the image of $u$, so $b_{a,o}\circ u$ is not proper.

\medskip

For all $a\in \Omega^\type_\rho$ the functions $b_{a,o}\circ u$ are convex and strictly convex at the critical point, which is therefore unique. The rest of the proof goes as in Lemma \ref{lem:FibrationCriticalPoints}.
\end{proof}

We define a map $\pi_u:\Omega^\type_u\to M$, using Lemma \ref{lem:TotallyGeodesicDoD}. We show that this map is a fibration. Recall that the regular base $\mathcal{B}^*_\type(\mathcal{P})$ defined in Section \ref{sec:Pencils} is a subset of the base $\mathcal{B}_\type(\mathcal{P})$ that is always a smooth codimension $d$ submanifold.

\begin{thm}
\label{thm:TotallyGeodesicFibration}
Let $a\in \Omega^\type_u$, the function $b_{a,o}\circ u$ admits a unique critical point denoted by $\pi_u(a)\in M$.
The map $\pi_u:\Omega_u^\type\to M$ is a smooth fibration, and the fiber of this map at $y\in M$ is the regular base $\mathcal{B}^*_\type(Tu(T_yM))$.
\end{thm}

\begin{proof}
By the same argument as for Theorem \ref{thm:FoliationDoD},
$\pi_u$ is a smooth submersion. 

\medskip

However we need to proceed differently to prove that this map is a fibration, since the fiber is not necessarily compact.
Let $g\in G$ be an element that stabilizes $u(M)\subset X$. The map $\pi_u$ is equivariant with respect to $g$, \ie for all $a\in \Omega^\type_\rho$:
$$\pi_u(g\cdot a)=g\cdot \pi_u(a).$$

 Let $y\in M$. Recall that the exponential map for the Lie group $G$ defines a map $\exp:T_{u(y)}\X\simeq \mathfrak{p}_{u(y)}\subset \mathfrak{g}\to G$. Moreover since $u(M)$ is totally geodesic any element of $\exp\left(\mathrm{d}u(T_yM)\right)$ is a transvection on this totally geodesic subspace hence it stabilizes $u(M)$. We consider the following map :
 \begin{align*} 
  \phi:T_yM\times \mathcal{B}^*_\type(Tu(T_yM))&\longrightarrow \Omega^\type_u \\
  (\mathrm{v},a)& \longmapsto \exp\left(\mathrm{d}u(\mathrm{v})\right)\cdot a 
  \end{align*}
  
This map is an immersion, between spaces of equal dimension. Moreover it is a bijection, hence it is a diffeomorphism. Through the identification $\exp:T_yM\to M$, this gives $\Omega^\type_u$ the structure of a fibration with projection $\pi_u$.
\end{proof}

Since the regular base is open, it is compact if and only if the regular points form a union of connected component of $\mathcal{B}_\type(\mathcal{P})$. This is for instance the case if $G=\PSL(3,\R)$ and $\mathcal{P}=\mathcal{P}_\text{irr}$ as in Example \ref{exmp:PencilsRegular}.

\begin{exmp}
\label{exmp:Nonregular}
Let $G=\SL(3,\R)$, and $\rho$ be a representation of the form $\rho=\iota_\text{irr}\circ \rho_0$ for some Fuchsian representation $\rho_0:\Gamma_g\to \SO(1,2)\simeq \PSL(2,\R)$ of a surface group and the natural inclusion $\iota_\text{irr}:\SO(1,2)\to \SL(3,\R)$. This representation admits a $\rho$-equivariant totally geodesic map $u:\widetilde{S_g}\to \X$ (see Section \ref{subsec:TotallyGeodesic} for more details). The pencil $\mathcal{P}=Tu(T_y\widetilde{S_g})$ for any $y\in \widetilde{S_g}$ is up to the action of $G$ equal to the pencil $\mathcal{P}_\text{irr}$ defined in Example \ref{exmp:PencilsRegular}.

\medskip

This pencil is not $\type_1$ regular, so $u$ is not $\type_1$-nearly geodesic. However because of Theorem \ref{thm:TotallyGeodesicFibration} the domain $\Omega_u^{\type_1}$ fibers over $\widetilde{S_g}$ with base $\mathcal{B}^*_{\type_1}(\mathcal{P}_\text{irr})$, which is a point in $\mathcal{F}_{\type_1}\simeq \mathbb{RP}^2$. This domain is the disk of positive vectors for the chosen bilinear form of signature $(1,2)$ on $\R^3$. In this example the regular points of $\mathcal{B}_{\type_1}(\mathcal{P}_\text{irr})$ form a connected component so the fibration is proper. If we consider a point $\ell\in \mathbb{RP}^2$ outside of the closure of this disk, the associated Busemann function is minimal in $u(\tilde{S_g})$ on a full geodesic line. If $\ell$ is in the boundary of the disk, the associated Busemann function is not bounded from below on $u(\tilde{S_g})$.

\end{exmp}

\subsection{Comparison with metric thickenings.}
\label{subsec:DoDKLP}
In this section we consider the case when $M=\widetilde{N}$ for some compact manifold $N$ with fundamental group $\Gamma$ and $u$ is equivariant with respect to a representation $\rho:\Gamma\to G$. In other words $\rho$ is a $\type$-nearly Fuchsian representation, as we defined in Section \ref{sec:AnosovNearlyFuch}. 

\medskip

 We show that if we have a $\type$-nearly Fuchsian $\rho$-equivariant map $u:\widetilde{N}\to \X$ for a representation $\rho:\Gamma\to G$ the domain $\Omega_u^\type$ coincides with a domain of discontinuity associated to Anosov representations constructed by Kapovich, Leeb, Porti \cite{KLP-DoD}.

\medskip

 The domain $\Omega^\type_\rho:=\Omega^\type_u$ depends on $\type$ and $\rho$ but not on $u$. Indeed a $1$-Lipshitz function on $\X$ is proper on the image of $u$ if and only if it is proper on any $\rho(\Gamma)$ orbit in $\X$.

 \medskip
 
 Even though this will be a consequence of Theorem \ref{thm:ComparaisonTotallyGeodesic}, one can easily check that the existence of the fibration of $\Omega^\type_\rho$ implies that it is a cocompact domain of discontinuity.

\begin{thm}
\label{thm:DoD}
Let $\rho:\Gamma\to G$ be a representation that admits an equivariant $\type$-nearly geodesic immersion $u:\widetilde{N}\to \X$. The action of $\Gamma_g$ on $\Omega_\rho^\type$ via $\rho$ is properly discontinuous and co-compact.
\end{thm}

\begin{proof}

Let $\pi_u$ be the fibration from Theorem \ref{thm:FoliationDoD}. Let $A$ be a compact subset of $\Omega^\type_\rho$. Its image $\pi_u(A)\subset \widetilde{N}$ is compact on $\widetilde{N}$. Since $\Gamma$ acts properly on $\widetilde{N}$, all but finitely many $\gamma\in \Gamma$ satisfy $\pi_u(A)\cap \gamma \pi_u(A)=\emptyset$. Hence for all but finitely many $\gamma\in \Gamma$ : $A\cap \rho(\gamma) A=\emptyset$. The action of $\Gamma$ via $\rho$ is therefore properly discontinuous.

\medskip

Let $D$ be a compact fundamental domain for the action of $\Gamma$ on $\widetilde{N}$, \ie a compact set that satisfies:
$$\bigcup_{\gamma\in \Gamma}D=\widetilde{N}. $$ 

The set $\pi_u^{-1}(D)$ is a fundamental domain for the action of $\Gamma$ on $\Omega^\type_\rho$ by $\rho$ by the equivariance of $\pi_u$. It is closed in $\Omega^\type_\rho$.
Moreover $\pi_u^{-1}(D)$ is closed in $\mathcal{F}_\type$. Indeed if we consider a sequence $(a_n)$ of elements of $\Omega^\type_\rho$ that converge to $a\in \mathcal{F}_\type$ such that $\pi_u(a_n)$ always belong to $D$, one can assume that $\pi_u(a_n)$ converges to $y_0\in D$ up to taking a subsequence. In the limit, one has $b_{a,o}\circ u(y_0)\leq b_{a,o}\circ u(y)$ for all $y\in \widetilde{N}$. Hence $y_0$ is a critical point for $b_{a,o}\circ u$ so by Lemma \ref{lem:FibrationCriticalPoints} $a\in \Omega^\type_\rho$.
 
 \medskip

Hence $\Omega^\type_\rho$ admits a compact fundamental domain for the action of $\Gamma$ via $\rho$, therefore this action is co-compact.

\end{proof}

We consider the domains of discontinuity constructed by \emph{metric thickenings}, which are particular instances of the domains of discontinuity associated with a \emph{Tits-Bruhat ideal} defined in \cite{KLP-DoD}.

\medskip

Let $(\type,\typeS)$ be a pair of elements in $\Sph\mathfrak{a}^+$. This pair will be called \emph{balanced} if $\typeS$ is $\type$-regular \ie:
$$\typeS\notin \bigcup_{w\in W}(w\cdot \type)^\perp .$$

Note that this is equivalent to $\type$ being $\typeS$-regular. Using the Tits angle $\Tangle:\partialvis \X ^2\to [0,\pi]$, see Section \ref{sec:symspaces}, we associate to any $b\in\partialvis \X$ a thickening $K_b\subset\mathcal{F}_\type$ defined as:
$$K_b=\left\lbrace a\in \mathcal{F}_\type|\Tangle(a,b)\leq \frac{\pi}{2} \right\rbrace.$$

Recall that the Tits angle was defined in Section \ref{sec:symspaces}, and is defined for points in $\partialvis \X$.

\begin{lem}
\label{lem:thickeningConnected}
Let $b_1,b_2$ belong to a common maximal facet in $\partialvis \X$, and suppose that their Cartan projections lie in the same connected component of:
 $$\Sph\mathfrak{a}^+\setminus \bigcup_{w\in W} (w\cdot \type)^\perp.$$
Then $K_{b_1}=K_{b_2}$.
\end{lem}

\begin{proof}
Let $f\in \mathcal{F}_\Delta$ be a maximal facet that contains $b_1$ and $b_2$. Let $a\in \mathcal{F}_\type$, there exist a  maximal flat of $\X$ such that $f$ and $a$ belong to its visual boundary. This flat can be identified with $\mathfrak{a}$ so that $f$ corresponds to $\partialvis \mathfrak{a}^+$. 

\medskip

Hence as long as $b$ lies in the visual boundary of this flat, $\Tangle(a,b)$ is equal to the Euclidean angle in the flat, so the sign of its cosine does not vary as long as the Cartan projection of $b$ does not lie in $(w\cdot \type)^\perp$ for any $w\in W$. Therefore if the Cartan projections of $b_1$ and $b_2$ are in the same connected component of the complement, $K_{b_1}=K_{b_2}$.
\end{proof}

 Recall that the flag manifold $\mathcal{F}_\typeS$ is $G$-equivariantly diffeomorphic to $\mathcal{F}_{\Theta(\typeS)}=G/P_{\Theta(\type_0)}$ where $\Theta(\typeS)$ is the set of simple root that do not vanish on $\typeS$. Hence given a flag $f\in \mathcal{F}_{\Theta(\type_0)}$ one can define $K_f^\typeS=K_b\subset \mathcal{F}_\type$ for the unique $b\in \mathcal{F}_\typeS$ corresponding to $f$.

\medskip

Given a pair $(\type,\typeS)$ and a $\Theta(\typeS)$-Anosov representation (see Definition \ref{defn:Anosov}), Kapovich, Leeb and Porti define a domain of discontinuity in $\mathcal{F}_\type$. 

\begin{thm}[{\cite{KLP-DoD}, Theorem 1.10}]
\label{thm:KLPDoD}
Let $(\type,\typeS)$ be a pair of elements of $\Sph\mathfrak{a}^+$. Let $\rho:\Gamma_g\to G$ be a $\Theta(\typeS)$-Anosov representation.  The following is a domain of discontinuity for $\rho$:
$$\Omega_\rho^{(\type,\typeS)}=\mathcal{F}_\type\setminus \bigcup_{\zeta\in \partial\Gamma_g}K^\typeS_{\xi_\rho^{\Theta}(\zeta)}.$$

Moreover if $(\type,\typeS)$ is balanced, then the action of $\Gamma_g$ via $\rho$ on $\Omega_\rho^{(\type,\typeS)}$ is co-compact.
\end{thm}

This theorem is a particular case of their result concerning Tits-Bruhat ideals. In \cite{KLP-DoD} it is explained how a pair $(\type, \typeS)$ yields a Tits-Bruhat ideal, defined via a \emph{metric thickening}. This ideal is balanced if and only if the pair is balanced in our sense.

\medskip

For a $\type$-nearly Fuchsian representation, the domain $\Omega_\rho^\type$ is always equal to some domain obtained by metric thickening. More precisely:

\begin{thm}
\label{thm:ComparaisonTotallyGeodesic}
Let $\rho$ be a $\type$-nearly Fuchsian representation of a non-elementary hyperbolic group. recall that $\sigma^\type_\rho$ and $\Theta(\sigma^\type_\rho)$ were defined in Section \ref{sec:AnosovNearlyFuch}.
Let $\typeS\in \sigma^\type_\rho$ be any element such that $\Theta(\typeS)=\Theta(\sigma^\type_\rho)$, whose existence is provided by Lemma \ref{lem:ConnectedRoot}. 
$$\Omega_\rho^\type=\Omega^{(\type,\typeS)}_\rho. $$
\end{thm}

 The Theorem \ref{thm:KLPDoD} from \cite{KLP-DoD} is a domain of discontinuity since $\rho$ is $\Theta(\sigma_\rho^\type)$-Anosov by Theorem \ref{thm:NearlyfuchsianAnosov}, and this domain is cocompact since the pair $(\type,\typeS)$ is balanced. 

\begin{proof}
Let us write $\Theta=\Theta(\sigma^\type_\rho)=\Theta(\typeS)$. Let $a\in\mathcal{F}_\type\setminus \Omega^\type_\rho$ and let $(y_n)_{n\in \N}$ be a diverging sequence of points in $M$ such that $\left(b_{a,o}(u(y_n))\right)_{n\in \N}$ is bounded from above. Up to taking a subsequence let us assume that it converges to a point $\zeta\in \partial\Gamma\simeq\partial \widetilde{N}$. We consider the geodesic segments $[o,u(y_n)]\subset \X$ for $n\in \N$. 
Since $\rho$ is $\type$-nearly Fuchsian it is a quasi isometric embedding by Proposition \ref{prop:QuasiGeod} so in particular the length of these segments goes to $+\infty$. Up to taking a subsequence, we can assume that these geodesic segments converge to a geodesic ray $\eta:\R_{\geq 0} \to \X$ with $\eta(0)=o$. Let $[\eta]\in \partialvis \X$ be the point corresponding to the class of $\eta$.

Busemann functions are convex so the function $b_{a,o}$ is bounded from above on all the geodesic segments $[o,u(y_n)]$ for $n\in \N$ and hence $b_{a,o}\circ \eta$ is bounded from above. Therefore $\Tangle (a,b)\leq \frac{\pi}{2}$ by Lemma \ref{lem:asymptoticsBusemann}, so $a\in K_{[\eta]}$.
Let $b\in \mathcal{F}_\typeS$ be an element such that $b$ and $[\eta]$ belong to a common maximal facet and whose Cartan projection lies in $\sigma_\rho^\type$. Lemma \ref{lem:thickeningConnected} implies that $K_{[\eta]}=K_b$.
 Theorem \ref{thm:AnosovLimitMap} implies that $K_{b}=K^\typeS_{\xi^{\Theta}_\rho(\zeta)}$. Therefore if $a\in \mathcal{F}_\type\setminus \Omega^\type_\rho$ then $a\in \mathcal{F}_\type\setminus \Omega^{(\type,\typeS)}_\rho$.

\medskip

Conversely let $a\in \Omega^\type_\rho$ and let $\zeta\in \partial\Gamma$. Consider a geodesic ray $\eta:\R_{>0}\to \widetilde{N}$ for the metric $u^*(g_\X)$ converging to $\zeta$. Theorem \ref{thm:MorseLemma} implies that there exist $D>0$ such that for all $t>0$, there exist a geodesic ray $\eta_t :\R_{>0}\to \X$ such that $\eta_t(0)= u\circ\eta(0)$, $\eta_t(t)$ is at distance at most $D$ of $u\circ \eta(t)$ and $[\eta_t]\in \partialvis\X$ belongs to a common maximal facet with $\xi_\rho^\Theta(\zeta)$.
Since $\mathcal{C}_\rho\subset \sigma_\rho^\type$, for all $t$ large enough, $K_{[\eta_t]}=K^\typeS_{\xi^\Theta_\rho(\zeta)}$.

The Busemann function $b_{a,o}$ is proper on $\eta$ hence for $t$ large enough $b_{a,o}\circ\eta_t(t)>b_{a,o}\circ\eta_t(0)$. Since $b_{a,o}$ is convex, this implies that $b_{a,o}$ is growing at least linearly on $\eta_t$, so $a\notin K_{[\eta_t]}$ by Lemma \ref{lem:asymptoticsBusemann}. Therefore $a\in \Omega_\rho^{(\type,\typeS)}$: this concludes the proof.
\end{proof}

\subsection{Invariance of the topology.}
\label{subsec:Ehresman Principle}

In this section we prove that the topology of the quotient of the domains of discontinuity considered by Kapovith-Leeb-Porti is not varying when the representation is deformed continuously. Guichard and Wienhard proved this already for the domaisn of discontinuity that they consider in \cite{GWDod}.

\medskip

Let $\Gamma$ be a torsion-free finitely generated group and $\mathcal{F}$ a $G$-homogeneous space. Let $(\rho_t)_{t\in [0,1]}$ a smooth family of representations from $\Gamma$ to $G$. Consider for every $t\in [0,1]$ an open $\rho_t(\Gamma)$-invariant domain $\Omega_t\subset \mathcal{F}$ .

\begin{lem}
\label{lem:EhresmannPrinciple}
Suppose that these domains are \emph{uniformly} co-compact domains of discontinuity for $(\rho_t)$, \ie the domain $\Omega=\lbrace(t,a)|a\in \Omega_t\rbrace\subset [0,1]\times \mathcal{F}$ is open and the action of $\Gamma$ via $\rho$ is properly discontinuous and co-compact where $\rho(\gamma)\cdot (t,a)=(t,\rho_t(\gamma)\cdot a)$. The quotient $\Omega_0/\rho_0(\Gamma)$ is diffeomorphic to $\Omega_1/\rho_1(\Gamma)$.
\end{lem}

\begin{proof}
The projection onto the first factor in $[0,1]\times \mathcal{F}$ descends to a submersion $p:\Omega/\rho(\Gamma)\to \R$. Since $\Omega/\rho(\Gamma)$ is compact and the base is connected, Ehresmann's fibration theorem implies that the proper submersion $p$ is a fibration. Hence $p^{-1}(0)$ and $p^{-1}(1)$ are diffeomorphic. 
\end{proof}

\begin{rems}
A concrete way to construct this diffeomorphism is to pick a Riemannian metric on $\Omega/\rho(\Gamma)$, and consider the flow of the gradient of $p$.
If we consider two different Riemannian metrics, the diffeomorphisms obtained are isotopic, since the space of Riemannian metrics is path connected. Hence this operation constructs a unique diffeomorphism up to isotopy.

\medskip

Note that a family of cocomact domains of discontinuity could be non-uniformly cocompact, for instance a family of representations so that $\rho_t:\mathbb{Z}\to \PSL(2,\R)$ is hyperbolic for $0\leq t<1$ and parabolic for $t=1$. These representations admit a unique maximal domain of discontinuity in $\mathbb{RP}^1$ with two connected components for $t<1$ and one for $t=1$. The quotient of the corresponding domain $\Omega$ is homeomorphic to the non compact space $S^1\times [0,1]\sqcup S^1 \times [0,1)$.
\end{rems}

In order to apply Lemma \ref{lem:EhresmannPrinciple}, we need a slight adaptation of Theorem \ref{thm:KLPDoD} from \cite{KLP-DoD}. Let $\Gamma$ be any Gromov-hyperbolic group. We check that the domains constructed by Kapovich, Leeb and Porti are uniformly cocompact domains of discontinuiy for any smooth path of Anosov representations.

\begin{prop}[{Adaptation of {\cite{KLP-DoD}, Theorem 1.10}}]
\label{prop:DoDUniform}
Let $(\type,\type_0)$ be a balanced pair as in Section \ref{subsec:DoDKLP}. Let $\rho:[0,1]\to \Hom(\Gamma,G)$, $t\mapsto \rho_t$ be a continuous path such that the family $(\rho_t)_{t\in [0,1]}$ consists only of $\Theta(\typeS)$-Anosov representations.

The family of domains $\Omega_t=\Omega_{\rho_t}^{(\type,\typeS)}$ for $t\in [0,1]$ are uniformly co-compact domains of discontinuity for the family of representations $\rho$.
\end{prop}

We check that the arguments from Kapovich, Leeb and Porti are uniform on neighborhoods of Anosov representations. The same proof holds if one considers more generally domains of discontinuity constructed with balanced Tits-Bruhat ideals as in \cite{KLP-DoD}.

\medskip

\begin{proof}

The domain $\Omega$ is the complement in $[0,1]\times \mathcal{F}_\type$ of:
$$K_\rho=\bigcup_{t\in [0,1]} \lbrace t\rbrace \times K_{\rho,t},$$
$$K_{\rho,t}=\bigcup_{x\in \partial\Gamma}K^{\typeS}_{\xi_{\rho_t}^{\Theta}(x)}.$$

Since the boundary maps $\xi_{\rho}^{\Theta}$ are continuous and vary continuously when $\rho$ varies continuously in the space of $\Theta$-Anosov representations (see \cite[Section 6]{BPS}), and since $K^{(\type,\typeS)}_{\xi_{\rho}^{\Theta}(x)}$ is compact, $K_\rho$ is compact so $\Omega=\lbrace(t,a)|a\in \Omega_t\rbrace\subset [0,1]\times \mathcal{F}$ is open.

\medskip

 Let us fix a Riemannian distance $d$ on $\mathcal{F}_\type$. Let $A=\lbrace (t,a)|t\in[0,1], a\in A_t\subset \Omega\rbrace$ be a compact set and let $(\gamma_n)\in \Gamma$ be a diverging sequence. \cite[Corollary 6.8]{KLP-DoD} implies that given $t\in [0,1]$, for any $\epsilon>0$ for all $n$ large enough if $d(a, K_{\rho,t})\geq \epsilon$, then $d(\rho_t(\gamma_n)\cdot a, K_{\rho,t})\leq \epsilon$, where the minimal value of $n$ needed depends on the constants $b,c$ that come into play in the definition of Anosov representations (Definition \ref{defn:Anosov}).

  Since we consider a compact set of Anosov representations, and since these constants can be chosen locally uniformly around a given Anosov representation (see \cite[Theorem 7.18]{KLP}), these constants can be chosen uniformly for all representations $(\rho_t)_{t\in [0,1]}$. Moreover the compact sets $A_t$ are at uniform distance from $K_{\rho,t}$ Therefore, for all $n\in \mathbb{N}$ large enough, for all $t\in \R$ $\rho_t(\gamma_n)\cdot A_t\cap A_t=\emptyset$, so for all $n$ large enough $\rho(\gamma_n)\cdot A\cap A=\emptyset$. We have proven that the action of $\rho$ on $\Omega$ is properly discontinuous. 

\medskip

In order to prove the cocompactness of the action on $\Omega$, we will check that the transverse expansion holds uniformly. It follows from \cite[Proposition 7.7]{KLP-DoD} that for every $t\in [0,1]$, the action of ${\rho_t}$ is \emph{transversely expanding} at the limit set of $\rho_t$ as in \cite[Definition 5.21]{KLP-DoD}, \ie for all $x\in \partial\Gamma$, there exist $\gamma\in \Gamma$, an open neighborhood $U$ of $K^{\typeS}_{\xi_{\rho_t}^\Theta(x)}$ in $\mathcal{F}_\type$ and a constant $\lambda>1$ such that for all $a\in U$ and $y\in \partialvis \Gamma$ that satisfy $K^{\typeS}_{\xi_{\rho_t}^\Theta(y)}\subset U$ one has:
$$d\left(\rho_t(\gamma)\cdot a,\rho_t(\gamma)\cdot K^{\typeS}_{\xi_{\rho_t}^\Theta(y)} \right) \geq \lambda d\left(a,K^{\typeS}_{\xi_{\rho_t}^\Theta(y)}\right).$$

\medskip

 Let $g\in G$, $f\in \mathcal{F}_\Theta$ $\lambda >1$ and $U\subset \mathcal{F}_\type$ be an open set. We say that $g$ \emph{is exapnding} at $K_f^\typeS$ over $U$ with factor $\lambda$ if for all $a\in U$:
$$d\left(g\cdot a, g\cdot K^\typeS_f\right)\geq \mu d\left(a,K^\typeS_f\right).$$

This property is open in the following sense: if $g$ is expanding at $K_f^\typeS$ over $U$ with factor $\lambda$, then for any $1<\lambda'<\lambda$ and any open subset $E$ such that $\overline{E}\subset U$ there exist a neighborhood $U_g$ of $g$ in $G$ and $U_f$ of $f$ in $\mathcal{F}_\Theta$ such that for all $g'\in U_g$ and $f'\in U_f$, $g'$ is is expanding at $K_{f'}^\typeS$ over $E$ with factor $\lambda'$.

\medskip

 This implies that the action of $\rho$ on $\mathcal{F}_\type\times [0,1]$ satisfies the transverse expansion property where $K_\rho$ is considered as a bundle over $[0,1]\times K$, \ie for all $t\in [0,1]$, $x\in \partial \Gamma$, there exist $\gamma\in \Gamma$, an open neighborhood $U_0$ of $K^{\typeS}_{\xi_{\rho_t}^\Theta(x)}$ in $\mathcal{F}_\type\times [0,1]$ and a constant $\lambda'>1$ such that for all $a\in U$ and $y\in \partial \Gamma$ that satisfy $K^{\typeS}_{\xi_{\rho_t}^\Theta(y)}\times \lbrace t\rbrace \subset U_0$ one has:
$$d\left(\rho(\gamma)\cdot a,\rho_(\gamma)\cdot K^{\typeS}_{\xi_{\rho_t}^\Theta(y)}\times \lbrace t\rbrace \right) \geq \lambda' d\left(a,K^{\typeS}_{\xi_{\rho_t}^\Theta(y)}\times \lbrace t\rbrace\right).$$
 
 Therefore by \cite[Proposition 5.26]{KLP-DoD} the action of $\rho$ is cocompact on $\Omega$.
\end{proof}

From Proposition \ref{prop:DoDUniform} and Lemma \ref{lem:EhresmannPrinciple} we get the following corollary.

\begin{cor}
\label{cor:EhresmanApplied}
Assume that $\Gamma$ is torsion-free. Let $\mathcal{C}\subset \Hom(\Gamma,G)$ be an open and connected set consisting only of $\Theta$-Anosov representations for some $\Theta\subset \Delta$. Let $(\type,\typeS)$ be a balanced pair such that $\forall \alpha\in \Delta\setminus \Theta$, $ \alpha(\typeS)=0$. The diffeomorphism type of $\Omega_{\rho}^{(\type,\typeS)}/\rho(\Gamma)$ is independent of $\rho\in \mathcal{C}$

\medskip

If moreover $\mathcal{C}$ is simply connected, the diffeomorphism provided by Lemma \ref{lem:EhresmannPrinciple} between $\Omega_{\rho_1}^{(\type,\typeS)}/\rho_1(\Gamma)$ and $\Omega_{\rho_2}^{(\type,\typeS)}/\rho_2(\Gamma)$ for $\rho_1,\rho_2\in \mathcal{C}$ is uniquely determined up to isotopy.
\end{cor}

Connectedness via paths that are smooth by parts for an open set in $\Hom(\Gamma,G)$ is equivalent to connectedness since $\Hom(\Gamma,G)$ is locally a real algebraic variety.

\section{Applications.}

In this section we apply our results to prove that all representations in some connected components of Anosov representations are the restricted holonomy of a geometric structure on a fiber bundle over a manifold. For these applications we only consider nearly geodesic surfaces that are totally geodesic. We will mostly focus on surface groups, but in Section \ref{subsec:HigherDimention} we also describe two applications for representations of fundamental groups of higher dimensional compact hyperbolic manifolds. 

\subsection{Totally geodesic immersions.}
\label{subsec:TotallyGeodesic}
Totally geodesic surfaces provide examples of $\type$-nearly geodesic surfaces if these surfaces are $\type$-regular (Proposition \ref{prop:TotGeodNearlyGeod}). The study of totally geodesic surfaces in $\X$ is related to the study of representation of semi-simple Lie algebras in $\mathfrak{g}$. We recall here a classical fact.

\begin{prop}
\label{prop:TotallyGeodesic}
Let $\mathfrak{h}\subset \mathfrak{g}$ be a semi-simple Lie subalgebra of non-compact type. Let $H$ be the closed Lie subgroup of $G$ with Lie algebra $\mathfrak{h}$ and $\mathbb{Y}$ it's associated symmetric space of non-compact type. There exist a $H$-equivariant and totally geodesic embedding $u_\mathfrak{h} :\mathbb{Y}\to \X$. The image of this embedding is unique up to the action of the centralizer:
$$\mathcal{C}_G(\mathfrak{h})=\lbrace g\in G|\forall \mathrm{h}\in \mathfrak{h},\,\Ad_g(\mathrm{h})=\mathrm{h}\rbrace.$$

If $y\in \mathbb{Y}$, let $K\subset G$ be the stabilizer of $u_\mathfrak{h}(y)$ in $G$. Every element in $\mathcal{C}_K\left(\mathfrak{h}\right)$ of $\mathfrak{h}$ in $G$ fixes $u_\mathfrak{h}(\mathbb{Y})$ pointwise. If the centralizer $\mathcal{C}_G(\mathfrak{h})$ in $G$ is compact, then $\mathcal{C}_K\left(\mathfrak{h}\right)=\mathcal{C}_G(\mathfrak{h})$ so the totally geodesic submanifold $u_\mathfrak{h}(\mathbb{Y})\subset \X$ is uniquely determined.

\end{prop}

\begin{proof}
 Let $\mathfrak{h}=\mathfrak{t}+\mathfrak{p}$ be a Cartan decomposition of $\mathfrak{h}$ associated with the Cartan involution $\theta_y$ for $y\in \mathbb{Y}$. Let $\widehat{\mathfrak{h}}=\mathfrak{t}+i\mathfrak{p}$ be the associated compact real form of $\mathfrak{h}\otimes \C$. Since $\widehat{\mathfrak{h}}$ is the Lie algebra of a semi-simple compact Lie group, it is the Lie algebra of a compact Lie subgroup of $G^\C$. In particular there exist a Cartan involution $\theta^\C$ of $\mathfrak{g}\otimes \C$. such that $\theta^\C(\mathrm{v})=\mathrm{v}$ for all $\mathrm{v}\in\widehat{\mathfrak{h}}$.

\medskip 

Let $\overline{\theta^\C}$ be the Cartan involution conjugate to $\theta^\C$ in $\mathfrak{g}\otimes \C$. Let $\theta$ be the Cartan involution of $\mathfrak{g}\otimes \C$ corresponding to the midpoint in the symmetric space associated with $\mathfrak{g}\otimes \C$ of $\overline{\theta^\C}$ and $\theta^\C$. It is invariant by conjugation, so it descends to a Cartan involution on $\mathfrak{g}$, associated with a point $x\in \X$.
There exist a transvection along the geodesic between $\overline{\theta^\C}$ and $\theta^\C$ in the symmetric space associated to $\mathfrak{g}\otimes \C$ that induces a unique inner automorphism $\phi$ of $\mathfrak{g}\otimes \C$ such that $\phi\theta^\C\phi^{-1}=\theta$, $\phi^2\theta^\C\phi^{-2}=\overline{\theta^\C}$. Moreover $\phi$ is symmetric positive with respect to $B^\C\left(\cdot,\theta^\C(\cdot)\right)$, where $B^\C$ is the Killing form on $\mathfrak{g}\otimes \C$ and the transvection $\phi^4$ is equal to a composition of symmetries $\theta^\C\overline{\theta^\C}$. The transvection $\phi^4$ stabilizes $\mathfrak{h}\otimes \C$, so $\phi$ also stabilizes $\mathfrak{h}\otimes \C$. Therefore $\theta=\phi\theta^\C\phi^{-1}$ stabilizes $\mathfrak{h}\otimes\C$. But $\theta$ is a real Cartan involution as it is preserved by conjugation, so it stabilizes $\mathfrak{h}$. Let $x\in \X$ be the point corresponding to $\theta$ in $\X$.

\medskip

The map $u_\mathfrak{h}:\mathbb{Y}\to \X$ such that for all $h\in H$, $u_\mathfrak{h}(h\cdot y)=h\cdot x$ is well defined since $\mathfrak{t}\subset \mathfrak{t}_x$ so the image of the stabilizer of $y\in \mathbb{Y}$ by $H$ lies in the stabilizer of $x\in \mathbb{X}$. Moreover it is totally geodesic, since $\eta(\mathfrak{p})\subset \mathfrak{p}_x$ (see \cite{Helgason} Ch 4, Section 7). This map is by definition $H$-equivariant.

\medskip
Suppose that there is an other $H$-equivaiant and totally geodesic embedding $u_\mathfrak{h}'$ such that $u_\mathfrak{h}'(y)=x'$. Let $\theta'$ be the corresponding Cartan involution of $\mathfrak{g}$, then $\theta' \circ \theta$ is the identity on $\mathfrak{h}$ and is equal to the adjoint action of $\exp(\mathrm{z})$ for some $\mathrm{z}\in \mathfrak{p}_{x}$. Therefore $\mathrm{z}$ is in the centralizer of $\mathfrak{h}$ in $\mathfrak{g}$. 
Hence $g=\exp(\frac{\mathrm{z}}{2})\in 
\mathcal{C}_G(\mathfrak{h})$ satisfies $\Ad_g\circ u'_\mathfrak{h}=u_\mathfrak{h}$.
Conversely let $g'\in \mathcal{C}_K\left(\mathfrak{h}\right)$, then $g'$ fixes $u_\mathfrak{h}(y)\in \X$ and it fixes $\mathfrak{h}$ so it fixes $\mathrm{d}u_\mathfrak{h}(T_y\mathbb{Y})$. Therefore it preserves and acts trivially on $u_\mathfrak{h}(\mathbb{Y})$. 

\medskip

Now assume that $\mathcal{C}_G(\mathfrak{h})$ is compact: there cannot be any element $\mathrm{z}\in \mathfrak{p}_x$ in the centralizer of $\mathfrak{h}$ in $\mathfrak{g}$, so the totally geodesic and $H$-equivariant embedding $u_\mathfrak{h}$ is unique and $\mathcal{C}_K\left(\mathfrak{h}\right)=\mathcal{C}_G(\mathfrak{h})$.
\end{proof}

\begin{rem}
A totally geodesic embedding $u:\mathbb{Y}\to \mathbb{X}$ can only be a $\type$-nearly geodesic immersions if $\rank(\mathbb{Y})=1$, because otherwise it cannot be $\type$-regular for any $\type\in \Sph\mathfrak{a}^+$ (see Proposition \ref{prop:Rank2}). When $\rank(\mathbb{Y})=1$, all the unit tangent vectors to this embedded surface have the same Cartan projection, so the embedding is $\type$-regular for all $\type$ in the complement in $\Sph\mathfrak{a}^+$ of a finite collection of hyperplanes.
\end{rem}

We illustrate Proposition \ref{prop:TotallyGeodesic} in the following example for some special $\mathfrak{sl}_2$ Lie subalgebras in $\mathfrak{sl}_n(\R)$.

\begin{exmp}
\label{exmp:PSL2N}
In this example we construct representations $\iota$ from $\SL(2,\R)$ into $\SL(n,\mathbb{K})$ that stabilize some totally geodesic hyperbolic planes inside  the symmetric space $\X=\mathcal{S}_n$ associated with $G=\SL(n,\K)$ for $\mathbb{K}=\R$ or $\C$.

\medskip
Let $\mathbb{K}=\mathbb{R}$ or $\mathbb{C}$. Let $V_n(\mathbb{K})$ be the space of homogeneous polynomial with coefficients in $\mathbb{K}$ of degree $n-1$ in two variables $X$ and $Y$. To an element $g\in \SL_2(\mathbb{R})$ one can associate an element $\iota_\text{irr}(g)\in \SL\left(V_n(\mathbb{K})\right)$ that acts by a change of variable on $V_n(\mathbb{K})$ i.e that associates to $P\in V_n(\mathbb{K})$ the polynomial $P\circ g^{-1}$. Let $\mathfrak{h}$ be the corresponding $\mathfrak{sl}_2$-Lie subalgebra of $\mathfrak{sl}_n(\R)$, note that $\iota_\text{irr}=\iota_\mathfrak{h}$. 

\medskip

Let $q$ the Euclidean or Hermitian metric on $V_n(\mathbb{K})$ such that for all $0\leq a,b\leq n-1$:
 $$q\left(X^aY^{n-1-a},X^bY^{n-1-b}\right)=\begin{cases}
    \binom{n-1}{a} ^{-1}& \text{if } a=b,\\
    0            & \text{otherwise}.
\end{cases}$$

Consider the following basis of the lie algebra $\mathfrak{sl}_2(\mathbb{R})$ :

$$\mathbf{h}=\begin{pmatrix}
1 & 0 \\
0 & -1 
\end{pmatrix}, \mathbf{f}=\begin{pmatrix}
0 & 1 \\
1 & 0 
\end{pmatrix} , \mathbf{g}=\begin{pmatrix}
0 & 1 \\
-1 & 0 
\end{pmatrix}$$

which satisfies $[\mathbf{g},\mathbf{h}]=-2\mathbf{f}$, $[\mathbf{h},\mathbf{f}]=2\mathbf{g}$, $[\mathbf{f},\mathbf{g}]=2\mathbf{h}$.
Fix the following orthonormal basis for $q$ :

 $$(e_a)_{0\leq a\leq n-1}=\left(X^aY^{n-1-a}\binom{n-1}{a}^{\frac{1}{2}}\right)_{0\leq a\leq n-1}$$.
 
 For $0\leq a\leq n-1$ one has :
$$\mathrm{d}\iota_\text{irr}(\mathbf{f})(e_a)=(2a-n+1)e_a$$
$$\mathrm{d}\iota_\text{irr}(\mathbf{e})(e_a)=aX^{a-1}Y^{n-a}\binom{n-1}{a}^{\frac{1}{2}} - (n-1-a)X^{a+1}Y^{n-2-a}\binom{n-1}{a}^{\frac{1}{2}}$$
$$\mathrm{d}\iota_\text{irr}(\mathbf{e})(e_a)=\sqrt{a(n-a))}e_{a-1} - \sqrt{(a+1)(n-1-a)}e_{a+1}$$

and moreover $\mathbf{g}=\frac{1}{2}[\mathbf{f},\mathbf{h}]$. In particular $\mathrm{d}\iota_\text{irr}(\mathbf{f})$ and $\mathrm{d}\iota_\text{irr}(\mathbf{g})$ are symmetric or Hermitian and $\mathrm{d}\iota_\text{irr}(\mathbf{h})$ is anti-symmetric or anti-Hermitian with respect to $q$. 

\medskip

Hence as in the proof of Propostion \ref{prop:TotallyGeodesic}, there is a totally geodesic map $u_\mathfrak{h}:\mathbb{H}^2\to \mathcal{S}_n$ such that the image of the Cartan involution $M\mapsto -M^t$ is equal to the point in $\mathcal{S}_n$ corresponding to $q$. 

\medskip

This representation $\iota_\text{irr}$ is the unique irreducible representation of $\SL(2,\R)$ into $\SL(n,\R)$ up to conjugation by elements of $\GL(n,\R)$. The image of $\iota_\text{irr}$ lies in the subgroup $\Sp(2k,\R)$ for some symplectic form on $\R^{2k}$ when $n=2k$ is even, and in $\SO(k,k+1)$ for some quadratic form on $\R^{2k+1}$ when $n=2k+1$ is odd.

\bigskip

One can construct other totally geodesic hyperbolic planes in $\mathcal{S}_n$ by considering representations of $\SL(2,\R)$ that can be decomposed into a direct sum of irreducible representations. Equivalently one can consider other reducible $\mathfrak{sl}_2$-subalgebras of $\mathfrak{sl}_n(\R)$.

\medskip

 For instance one can define $\iota_\text{red}:\SL(2,\R)\to \SL(2n,\R)$ that associates to a matrix $M\in \SL(2,\R)$ the block diagonal matrix:
 $$\iota_\text{red}(M)=\begin{pmatrix}
M & &\\
 & \ddots&\\
 &&M
 \end{pmatrix} .$$

\medskip

The image of $\iota_\text{red}$ lies in $\Sp(2n,\R)$ for some symplectic form $\omega$ on $\R^{2n}$.

\end{exmp}

\subsection{Geometric structures on fiber bundles.}

Using the projection defined by Busemann functions from Theorem \ref{thm:FoliationDoD}, one can show that Anosov deformations of representation that admit an equivariant nearly geodesic immersion are holonomies of $(G,X)$ structures on a fiber bundle over $S_g$. 

\medskip

A $(G,X)$-structure on a manifold $N$, for a Lie group $G$ and a $G$-homogeneous space $X$ on which $G$ acts faithfully, is a maximal atlas of charts of $M$ valued in $X$ whose transition functions are the restriction of the action of some elements in $G$. A more developed introduction to this notion can be found in \cite{GXStruct}.

To a $(G,X)$-structure, one can associate a \emph{developing map} $\dev:\widetilde{N}\to X$, that is a local diffeomorphism compatible with the atlas defining the $(G,X)$-structure on $N$, and a \emph{holonomy} $\hol:\pi_1(N)\to G$ so that $\dev$ is $\hol$-equivariant. This pair is unique up to the action of $G$ by conjugation of the holonomy and post-composition of the developing map. 

\medskip

Let $N$ be a manifold and $\Gamma$ its fundamental group. In what follows, we say that a $(G,X)$-structure \emph{on a fiber bundle} $F$ over $N$ is a $(G,X)$-structure on $F$ for which the fundamental group of the fibers is included in the kernel of the holonomy. Hence one can define the \emph{restricted holonomy} of the structure as the quotient map $\rho:\pi_1(N)\to G$ induced by the holonomy. 

\medskip

Constructing domains of discontinuity allows us to construct geometric structures.

\begin{prop} 
Let $\rho:\Gamma\to G$ be a representation and $\Omega\subset X$ a co-compact non-empty domain of discontinuity which fibers $\rho$-equivariantly over $\widetilde{N}$. Any connected component of the quotient $\Omega/\rho(\Gamma)$ inherits a $(G,X)$-structure \emph{on a fiber bundle}, with restricted holonomy $\rho$.
\end{prop}

 Note that even if $\Omega$ is disconnected, the quotient $\Omega/\rho(\Gamma)$ can be connected. 
 
\medskip
 
From now on, we assume that $G$ is center-free, so it acts faithfully on its flag manifolds. Let $N$ be a compact manifold whose fundamental group $\Gamma$ is Gromov hyperbolic and torsion-free. Let $\tau \in \Sph\mathfrak{a}^+$.

\begin{thm}
\label{thm:Application}
Let $\rho_0:\Gamma\to G$ be a representation that admits an equivariant $\tau$-nearly geodesic surface $u:\widetilde{N}\to \X$ such that $\Omega^\type_{u}\neq \emptyset$. Let $\mathcal{C}$ be the connected component of the space of $\Theta(\sigma^\type_{\rho_0})$-Anosov representations in $\Hom(\Gamma,G)$ containing $\rho_0$. Every representation in $\mathcal{C}$ is the restricted holonomy of a $(G,\mathcal{F}_\type)$-structure \emph{on a fiber bundle} $F$ over $N$.
\end{thm}

\begin{proof}
Theorem \ref{thm:DoD} implies that the domains $\Omega^\type_u$ admits a $\rho_0$-equivariant fibration over $\widetilde{N}$. The domain $\Omega_{\rho_0}^\type$ coincides with a domain obtained as a metric thickening $\Omega^{(\type,\type_0)}_{\rho'}$ for some $\type_0\in \Sph\mathfrak{a}^+$ such that $\Theta(\type_0)=\Theta(\sigma^\type_{\rho_0})$. Let $F$ be a connected component of $\Omega_{\rho_0}^\type/\rho_0(\Gamma)$, Corollary \ref{cor:EhresmanApplied} implies that for every representation $\rho\in\mathcal{C}$, a connected component $F_\rho$ of $\Omega_\rho^{(\type,\type_0)}/\rho(\Gamma)$ is diffeomorphic to $F$, which is a fiber bundle over $M$. The covering map $\Omega_\rho^{(\type,\type_0)}\to\Omega_\rho^{(\type,\type_0)}/\rho(\Gamma)$ induces the covering $\widehat{F_\rho} \to F_\rho\simeq F$ associated to the subgroup of the fundamental group of $F$ corresponding to the fundamental group of the fiber of $F$ over $M$.

\medskip

Note that $F$ is assumed to be non-empty. The $(G,\mathcal{F}_\type)$-structure on $F_\rho\simeq F$ is such that the holonomy of the fundamental group of each fiber is trivial. Indeed the developing map $\dev:\widetilde{F_\rho}\to \mathcal{F}_\type$ descends to the inclusion $\widehat{F_\rho}\to \mathcal{F}_\type$, so the fundamental group of the fiber belongs to the kernel of the holonomy.
\end{proof}

Let $M=S_g$ be a closed orientable surface of genus $g$ and $\Gamma=\Gamma_g$ be its fundamental group. We apply Theorem \ref{thm:Application} in cases when the nearly geodesic surface is totally geodesic and we describe the fibration that is obtained.

\medskip

Let $\mathfrak{h}\subset \mathfrak{g}$ be a $\mathfrak{sl}_2$ Lie subalgebra, \ie a Lie subalgebra isomorphic to $\mathfrak{sl}_2(\R)$. Note that if $G$ is a quotient of its adjoint form, which is the case since we assumed that $G$ is center-free, the corresponding Lie group $H$ is isomorphic to $\SL(2,\R)$ or $\PSL(2,\R)$. We write $\iota_\mathfrak{h}:\SL(2,\R)\to G$ the corresponding Lie group representation, and $u_\mathfrak{h}:\mathbb{H}^2\to \X$ a corresponding equivariant totally geodesic embedding.

\begin{defn}
\label{defn:GeneralFuchs}
We say that a representation $\rho:\Gamma_g\to G$ is $\mathfrak{h}$-\emph{generalized Fuchsian} if it preserves and acts cocompactly on $u_\mathfrak{h}\left(\mathbb{H}^2\right)$. 
\end{defn}

If $G$ is center-free and, a $\mathfrak{h}$-generalized Fuchsian representation can be written $\gamma\mapsto \iota_\mathfrak{h}(\rho_0(\gamma))\times \chi(\gamma)$ for some Fuchsian representation $\rho_0:\Gamma_g\to\SL(2,\R)$ and some \emph{associated character} $\chi:\Gamma_g\to \mathcal{C}_K(\mathfrak{h})$.
\medskip

Let $y\in \mathbb{H}^2$ be a base-point, and let $K$ be the stabilizer in $G$ of $u_\mathfrak{h}(y)$. We write $\mathcal{B}^\type(\mathfrak{h})$ for the $\type$-base of the pencil of tangent vectors $\mathrm{d}u_\mathfrak{h}(T_y\mathbb{H}^2)$. Note that this pencil is stabilized by $\iota_\mathfrak{h}\left(\SO(2,\R)\right)\times \mathcal{C}_K(\mathfrak{h})$.

The quotient map  $\SL(2,\R)\to \SL(2,\R)/\SO(2,\R)\simeq \mathbb{H}^2$ defines a principal $\SO(2,\R)$-bundle $\widetilde{P}$ over $\mathbb{H}^2$. Let $P_{S_g}$ be its quotient via some Fuchsian representation. Given a character $\chi:\Gamma_g\to \mathcal{C}_K(\mathfrak{h})$, let $P_{S_g,\chi}\to S_g$ be the $\SO(2,\R)\times \mathcal{C}_K(\mathfrak{h})$-principal bundle obtained as the product of $P_{S_g}$ and the flat $\mathcal{C}_K(\mathfrak{h})$-bundle associated to $\chi$. 

\begin{thm}
\label{thm:HighTeich}
Suppose that $\mathfrak{h}$ is \emph{$\type$-regular}, \ie $u_\mathfrak{h}$ is $\type$-regular, and that $\Omega^\type_{u_\mathfrak{h}}\neq \emptyset$. Let $\rho:\Gamma_g\to G$ be a $\mathfrak{h}$-generalized Fuchsian representation with associated character $\chi$ and let $\type\in \Sph\mathfrak{a}^+$.

Let $\mathcal{C}$ be the connected component of the space of $\Theta(\sigma^\type_\rho)$-Anosov representations that contains $\rho$. Every representations in $\mathcal{C}$ is the restricted holonomy of a $\left(G,\mathcal{F}_{\type}\right)$-structure \emph{on a fiber bundle} $F$ over $S_g$.

\medskip

The fiber bundle $F$ is a connected component of the reduction of the $\SO(2,\R)\times \mathcal{C}_K(\mathfrak{h})$-principal bundle $P_{S_g,\chi}$ over $S_g$ via the action of $\iota_\mathfrak{h}\left(\SO(2,\R)\right)\times \mathcal{C}_K(\mathfrak{h})$ on $\mathcal{B}^\type(\mathfrak{h})$. 
\end{thm}

Note that when $\Theta\subset \Delta$ is a Weyl orbit of simple roots and $\type=\type_\Theta$, a $\type_\Theta$-regular immersion is just a $\Theta$-regular immersion and $\Theta(\sigma^\type_\rho)=\Theta$.

\begin{rem}
If we replace $G$ and $\mathcal{F}_\type$ by finite covers $\widehat{G}$ and $\widehat{\mathcal{F}_\type}$  so that $\widehat{G}$ acts faithfully on $\widehat{\mathcal{F}_\type}$, Theorem \ref{thm:HighTeich} still applies. In this case one should replace $\mathcal{B}^\type(\mathfrak{h})$ by its pre-image $\widehat{\mathcal{B}^\type(\mathfrak{h})}$ by the covering map $\widehat{\mathcal{F}_\type}\to \mathcal{F}_\type$.
\end{rem}

The only part of Theorem \ref{thm:HighTeich} that is not already contained in \ref{thm:Application} is the description of the fibration.

\medskip

\begin{proof}

The embedding $u_\mathfrak{h}$ is totally geodesic. Moreover it is assumed to be $\type$-regular. Hence $u_\mathfrak{h}$ is $\type$-nearly geodesic, so one can apply Theorem \ref{thm:Application}. It only remains to describe the fibration.

\medskip
 The fibration $\pi:\Omega_{u_\mathfrak{h}}^\type\to \mathbb{H}^2$ is $\iota_\mathfrak{h}\left(\SL(2,\R)\right)$-equivariant. The corresponding fiber bundle can be identified in a $\SL(2,\R)$-equivariant way as the reduction of the $\SO(2,\R)$-principal bundle $\widetilde{P}$ by the action of $\iota_\mathfrak{h}\left(\SO(2,\R)\right)$ on the fiber $\mathcal{B}^\type(\mathfrak{h})$. The quotient of this fiber bundle by any Fuchsian representation $\rho_0:\Gamma_g\to \SL(2,\R)$ is hence the bundle induced by the principal bundle $P_{S_g}$.
Once we quotient $\Omega_{u_\mathfrak{h}}^\type$ by $\rho=\iota_\mathfrak{h}\circ \rho_0 \times \chi$ the quotient becomes the twisted fiber bundle induced by $P_{S_g,\chi}$.
\end{proof}

\subsection{Higher rank Teichmüller spaces.}
\label{subsec:ThetaPositive}

In this section we apply Theorem \ref{thm:HighTeich} to Hitchin representations in $\PSL(n,\R)$ and to maximal representations in $\Sp(2n,\R)$. Then we explain how in general one can apply it to the connected components of $\Theta$-positive representations containing at least one generalized Fuchsian representation associated to a $\Theta$-principal $\mathfrak{sl}_2$-Lie subalgebra. 

\subsubsection{Positive representations.}

Let $G$ be a connected simple Lie group of non-compact type with trivial center and $\Theta$ a set of its simple roots such that the pair $(G,\Theta)$ admits a notion of $\Theta$-positivity in the sense of \cite{GWpositivity22}. We moreover assume that $G$ is not locally isomorphic to $\PSL(2,\R)$. This means that either :
\begin{itemize}
\item[(i)] $G$ is split real and $\Theta=\Delta$,
\item[(ii)] $G$ is Hermitian of tube type,
\item[(iii)] $G$ is locally isomorphic to $\SO(p,q)$ and $\Theta=\lbrace\alpha_1,\cdots, \alpha_p\rbrace$
\item[(iv)] $G$ is a real form of the complex Lie group with Dynkin diagram $F_4$, $E_6$, $E_7$ or $E_8$ with restricted Dynkin diagram $F_4$ and $\Theta$ consists of the $2$ larger roots.
\end{itemize}

For any of these pairs, Guichard and Wienhard \cite{GWpositivity22} constructed a connected component $U$ in the space and transverse triples of elements in $G/P_\Theta$. They call such triples \emph{positive triples}, and a representation $\rho:\Gamma_g\to G$ is called $\Theta$-positive if it admits a continuous and $\rho$-equivariant map $\xi:\partial\Gamma\to G/P_\Theta$ so that for all distinct triple of points $(x,y,z)\in \partial\Gamma^{(3)}$, $(\xi(x),\xi(y),\xi(z))\in U$. They prove in particular that such representations are $\Theta$-Anosov. Moreover the space of $\Theta$-positive representation is closed in the space of representation that do not virtually factor through a parabolic subgroup, by a work of Guichard, Labourie and Wienhard \cite{positivityLabourie}.

\medskip

A $\Theta$-principal Lie subalgebra $\mathfrak{h}_\Theta$ for a pair $(G,\Theta)$ that admits a notion of $\Theta$-positivity is a principal subalgebra of the split Lie subalgebra $\mathfrak{g}_\Theta\subset \mathfrak{g}$ generated by all the rootspaces associated to $\Theta$, see \cite{GWpositivity22}. These Lie subalgebra were introduced by Bradlow, Collier Gothen and Garci-Prada as \emph{magical triples} in \cite{magical}. They proved the following. Let $\mathfrak{h}_\Theta$ be a $\Theta$-principal $\mathfrak{sl}_2$ Lie subalgebra of $\mathfrak{g}$.

\begin{thm}{\cite[Theorem 8.8]{magical}}
The exist a union of connected components of $\rho:\Gamma\to G$, the \emph{Cayley components}, consisting only of representations that do not factor through any parabolic subgroup. All $\mathfrak{h}_\Theta$-generalized Fuchsian representations with respect to the principal $\mathfrak{sl}_2$ Lie subalgebra $\mathfrak{h}_\Theta$ lie in some Cayley component.
\end{thm}

Cayley components are conjectured to be all the connected components of $\Theta$-positive representations \cite{magical}, but there exist components consisting of positive representations that do not contain $\mathfrak{h}_\Theta$-generalized representations, for instance the Gothen components for $G=\Sp(4,\R)$. The results of Guichard, Labourie implies the following \cite{positivityLabourie}. 

\begin{cor}
Every connected component of representations $\rho:\Gamma_g\to G$ containing a $\mathfrak{h}_\Theta$-generalized Fuchsian representation consist only of $\Theta$-positive representations. 
\end{cor}

The sets of simple roots $\Theta\subset \Delta$ that admit a notion of $\Theta$-positivity aways admit one or two subset which are Weyl orbits of simple roots, see Figure \ref{fig:TableWeylOrbit}. Let $\Theta'\subset \Theta$ be a Weyl orbit of simple roots. Let $G\neq \PSL(2, \R)$, Theorem \ref{thm:HighTeich} implies:

\begin{cor}
\label{cor:ThetaStruct}
Let $\rho:\Gamma_g\to G$ be a representation in a connected component of $\Theta$-positive representations containing a $\mathfrak{h}_\Theta$-generalized Fuchsian representation. It is the restricted holonomy of a $(G, \mathcal{F}_{\type_{\Theta'}})$-structure on a fiber bundle $F$ over $S_g$. 

\medskip

Let $\chi$ be the character associated to one of the $\mathfrak{h}_\Theta$-generalized Fuchsian representations in the Cayley component. This fiber bundle is diffeomorphic to the reduction of the $\SO(2,\R)\times \mathcal{C}_K(\mathfrak{h})$-bundle $P_{S_g,\chi}$ via its action on the base $\mathcal{B}^\type(\mathfrak{h})$ of the pencil of tangent vectors associated to $\mathfrak{h}$.
\end{cor}

The proof of the fact that the associated domains are non empty as soon as $G$ is not isomorphic to $\PSL(2, \R)$ is delayed to Section \ref{subsec:CheckNonEmpty}.

\subsubsection{Hitchin representations in \texorpdfstring{$\PSL(n,\R)$}{PSL(n,R)}.}
\label{sec:TeichSL3}
Let $\mathfrak{h}$ be a principal $\mathfrak{sl}_2$ Lie sugalgebra in $\mathfrak{sl}_n(\R)$. The associated representation $\iota_\mathfrak{h}$ is the representation $\iota_\text{irr}$ from Example \ref{exmp:PSL2N}.

\begin{defn}
A representation $\rho:\Gamma_g\to \PSL(n,\R)$ is \emph{Hitchin} if it is a deformation in $\Hom(\Gamma_g,G)$ of a $\mathfrak{h}$-generalized Fuchsian representaton.
\end{defn}

 The centralizer of $\mathfrak{h}$ in $\PSL(n,\R)$ is trivial, so $\mathfrak{h}$-generalized Fuchsian representations can be written $\iota_\text{irr}\circ \rho_0$.

\medskip

Hitchin proved that the quotient of the space of Hitchin representations by conjugation in $\PGL(n,\R)$ is a ball of dimention $(2g-2)(n^2-1)$ \cite{Hitchin}. Labourie proved that Hitchin representations are Borel Anosov, \ie $\Delta$-Anosov \cite{Labourie}. 

\medskip

The unique Weyl orbit of simple roots for $G=\PSL(n,\R)$ is $\Delta$. The flag manifold $\mathcal{F}_{\type_\Delta}$ can be identified with the flag manifold $\mathcal{F}_{1,{n-1}}$ consisting of 
pairs of subspaces $(\ell,H)$ where $\ell\subset H\subset \R^2$, $\dim(\ell=1)$, $\dim(H)=n-1$. The $\mathfrak{sl}_2$ Lie subalgebra $\mathfrak{h}$ is $\Delta$-regular. Theorem \ref{thm:HighTeich} implies the following:

\begin{cor}
Let $n\geq 3$, every Hitchin representation $\rho:\Gamma_g\to \PSL(n,\R)$ is the restricted holonomy of a $(\PSL(n,\R),\mathcal{F}_{1,n-1})$-structure on a fiber bundle over $S_g$.
\end{cor}

\medskip

When $n=3$, the boundary map of any $\mathfrak{h}$-generalized Fuhsian representation is an ellipsoid $\mathcal{E}\subset \mathbb{RP}^2$. The domain $\Omega^{\type_\Delta}_\rho\subset\mathcal{F}_{1,2}$ admits $3$ connected components: the set of $(\ell,H)\in \mathcal{F}_{1,2}$ with $\ell$ in the inside of the ellipsoid $\mathcal{E}$, with $H$ completely outside of the ellipsoid, and finally with $\ell$ outside the ellipsoid and $H$ crossing the ellipsoid in two points.
 One can see that the quotient of this domain is the union of three copies of the projectivization of the tangent bundle of $S_g$. If we apply Theorem \ref{thm:FoliationDoD} to $u_\mathfrak{h}$, we obtain a fibration where the model fiber $\mathcal{B}^\type(\mathfrak{h})$ is the union of $3$ circles described in Example \ref{exmp:PencilsRegular}. 
 Also when $n=3$ one can get a domain in projective space, as described in Example \ref{exmp:Nonregular} that is the interior of an ellipse.

\medskip

When $n$ is even, Hitchin representations are also the holonomy of projective structures. We consider these structures in Section \ref{appendix:CompPSL2N}.


\medskip

In general for any split simple Lie group $G$, Fock and Goncharov proved that all $\Delta$-positive representations can be deformed into a $\mathfrak{h}_\Delta$-generalized Fuchsian representation \cite{FockGoncharov}, \ie lie in a Hitchin component, so our method always applies.

\subsubsection{Maximal representations in \texorpdfstring{$\PSp(2n,\R)$}{PSp(2n,R)}.}
Given an orientation of the surface $S_g$, one can define the \emph{Toledo invariant} $\Tol:\Hom(\Gamma_g, \PSp(2n,\R))\to \mathbb{Z}$. This continuous map can de defined as the pullback by $\rho$ of an element of the continous group cohomology $H^2(G,\mathbb{Z})$ of $G$ by $\rho$ \cite{BurgerIozziWienhard}. Reversing the orientation of $S_g$ reverses the sign of the Toledo invariant.

\medskip

A representation $\rho:\Gamma_g\to \PSp(2n,\R)$ is called \emph{maximal} if its Toledo invariant is maximal among all representations, \ie if $\Tol(\rho)=n(2g-2)$. A way to construct maximal representation is to use the representation $\iota_\text{red}:\PSL(2,\R)\to \PSp(2n,\R)$. Burger, Iozzi, Labourie and Wienhard proved that maximal representations are $\lbrace \alpha_n\rbrace $-Anosov \cite{BurgerIozziWienhardLabourie}.

Let $\mathfrak{h}$ be the $\mathfrak{sl}_2$ Lie subalgebra of $\mathfrak{sp}_{2n}(\R)$ which is the image of $\mathrm{d}\iota_\mathfrak{h}$. Every $\mathrm{h}$-generalized representation is maximal for one of the two orientations of the surface $S_g$.

\begin{thm}[{\cite{GuichardWienhardMaximality}}]
If $n\geq 3$, every maximal representation $\rho:\Gamma_g\to \PSp(2n,\R)$ can be deformed into a $\mathfrak{h}_{\lbrace\alpha_n\rbrace}$-generalized Fuchsian representation in the sense of Definition \ref{defn:GeneralFuchs}.
\end{thm}
 Theorem \ref{thm:HighTeich} implies :

\begin{cor}
\label{cor:Maximal}
Let $n\geq 2$, every maximal representation $\rho:\Gamma_g\to \PSp(2n,\R)$ is the holonomy of a contact projective structure \ie a $\left(\PSp(2n,\R), \mathbb{RP}^{2n-1}\right)$-structure on a fiber bundle .
\end{cor}

Indeed one can consider the Weyl orbit of simple roots  $\Theta=\lbrace \alpha_n\rbrace $. The corresponding flag manifold $\mathcal{F}_{\tau_{\lbrace\alpha_n\rbrace}}$ can be identified with $\mathbb{RP}^{2n-1}$. The Lie subalgebra $\mathfrak{h}$ is $\lbrace\alpha_n\rbrace$-regular, so one can apply Theorem \ref{thm:HighTeich}.

\medskip

It is not clear if our method applies to the \emph{Gothen components} of representations $\rho:\Gamma_g\to \PSp(4, \R)$ that contain only Zariski-dense representations \cite{Gothen}. However since $\PSp(4,\R)\simeq \SO_o(2,3)$, the case $n=2$ of Corollary \ref{cor:Maximal} is a consequence of the work of Collier, Tholozan, Toulisse \cite{Collier_2019}. We discuss in more depth the relation between their consruction and our techniques in Appendix \ref{appendix:CCT}.

\medskip

The fiber obtained for Hitchin representations, that are also ${\lbrace\alpha_n\rbrace}$-positive, is a union of connected components of $\mathcal{B}^{\type_{\lbrace\alpha_n\rbrace}}(\mathfrak{h}_\Delta)$ and for maximal representations one gets a union of connected components of $\mathcal{B}^{\type_{\lbrace\alpha_n\rbrace}}(\mathfrak{h}_\Theta)$. These two submanifolds of $\mathbb{RP}^{2n-1}$ are diffeomorphic to the same Stiefel manifold, which is connected if $n\geq 3$, see the Appendix \ref{appendix:CompPSL2N}.

\subsubsection{Positive representations in \texorpdfstring{$\PSO(p,q)$}{PSO(p,q)}.}
\label{sec:PSO(p,q)}

Let $G=\PSO(p,q)$ with $q>p$ and $\Theta=\Delta\setminus \lbrace\alpha_p\rbrace$. This pair admits a notion of $\Theta$-positivity. The corresponding flag manifold $\mathcal{F}_{\tau_\Theta}$ can be identified with the Grassmanian of isotropic planes in $\R^{p,q}$. 

\medskip

Representations satisfying the $\Theta$-positive property were studied by Beyer and Pozzetti \cite{BeyerPozzetti}. In particular they show that all $\Theta$-positive representations $\rho:\Gamma_g\to \PSO(p,q)$ can be deformed to a $\mathfrak{h}_\Theta$-generalized Fuchsian representation when $q>p+1$, so in this case Corollary \ref{cor:ThetaStruct} applies to all $\Theta$-positive representations.
However when $q=p+1$, there are connected components of $\Theta$-positive representations that are conjectured to contain only Zariski-dense representations, it is not clear if our techniques can be applied to these components.

\subsection{Non-empty domains.}
\label{subsec:CheckNonEmpty}
In order to get a geometric structure associated to a domain of discontinuity, one need to ensure that the domain is non-empty. Kapovich, Leeb and Porti have a condition that ensures that there exist a thickening such that the domain is not empty \cite{KLP-DoD}, and Guichard and Wienhard proved that the domains they considered were not empty by computing the dimension of their complement. 

\medskip

We will use the following criterion to prove that some domains of discontinuity for surface groups are non-empty. Remember that the groups $\Gamma_g$ that we consider here are surface groups.

\begin{lem}
\label{lem:Non-emptyDomain}
Let $\rho:\Gamma_g\to G$ be a $\Theta$-Anosov representation and $(\type, \type_0)$ be a balanced type pair of elements in $\Sph\mathfrak{a}^+$, \ie such that $\typeS$ is $\type$-regular, satisfying $\Theta(\type_0)=\Theta$. Suppose that $\mathcal{F}_\type$ has finite fundamental group. The domain of discontinuity $\Omega_\rho^{(\type,\typeS)}$ is non-empty.
\end{lem}

 Note that it is however not a necessary condition to have a non-empty domain, as we see later for $\SO(2,3)$. This lemma is very similar to Proposition \ref{prop:NonEmptyPencils}.
 
\medskip 
 
\begin{proof}
 
If the domain is empty, it means that the flag manifold $\mathcal{F}_\type$ can be written as the union for $x\in \partial \Gamma_g$ of the thickenings $K^{\type_0}_{\xi_\rho^\Theta(x)}$. This uniojn is disjoint since the limit map is transverse. Moreover the limit map $\xi^\Theta_\rho$ is continuous, see \cite{BPS} for instance, which implies that $\mathcal{F}_\type$ fibers over the circle with a compact base. As in Proposition \ref{prop:NonEmptyPencils}, this implies that $\mathcal{F}_\type$ has infinite fundamental group.
\end{proof}

We prove that some domains of discontinuity are non-empty for representations of a surface group $\Gamma_g$. Let $(G, \Theta)$ be a pair that admits a notion of positivity and let $\Theta'\subset \Theta$ is a Weyl orbit of simple roots.

\begin{prop}
\label{prop:CheckNonEmptyPositivity}
The flag manifold $\mathcal{F}_{\type_{\Theta'}}$ has finite fundamental group, except if $G$ is locally isomorphic to $\PSL(2, \R)$ or if $G$ is locally isomorphic to $\SO_o(2,3)$, $\Theta=\Delta$, and $\Theta'=\lbrace\alpha_2\rbrace$.
\end{prop}

If $G$ is locally isomorphic to $\SO_o(2,3)$, $\Theta=\Delta$, and $\Theta'=\lbrace\alpha_2\rbrace$ we work by hand using the notations from Section \ref{subsec:ExamplesRoots}. The domain of discontinuity associated to a $\lbrace \alpha_2\rbrace$-Anosov representation $\rho:\Gamma_g\to \SO_o(2,3)$ obtained by metric thickening for any $\lbrace \alpha_2\rbrace$-regular $\typeS\in \Sph\mathfrak{a}^+$ is:
$$\Omega_\rho^{(\type_{\lbrace \alpha_2\rbrace},\typeS)}=\Ein(\R^{2,3})\setminus \bigcup_{x\in \partial\Gamma_g}\mathbb{P}\left(\xi^2_\rho(x)^\perp\right)^ .$$

For any $x\in \partial\Gamma_g$, The submanifold $\mathbb{P}\left(\xi^2_\rho(x)^\perp\right)$ has dimension $1$ in $\Ein(\R^{2,3})$ which has dimension $3$. Therefore the domain $\Omega_\rho^{(\type_{\lbrace \alpha_2\rbrace},\typeS)}$ is non-empty.

\medskip

Lemma \ref{lem:Non-emptyDomain} together with Proposition \ref{prop:CheckNonEmptyPositivity} implies the following.

\begin{cor}
\label{cor:CheckNonEmptyPositivity}
If $G$ is not locally isomorphic to $\PSL(2,\R)$, the domain of discontinuity obtained by metric thickening $\Omega_\rho^{(\type_{\Theta'},\typeS)}\subset \mathcal{F}_{\type_{\Theta'}}$ for any $\Theta'$-regular vector $\typeS\in \Sph\mathfrak{a}^+$ and any $\Theta$-Anosov representation $\rho:\Gamma_g\to G$ is non-empty.
\end{cor}

\medskip

The proof of Proposition \ref{prop:CheckNonEmptyPositivity} relies on a description of the fundamental groups of flag manifolds associated to real Lie groups.

\medskip

Let us write $\Delta_0\subset \Delta$ the set of roots whose associated root-space is a line. Let $\alpha,\beta\in \Delta$. We define $\epsilon(\alpha,\beta)=(-1)^{(\alpha, \beta^\vee)}$ \ie $\epsilon(\alpha, \beta)=1$ if $\alpha$ and $\beta$ are linked by no edge or if they are liked by two edges and $\alpha$ is the longest root, and else $\epsilon(\alpha, \beta)=-1$.

\begin{thm}{\cite[Theorem 1.1]{ThesisTopologyFlag}}
Let $A\subset \Delta$  be a set of simple roots.  The fundamental group of $\mathcal{F}_A=G/P_A$ is the group generated by $(t_\alpha)_{\alpha\in \Delta_0}$, defined by the relations $t_\beta t_\alpha=t_\alpha \left(t_\beta\right)^{\epsilon(\alpha, \beta)}$ for $\alpha,\beta\in \Delta_0$ and $\alpha\neq \beta$, and $t_\alpha=e$ if $\alpha\in \Delta_0\setminus A$. 
\end{thm}

The following lemma will deal with most cases.

\begin{lem}
\label{lem:DynkinFlag}
If the Dynkin Diagram restricted to $\Delta_0$ of the restricted root system of $G$ has no connected component of type $C_n$ or $A_1$, every flag manifold of $G$ has finite fundamental group.
\end{lem}

\begin{proof}
Let $A\subset \Delta$. Using the relations, one write any element of $\pi_1(\mathcal{F}_A)$ as a product of powers of generators so that each generator appears at most once. Therefore $\pi_1(\mathcal{F}_A)$ is a finite group if and only if for every $\alpha \in A$, $t_\alpha$ has finite order.

\medskip

If no connected component of the Dynkin Diagram restricted to $\Delta_0$ is of type $C_n$ or $A_1$, then every $\alpha\in \Delta_0$ belongs to a sub-diagram inside the Dynkin diagram of $\Delta_0$ of type $A_2$, $G_2$ or $B_3$. In each of this cases the relations between the generators imply that they have finite order: indeed the flag manifold associated to the Borel subgroup of $\SL(2, \R)$, $\SO(3,4)$ and the real split Lie group associated to $G_2$ have finite fundamental group, see \cite{ThesisTopologyFlag}.
\end{proof}

We now prove Proposition \ref{prop:CheckNonEmptyPositivity}. 

\medskip

\begin{proof}

If $G$ is a split Lie group $\Delta_0=\Delta$. If $G$ is locally isomorphic to $\SO(p,q)$ with $p\geq 3$ and $q>p+1$, the Dynkin diagram restricted to $\Delta_0$ is of type $A_{p-1}$. Finally if $G$ is the real forms of the complex Lie group associated to $E_6$, $E_7$ or $E_8$ whose restricted root system is of type $F_4$, the Dynkin diagram restricted to $\Delta_0$ is of type $A_2$. This follows from \cite[Table 9]{LieGroupTable}.

\medskip
 
 Therefore if $(G, \Theta)$ is a pair that admits a notion of $\Theta$-positivity Lemma \ref{lem:DynkinFlag} applies and every flag manifold associated to $G$ has finite fundamental group except in the following two cases: if $G$ is of Hermitian type and of tube type with $\Theta$ consisting of only the longest simple root, or if $G$ is split of type $A_1$ or $C_n$ and $\Theta=\Delta$.
 
\medskip

If $G$ is of Hermitian type the Dynkin diagram of the associated root system is $C_n$, $n\geq 2$. Suppose that $\Theta'=\lbrace\beta_n\rbrace$, then $\mathcal{F}_{\type_{\Theta'}}=\mathcal{F}_{\lbrace\beta_1\rbrace}$ as in Figure \ref{fig:TableWeylOrbit}. Either $\beta_1\notin \Delta_0$ in which case $\mathcal{F}_{\lbrace\beta_1\rbrace}$ is trivial, or $G$ is split, by \cite[Table 9]{LieGroupTable}.
If $G$ is split, $\beta_1,\beta_2\in \Delta_0$, so the generator $t_{\beta_1}$ of $\pi_1(\mathcal{F})$ satisfies the relation $t_{\beta_1} t_{\beta_2}=t_{\beta_2} \left(t_{\beta_1}\right)^{-1}$, and $t_{\beta_2}=e$ so $t_{\beta_1}^2=e$. Therefore $\mathcal{F}_{\type_{\Theta'}}$ has finite fundamental group.

\begin{center}
\dynkin[labels={\beta_1,\beta_2, \beta_{n-2}, \beta_{n-1}, \beta_n},
edge length=.75cm]C{}
\end{center}

It remains to consider the case where $G$ is split with root system $C_n$ for $n\geq 3$, $\Theta=\Delta$ and $\Theta'=\lbrace\alpha_1, \cdots ,\alpha_n\rbrace$. But in this case $\mathcal{F}_{\type_{\Theta'}}$ is the flag manifold associated to the root $\beta_2$, as shown in Figure \ref{fig:TableWeylOrbit}. The root $\beta_2$ belongs to a subdiagram of type $A_2$, so the fundamental group of $\mathcal{F}_{\Theta'}$ is finite.
\end{proof}

\subsection{Other applications.}
\label{subsec:HigherDimention}
Theorem \ref{thm:Application} can also be applied to Gromov hyperbolic groups that are not surface groups. In this subsection we consider representations of fundamental groups of hyperbolic manifolds.

\medskip

For instance one can consider a compact hyperbolic $3$-manifold $M$ with fundamental group $\Gamma$ and holonomy $\rho_0:\Gamma\to \PSL(2,\C)$. Let $n\geq 3$ and let $\iota_\text{irr}:\PSL(2,\C)\to \PSL(n,\C)$ be the irreducible representation as in Example \ref{exmp:PSL2N} and let $\mathfrak{h}=\mathrm{d}\iota_\text{irr}(\mathfrak{sl}_2(\C))$. As for the real case, $\mathfrak{h}$ is $\Delta$-regular, so in particular $\iota_\text{irr}\circ \rho_0$ is $\Delta$-Anosov. Here $\Delta$ is the only Weyl orbit of simple roots, and the corresponding flag manifold $\mathcal{F}_{\tau_\Delta}$ can be identified with the space $\mathcal{F}_{1,n-1}^\C$ of pairs $(\ell,H)$ where $\ell\subset H\subset \C^n$, $\ell$ is a line and $H$ a hyperplane. 

\medskip

The domain of discontinuity associated to a $\Delta$-Anosov representation $\rho$ for one and hence any balanced pair of the form $(\type_\Delta, \type_0)$ is the following, where the limit map of $\rho$ decomposes as $\xi_\rho^\Theta=\left(\xi^1_\rho,\cdots ,\xi^{n-1}_\rho\right)$:

$$\mathcal{F}_{1,n-1}^\C\setminus\bigcup_{x\in \partial \Gamma} \lbrace(\ell,H)|\exists 1\leq k\leq n-1 \,\text{s.t. }\ell\subset \xi^k_\rho\subset H \rbrace.$$

The topological dimension of the thickening $K_f^\type\subset \mathcal{F}_{1,n-1}^\C$ for any flag $f\in G/P_\Delta$ is the maximum for $1\leq k\leq n-1$ of the real dimension of $\lbrace (\ell,H)\in \mathcal{F}_{1,n-1}^\C|\ell\subset E\subset H\rbrace $ for some $E\subset \C^n$ of dimension $k$. The dimension of the thickening equals $2n-4$ and the dimension of the flag manifold equals $4n-6$ so for $n\geq 3$ the domain is non-empty. Theorem \ref{thm:Application} implies therefore the following:

\begin{cor}
The representation $\iota_\text{irr}\circ \rho_0:\Gamma \to \PSL(n,\C)$ is the restricted holonomy of a $(\PSL(n,\C),\mathcal{F}_{1,n-1}^\C)$-structure \emph{on a fiber bundle} over $M$.
\end{cor}



One can also consider a hyperbolic $n$-manifold $M$ for $n\geq 2$ with fundamental group $\Gamma$ and holonomy $\rho_0:\Gamma\to \SO_o(1,n)$. Let $\iota:\SO(1,n)\to \SO(p,np)$ be the diagonal representation for $n\geq 1$ and let $\mathfrak{h}$ be the image of $\mathrm{d}\iota$. Here $\mathfrak{h}$ is $\lbrace\alpha_p\rbrace$-regular, so in particular $\iota_\text{irr}\circ \rho_0$ is $\lbrace\alpha_p\rbrace$-Anosov. Note that $\lbrace\alpha_p\rbrace$ is a Weyl orbit of simple roots, and the corresponding flag manifold $\mathcal{F}_{\tau_{\lbrace\alpha_p\rbrace}}$ can be identified with the set of isotropic lines $\mathcal{I}\subset \mathbb{P}\left(\R^{p,np}\right)$.

\medskip

The domain of discontinuity associated to a $\lbrace \alpha_p\rbrace$-Anosov representation $\rho$ for one and hence any balanced pair of the form $(\type_{\lbrace \alpha_p\rbrace}, \type_0)$ is the following:

$$\mathcal{I}\setminus\bigcup_{x\in \partial \Gamma} \lbrace \ell|\ell\subset \xi_\rho^{\lbrace \alpha_p\rbrace}(x)\rbrace.$$

The dimension of the complement equals $(p-1)+n-1$ and the dimension of the flag manifold equals $n(p+1)-2$ so for $n\geq 2$ and $p\geq 2$ the domain is non-empty. Theorem \ref{thm:Application} implies therefore the following:

\begin{cor}
Let $\mathcal{C}$ be the connected component of $\iota\circ \rho_0$ in the space of $\lbrace\alpha_p\rbrace$-Anosov representations $\rho:\Gamma\to \PSL(n,\C)$. Every representation in $\mathcal{C}$ is the restricted holonomy of a $(\SO(p,pn),\mathcal{I})$-structure \emph{on a fiber bundle} over $M$.
\end{cor}

The fibers of this fiber bundle can be described as the set of isotropic lines in the intersection of $p$ quadrics in $\R^{p,pn}$.

\medskip

 The centralizer of $\iota\left(\SO(1,n-1)\right)$ in $\SO(p,pn)$ has a larger dimension than the centralizer of $\iota\left(\SO(1,n)\right)$. Indeed let $E\subset \R^{p,np}$ be the $p$-dimensional subspace preserved by $\iota\left(\SO(1,n-1)\right)$. Any element $g\in \SO(p,pn)$ acting trivially on $E^\perp$ centralizes $\iota\left(\SO(1,n-1)\right)$ but only finitely many such elements centralize $\iota\left(\SO(1,n)\right)$. Hence if the hyperbolic manifold contains a totally geodesic embedded hypersurface, there exist non-trivial deformations of $\iota\circ \rho$ that one can construct using bending.

\newpage 

\appendix
\section{Maximal representations of surface groups into \texorpdfstring{$\SO_o(2,n)$}{SO(2,n)}.}
\label{appendix:CCT}
In \cite{Collier_2019} Collier, Tholozan and Toulisse studied maximal representations into $\SO_o(2,n)$, the connected component of the identity of $\SO(2,n)$, for $n\geq 3$. In particular they showed that maximal representations are the restricted holonomies of photon structures which are geometric structures on a fiber bundle. They also show a similar result for Einstein structures and Hitchin representations in $\SO_o(2,3)\simeq \PSp(4,\R)$. In this appendix we compare the fibration of the domain of discontinuity that they obtain with the one coming from the projection defined in Section \ref{sec:Fibration}, when the unique minimal equivariant immersion satisfies the corresponding nearly geodesic property.

\medskip

Maximal representations are $\lbrace\alpha_1\rbrace$-Anosov, where $\lbrace\alpha_1\rbrace$ is the root such that $P_{\lbrace\alpha_1\rbrace}$ is the stabilizer of an isotropic line in $\R^{2,n+1}$. 
The set $\lbrace\alpha_1\rbrace\subset \Theta$ is a Weyl orbit of simple roots. We will denote by $\Omega^{\type_\Theta}_\rho$ this domain associated to $\rho$ which can be seen as an open subspace of the Grassmanian of isotropic planes on $\R^{2,n+1}$.

\medskip

Let $\R^{2,n+1}$ be a $n+3$-dimensional vector space equipped with a symmetric bilinear form $q$ of signature $(2,n+1)$. Let $\PSO_o(2,n)\subset \PSL(n+3,\R)$ be the connected component of the identity of the isometry group of $q$. Let $\mathbb{H}^{2,n}\subset \mathbb{P}\left(\R^{2,n+1}\right)$ be the space of \emph{timelike} lines, \ie lines on which $q$ is negative definite. The bilinear form $q$ induces a pseudo-Riemannian metric of signature $(2,n)$ on $\mathbb{H}^{2,n}$. 

For a subspace $V\subset \R^{2,n+1}$, we denote by $\Pho(V)$ the set of \emph{photons} of $V$, \ie the set of planes $P\subset V$ such that $q_{|P}=0$.  The flag manifold $\mathcal{F}_{\type_{\lbrace\alpha_2\rbrace}}$ can be identified with the space $\Pho(\R^{2,n+1})$. We denote by $\Ein(V)$ the space of isotropic lines in $V$. When $n=2$ the flag manifold $\mathcal{F}_{\type_{\lbrace\alpha_1\rbrace}}$ can be identified with the space $\Ein(\R^{2,3})$.

Fix an orientation of $S_g$. A maximal representation $\rho:\Gamma\to \PSO_o(2,n+1)$ is a representation whose Toledo invariant is maximal. Collier, Tholozan and Toullisse proved the following:

\begin{thm}[\cite{Collier_2019}]
\label{thm:CollierTholozanToulisse}
Let $\rho:\Gamma\to \SO_o(2,n)$ be a maximal representation. There exists a unique  $\rho$-equivariant maximal space-like embedding $\ell:\widetilde{S_g}\to \mathbb{H}^{2,n}$ and a unique minimal $\rho$-equivariant embedding $u:\widetilde{S_g}\to \X$. 
The following cocompact domain of discontinuity  fibers over the surface:
$$\Omega_\rho= \Pho(\R^{2,n})\setminus \bigcup_{x\in \partial \Gamma_g}\lbrace P\in  \Pho(\R^{2,n})|P\perp \xi^1_\rho(x)\rbrace.$$ 
The fibration $\pi:\Omega_\rho\to \widetilde{S_g}$ is such that for $x\in \widetilde{S}$: 
$$\pi^{-1}(x)=\Pho(\ell(x)^\perp).$$

\medskip

Let $\rho:\Gamma\to \SO_o(2,3)$ be a maximal Hitchin representation. There exists a unique  $\rho$-equivariant maximal space-like embedding $\ell:\widetilde{S_g}\to \mathbb{H}^{2,n}$ and a unique minimal $\rho$-equivariant embedding $u:\widetilde{S_g}\to \X$. 
The following cocompact domain of discontinuity  fibers over the surface:
$$\Omega_\rho= \Ein(\R^{2,n})\setminus\bigcup_{x\in \partial \Gamma_g}\lbrace \ell\in  \Pho(\R^{2,n})|\ell\perp \xi^2_\rho(x)\rbrace.$$ 
The fibration $\pi:\Omega_\rho\to \widetilde{S_g}$ is such that for $x\in \widetilde{S}$:
$$\pi^{-1}(x)=\Ein(u(x)\oplus \ell(x)).$$
\end{thm}

Theorem \ref{thm:ComparaisonPhotons} provides an other descitpion of the fibration from Theorem \ref{thm:CollierTholozanToulisse}: its purpose is to describe the difference with the one constructed in Theorem \ref{thm:FoliationDoD}.

\medskip

 The symmetric space $\X$ associated to $\PSO_o(2,n+1)$ can be identified with one connected component of the Grassmanian of oriented planes $U\subset \R^{2,n}$ that are \emph{spacelike}, \ie on which $q$ is positive. Moreover it of \emph{Hermitian type}, so it admits a complex structure $J$. 
 
 Let $\mathcal{P}$ be an oriented $2$-pencil of tangeant vectors based at $x\in \X$ \i.e. a subspace of dimension $2$ of $T_x\X$. It is  \emph{holomorphic}, respectively \emph{anti-holomorphic}, if $J_x$ fixes $P$ and for all non-zero $v\in P$, $\left(v,J_x(v)\right)$ is a positively oriented basis of $P$, respectively negatively oriented. 

To an oriented $2$-pencil $\mathcal{P}$  at $x\in \X$ that is not anti-holomorphic, with some orthonormal and oriented basis $(e,f)$, one can associate a holomorphic pencil:
 $$\phi(\mathcal{P})=\langle e-J_x(f),f+J_x(e)\rangle.$$
 
 Note that the pencil $\phi(\mathcal{P})$ does not depend on the positively oriented orthonormal basis $(e,f)$.  Similarly one can associate to a a $2$-pencil that is not holomorphic the anti-holomorphic pencil:
$$\overline{\phi}(\mathcal{P})=\langle e+J_x(f),f-J_x(e)\rangle.$$ 
 
 \begin{rem}
The pencils $\phi(\mathcal{P})$ and $\overline{\phi}(\mathcal{P})$
can be interpreted as the holomorphic and anti-holomorphic 
part of the complex $2$-dimensional subspace generated by $\mathcal{P}$ in $T_x\X$. In particular they are orthogonal to one another.
 \end{rem}

The fibers of the fibration from Theorem \ref{thm:CollierTholozanToulisse} can be described as bases of pencils of tangent vectors in $\X$. 

\begin{thm}
\label{thm:ComparaisonPhotons}
Let $\rho:\Gamma_g\to \PSO_o(2,n)$ be a maximal representation. Let $u:\widetilde{S_g}\to \X$ be the unique minimal $\rho$-equivariant embedding and $\ell:\widetilde{S_g}\to \mathbb{H}^{2,n}$ the unique maximal $\rho$-equivariant spacelike embedding. For all $x\in \widetilde{S_g}$ the submanifold:
$$\Pho(\ell(x)^\perp)\subset\Pho(\R^{2,n})\simeq \mathcal{F}_{\type_{\lbrace\alpha_1\rbrace}}$$
is equal to the $\type_{\lbrace\alpha_1\rbrace}$-base of the holomorphic pencil of tangent vectors:
$$\phi(\mathrm{d}u(T_x\widetilde{S_g})).$$

\medskip

Let $\rho:\Gamma_g\to \PSO_o(2,3)$ be a maximal Hitchin representation. Let $u:\widetilde{S_g}\to \X$ be the unique minimal $\rho$-equivariant embedding and $\ell:\widetilde{S_g}\to \mathbb{H}^{2,n}$ the unique maximal $\rho$-equivariant spacelike embedding. For all $x\in \widetilde{S_g}$ the submanifold:
$$\Ein(\ell(x)\oplus u(x))\subset\Ein(\R^{2,n})\simeq\mathcal{F}_{\type_{\lbrace\alpha_2\rbrace}}$$ 

is equal to the $\type_{\lbrace\alpha_2\rbrace}$-base of the anti-holomorphic pencil of tangent vectors :
$$\overline{\phi}(\mathrm{d}u(T_x\widetilde{S_g})).$$

\end{thm}

One can compare this fibration with the one provided by Theorem \ref{thm:FoliationDoD}, if $u$ is $\type$-nearly geodesic, for which the fiber at $y\in \widetilde{S_g}$ was the $\type$-base of the pencil $\mathrm{d}u(T_y\widetilde{S_g})$. In particular in the case of photon structures the two fibration coincide if $u$ is holomorphic, which happens when $\rho:\Gamma_g\to\PSO_o(2,n)$ factors through the reducible representation $\iota:\PSO_o(2,1)\to \PSO_o(2,n)$. However they differ in general and in particular if $\rho$ factors through the irreducible representation $\iota:\PSO_o(2,1)\to \PSO_o(2,3)$.

\medskip

From now on, this section is devoted to proving Theorem \ref{thm:ComparaisonPhotons}. We will use the notations from \cite[{Section 3.3}]{Collier_2019}. We first need a preliminary result about Hitchin representations.

\begin{lem}
\label{lem:ApeendixHitchinMax}
Let $\rho:\Gamma_g\to \SO_o(2,3)$ be a maximal and Hitchin representation. Its associated minimal immersion $u:\widetilde{S_g}\to \X$ is never tangent to a totally geodesic sumbanifold $\X_0\subset \X$ of the form $\X_0=\lbrace U\in \X|U\subset E\rbrace$ for any hyperplane $E\subset \R^{2,3}$.
\end{lem}

The proof of this fact relies on the parametrization of the minimal surface using Higgs bundles. We use the notations from \cite[Section 4.3]{Collier_2019}.

\medskip

\begin{proof}
If this was the case at a point $x\in \SO_o(2,3)$, it would imply that the associated Higgs field $\Phi_x$ at some point $x\in \widetilde{S_g}$ preserves the complexification of a real hyperplane of the fiber at $x$ of the Higgs bundle:
 $$E_\rho\simeq K^2\oplus K^1\oplus \mathcal{O}\oplus K^{-1}\oplus K^{-2}.$$
Recall that $\Phi$ is the following section of $K\otimes\End\left(E_\rho\right)$.

$$\Phi_x=\begin{pmatrix}
0& 0 & 0& \gamma& 0 \\
1 & 0 & 0& 0& \gamma  \\
0& \beta & 0&0 & 0  \\
 0& 0 &\beta&0 &0  \\
  0&0    & 0& 1&0 \end{pmatrix}.$$
  
Since $\Phi_x$ preserves a non-degenerated bilinear form and a real hyperplane, it must admit a real eigenvector $(z_1,z_2,z_3,z_4,z_5)$.
Note that the real from we consider on $K^2\oplus K^1\oplus \mathcal{O}\oplus K^{-1}\oplus K^{-2}$ is not the standard one, but the following
$(z_1,z_2,z_3,z_4,z_5)\mapsto (h_1^{-1}\overline{z_5},h_2^{-1}\overline{z_4},\overline{z_3},h_2\overline{z_2},h_1\overline{z_1})$ for some sections $h_1, h_2$ of $\overline{K^2}K^{-2}$ and  $\overline{K^1}K^{-1}$ respectively. Therefore a real eigenvector would be of the form $(z_1,z_2,z_3=\overline{z_3}, h_2\overline{z_2}, h_1\overline{z_1})$ and would satisfy for some $\lambda\in \C$:
$$\lambda z_1=\gamma(\partial_z) h_2\overline{z_2},$$
$$\lambda z_2=1(\partial_z)z_1+\gamma(\partial_z)h_1\overline{z_1},$$
$$\lambda z_3=\beta(\partial_z) z_2,$$
$$\lambda h_2\overline{z_2}=\beta(\partial_z) z_3 ,$$
$$\lambda h_1\overline{z_1}=1(\partial_z)h_2\overline{z_2} .$$

Moreover $\lVert\gamma\rVert<\lVert 1\rVert$ for the metric $h_1$ on $K^2$ and $\beta\neq 0$, see \cite[Section 4.3]{Collier_2019}. Therefore $\Phi_x$ cannot admit a non-zero real eigenvector. Indeed if $\lambda=0$, $z_2,z_3=0$ and $0=1(\partial_z)z_1+\gamma(\partial_z)h_1\overline{z_1}$ so $z_1=0$. If $\lambda\neq 0$, $1(\partial_z)h_2\overline{z_2}=h_1\gamma(\partial_z) h_2\overline{z_2}$ implies that $z_2=0$ and therefore $z_1,z_3=0$.

\end{proof}

We now prove Theorem \ref{thm:ComparaisonPhotons}.

\medskip

\begin{proof}
Let $\mathcal{G}(\mathbb{H}^{2,n})$ be the space of triples $(U,L,N)$ where $\R^{2,n}=U\oplus L  \oplus N$ is an orthogonal splitting, $L$ is a timelike line and $U$ a spacelike plane. The tangent space of $\mathcal{G}(\mathbb{H}^{2,n})$ can be decomposed as follows:
$$T_{(U,L,N)}\mathcal{G}\left(\mathbb{H}^{2,n}\right)\simeq \Hom\left(L,U\right)\oplus\Hom\left( L,N\right) \oplus \Hom\left(U,N\right).$$

The generalized Gauss map $\mathcal{G}:\widetilde{S} \to \mathcal{G}(\mathbb{H}^{2,n})$ associated to a spacelike immersion $\ell:\widetilde{S_g}\to \mathbb{H}^{2,n}$ is the unique map $x\mapsto (u(x),\ell(x),N(x))$ such that:
$$\mathrm{d}\mathcal{G}=(\mathrm{d}\ell ,0,A_\ell).$$
 Here $\mathrm{d}u$ is the differential of $\ell$:
$$\mathrm{d}\ell:T_x\widetilde{S_g}\to \Hom(\ell(x),u(x))\subset \Hom(\ell(x),\ell(x)^\perp)\simeq T_{\ell(x)}\mathbb{H}^{2,n}.$$
And $A_\ell$ is related to the second fundamental form $\II_\ell$ of the immersion $\ell$:
$$A_\ell:T_x\widetilde{S_g}\to \Hom(u(x),N(x))\subset \Hom(u(x),\ell(x)^\perp)\simeq T_{\ell(x)}\mathbb{H}^{2,n}\otimes\left( \mathrm{d}\ell \left(T_x\widetilde{S_g}\right)\right)^*,$$
$$\II_\ell(\mathrm{v},\mathrm{w})=A_\ell(\mathrm{v})\circ\mathrm{d}\ell(\mathrm{w})\in T_{\ell(x)}\mathbb{H}^{2,n}\simeq \Hom\left(\ell(x),\ell(x)^\perp\right).$$

If $\ell$ is the unique equivariant maximal spacelike surface for a maximal representation $\rho:\Gamma_g\to \PSO_o(2,n)$, the unique $\rho$-equivariant minimal surface $u:\widetilde{S_g}\to \X$ is the map coming from the Gauss map of $\ell$ \cite{Collier_2019}. The differential of $u$ in a direction $\mathrm{v}\in T_x\widetilde{S_g}$ is determined by $\mathrm{d}\mathcal{G}$, and is equal to the element $\mathrm{d}u(\mathrm{v})\in \Hom\left(u(x),u(x)^\perp\right)\simeq T_{u(x)}\X$ such that:
$$\mathrm{d}u(\mathrm{v})=A_\ell(\mathrm{v})-\transpose{\left(\mathrm{d}\ell(\mathrm{v})\right)}.$$

In this expression $\transpose{\left(\mathrm{d}\ell(\mathrm{v})\right)}\in \Hom\left(u(x),\ell(x)\right)$ is obtained as the adjoint of the endomorphism $\mathrm{d}\ell(\mathrm{v})$, for the restriction of the bilinear form $q$ on $\ell(x)$ and $u(x)$. 
\medskip

The endomorphism $J=J_{u(x)}$ coming from the complex structure on $\X$ acts on $T_{u(x)}\X\simeq \Hom\left(u(x),u(x)^\perp\right)$ by pre-composition by the rotation of $u(x)$ of angle $-\frac{\pi}{2}$, with the orientation chosen on $u(x)$. The choice of an orientation for every spacelike plane is equivalent to the choice of a complex structure on $\X$.
Since $\rho$ is maximal for the orientation we fixed on $S_g$,  this orientation is preserved by $\mathrm{d}u$. We write $R$ the corresponding rotation of $T_x\widetilde{S_g}$, which is a rotation of angle $-\frac{\pi}{2}$.

\medskip

Let $(\mathrm{v},\mathrm{w})$ be a positively oriented basis of $T_x\widetilde{S_g}$. The endomorphsim $f=\mathrm{d}u(\mathrm{v})-J(\mathrm{d}u(\mathrm{w}))\in \Hom(u(x),u(x)^\perp)\simeq T_{u(x)}\X$ satisfies :
$$f\left(\mathrm{d}\ell(\mathrm{v})\right)= \II_\ell(\mathrm{v},\mathrm{v})-\transpose{\left(\mathrm{d}\ell(\mathrm{v})\right)}\left(\mathrm{d}\ell(\mathrm{v})\right)- 
\II_\ell(\mathrm{w},R\cdot\mathrm{v})+\transpose{\left(\mathrm{d}\ell(\mathrm{w})\right)}\left(\mathrm{d}\ell\left(R\cdot\mathrm{v}\right)\right), $$
$$f\left(\mathrm{d}\ell(\mathrm{w})\right)= \II_\ell(\mathrm{v},\mathrm{w})-\transpose{\left(\mathrm{d}\ell(\mathrm{v})\right)}\left(\mathrm{d}\ell(\mathrm{w})\right)- 
\II_\ell(\mathrm{w},R\cdot\mathrm{w})+\transpose{\left(\mathrm{d}\ell(\mathrm{w})\right)}\left(\mathrm{d}\ell(R\cdot\mathrm{w})\right). $$

One has $R\cdot\mathrm{v}=-\mathrm{w}$, $R\cdot\mathrm{w}=\mathrm{v}$. Moreover $\II_\ell(\mathrm{v},\mathrm{w})=\II_\ell(\mathrm{w},\mathrm{v})$, and since $\ell$ is maximal $\II_\ell(\mathrm{v},\mathrm{v})+\II_\ell(\mathrm{w},\mathrm{w})=0$. So the image of $f$ is included in $\ell$, and $f\left(\mathrm{d}\ell(\mathrm{v})\right)\neq 0$.

\medskip

This holds for any orthonormal positively oriented basis  $(\mathrm{v},\mathrm{w})$, so $\phi(\mathrm{d}u(T_x\widetilde{S_g}))$ can be identified with the dimention $2$ subspace:
$$\Hom\left(u(x),\ell(x)\right)\subset\Hom\left(u(x),u(x)^\perp\right)\simeq T_{u(x)}\X.$$

It suffices now to show that for any $(U,L,N)\in \mathcal{G}(\mathbb{H}^{2,n})$, the $\type_{\lbrace\alpha_1\rbrace}$-base of the holomorphic pencil $\mathcal{P}\simeq\Hom\left(U,L\right)\subset\Hom(U, U^\perp)\simeq T_U\X$ is the photon $\Pho(L^\perp)$, independently of $U$.

\medskip

Let $a\in \Pho(\R^{2,n})$ be a photon. The vector $\mathrm{v}_{a,U}\in T_{U}\X$ pointing towards $a$ can be identified with the homomorphism $f_a:U\to U^\perp$ whose graph $\lbrace x+f_a(x)|x\in U\rbrace$ is the photon $a$.

\medskip

The Riemannian metric on $\X$ induces the following scalar product on $T_{u(x)}\X\simeq \Hom\left(U,U^\perp\right)$:
$$\langle f, g\rangle=\tr(g^\perp\circ f).$$

For this scalar product, $f_a$ is orthogonal to $\Hom\left(U,L\right)\subset \Hom\left(U, U^\perp\right)$ if and only if its image is orthogonal to $L$, which is orthogonal to $L$ if and only if $a\in \Pho(L^\perp)$.

\medskip

Now if $n=2$ and $\rho$ is maximal and Hitchin, the pencil $\overline{\phi}\left(\mathrm{d}u(T_x\widetilde{S_g})\right)$ is orthogonal to $\Hom\left(u(x),\ell(x)\right)\simeq\phi\left(\mathrm{d}u(T_x\widetilde{S_g})\right)$. But every unit vector $f\in\Hom\left(u(x),\ell(x)\right)$ is identified to $\mathrm{v}_{a,u(x)}\in T_{u(x)}\X$ where $a\in \Ein(\R^{2,n})$ is the graph of $\transpose{f}:\ell(x)\to u(x)$. Therefore $\Ein\left(u(x)\oplus \ell(x)\right)$ belongs to the $\type_{\lbrace \alpha_2\rbrace}$-base of $\overline{\phi}\left(\mathrm{d}u(T_x\widetilde{S_g})\right)$. It only remains to check that this base is not larger. This is due in some sense to the  \emph{genericity} of $\overline{\phi}\left(\mathrm{d}u(T_x\widetilde{S_g})\right)$.

\medskip

Suppose that there exist $\ell_0\in \Ein\left(\R^{2,n}\right)\setminus \Ein\left(u(x)\oplus \ell(x)\right)$ that belongs to the $\type_{\lbrace \alpha_2\rbrace}$-base of $\overline{\phi}\left(\mathrm{d}u(T_x\widetilde{S_g})\right)$. The associated element : 
$$\mathrm{v}_{\ell_0,u(x)}\in T_x\X\simeq \Hom(u(x), u(x)^\perp),$$
corresponds to a rank one  homomorphism $\mathrm{v}\in u(x)\mapsto \langle \mathrm{v},\mathrm{u}_0 \rangle \mathrm{n}_0$ where $\mathrm{u}_0\in u(x)$, $\mathrm{n}_0\in u(x)^\perp$ and $\mathrm{u}_0+\mathrm{n}_0$ generates $\ell_0$. Since $\ell_0$ belongs to the $\type_{\lbrace \alpha_2\rbrace}$-base of $\overline{\phi}\left(\mathrm{d}u(T_x\widetilde{S_g})\right)$ one has $\text{Im}(f)\perp \mathrm{n}_0$ for any : 
$$f\in \overline{\phi}\left(\mathrm{d}u(T_x\widetilde{S_g})\right)\subset T_x\X\simeq \Hom(u(x), u(x)^\perp).$$

\medskip

Let $\mathrm{n}_0'$ be a non-zero vector in $\left(u(x)\oplus \ell(x)\oplus \ell'(x)\right)\cap N(x)\subset \R^{2,3}$. It is orthogonal to the image of every element in $\overline{\phi}\left(\mathrm{d}u(T_x\widetilde{S_g})\right)$ since $\overline{\phi}\left(\mathrm{d}u(T_x\widetilde{S_g})\right)$ is orthogonal to $\phi\left(\mathrm{d}u(T_x\widetilde{S_g})\right)$. Moreover the image of every element in $\phi\left(\mathrm{d}u(T_x\widetilde{S_g})\right)$ is orthogonal to $N(x)$, so $\mathrm{n_0}'$ is orthogonal to the image of every element in $\mathrm{d}u(T_x\widetilde{S_g})\subset T_x\X\simeq \Hom\left(u(x),u(x)^\perp\right)$. This cannot happen if $\rho$ is Hitchin because of Lemma \ref{lem:ApeendixHitchinMax}. Therefore $\Ein\left(u(x)\oplus \ell(x)\right)$ is the $\type_{\lbrace \alpha_2\rbrace}$-base of $\overline{\phi}\left(\mathrm{d}u(T_x\widetilde{S_g})\right)$.

\end{proof}

\newpage

\section{Projective Structures with (Quasi-)Hitchin Holonomy.} 
\label{appendix:CompPSL2N}

In this appendix we compare the fibration obtained using Theorem \ref{thm:Application} with the one constructed in \cite{ADL}.

\medskip

Let $N\geq 3$ be an integer, and let $\mathfrak{h}$ be a principal $\mathfrak{sl}_2$-lie subalgebra of $\mathfrak{sl}(N,\R)$. The associated totally geodesic immersion $u_\mathfrak{h}:\mathbb{H}^2\to \X$ is $\type_1$-regular if and only if $N=2n$ is even. Recall that $\type_1\in\Sph\mathfrak{a}^+$ is such that $\mathcal{F}_{\type_1}$ is identifiable with $\mathbb{RP}^{2n-1}$.

\medskip

Any $\lbrace\alpha_n\rbrace$-Annsov representation $\rho$ admits a cocompact domain of discontinuity:
 $$\Omega_\rho=\mathbb{RP}^{2n-1}\setminus \bigcup_{x\in \partial\Gamma_g}\mathbb{P}\left(\xi_\rho^{\lbrace\alpha_n\rbrace}\right).$$
This is one of the domain constructed by Guichard and Wienhard \cite{GWDod}, and is also the domain $\Omega_\rho^{(\type_1,\typeS)}$ as in Theorem \ref{thm:KLPDoD} for any symmetric $\typeS\in \Sph\mathfrak{a}^+$ that is $\type_1$-regular. For any Fuchsian representation $\rho_0:\Gamma_g\to \PSL(2,\R)$ the domain associated with $\iota_\mathfrak{h}\circ \rho_0$ coincides with the domain $\Omega^{\type_1}_{u_\mathfrak{h}}$. This domain is non-empty for $n\geq 2$.

\medskip

To understand the topology of the quotient, it is sufficient to understand it for generalized Fuchsian representations, so let $u_\mathfrak{h}:\mathbb{H}^2\to \X$ be the totally geodesic embedding that is $\iota_\text{irr}$-equivariant, as defined in Example \ref{exmp:PSL2N}.

\medskip

We will compare two $2$-pencils of quadrics. For this we will use an invariant associated to a regular $2$-pencil of quadrics, which is called the Segre symbol \cite{fevola2020pencils}.

\medskip

A pencil of quadrics is \emph{regular} if it admits a non-degenerate element. Let $\mathcal{P}$ be a $2$-pencil with a basis $(q_1,q_2)$ such that $q_2$ is non-degenerate. There exists a unique $q_2$-symmetric or $q_2$-Hermitian endomorphism  $A_{q_1,q_2}$ of $V$ such that $q_1(\cdot,\cdot)=q_2(A\cdot,\cdot)$.

\begin{lem}

\label{lem:PSL2ActionPencilBasis}
Let $\mathcal{P}$ be a $2$-pencil of quadrics and let $(q_1,q_2)$ and $(q'_1,q'_2)$ be two basis of the pencil with $q_1,q_2,q'_1,q'_2$ non degenerate quadratic forms. Let $A=A_{q_1,q_2}$ and let $B=A_{q'_1,q'_2}$  and let $a,b,c,d\in \mathbb{R}$ be such that $q'_1=aq_1+bq_2$ and $q'_2=cq_1+dq_2$. Then $cA+d\text{Id}$ is invertible and $B=(aA+b\text{Id})(cA+d\text{Id})^{-1}$.
\end{lem}

\begin{proof}

Let  $M=aA+b\text{Id}$ and $N=cA+d\Id$. Then $q'_2(\cdot,\cdot)=q_2(N\cdot,\cdot)$. But since $q_2$ and $q'_2$ are non-degenerate, then $N=A_{q'_2,q_2}$ is invertible.

\medskip

In order to complete the proof of the lemma, we need to check that $q'_1(N\cdot,\cdot)=q'_2(M\cdot,\cdot)$. Let $v,w\in V$, then :
\begin{equation*} 
\begin{split}
q'_1(Nv,w) & = aq_1(Nv,w)+bq_2(Nv,w) \\
 & = adq_1(v,w)+bcq_2(Av,w)+acq_1(Av,w)+bdq_2(v,w)
\end{split}
\end{equation*}
\begin{equation*} 
\begin{split}
q'_2(Mv,w) & = cq_1(Mv,w)+dq_2(Mv,w) \\
 & = caq_1(Av,w)+cbq_1(v,w)+daq_2(Av,w)+dbq_2(v,w)
\end{split}
\end{equation*}

which concludes the proof since $q_1(\cdot,\cdot)=q_2(A\cdot,\cdot)$.
\end{proof}

We can associate to a $2$-pencil of quadrics an equivalence class of Segre symbols in the following way.

\begin{defn}
The \emph{Segre symbol} associated to a pair $(q_1,q_2)$ of non degenerate symmetric or Hermitian forms is the set $S$ of characteristic values of the endomorphism $A_{q_1,q_2}$, with the weights $w(s)$ given by the multiplicity of the characteristic values and the partition $p(s)$ given by the sizes of the blocks of the Jordan block decomposition of $A_{q_1,q_2}$.
\end{defn}

To a non-degenerate $2$-pencil of quadrics $\mathcal{P}$ we can associate the class of all Segre symbols associated to any basis $(q_1,q_2)$ of $\mathcal{P}$ where $q_1,q_2\in \mathcal{P}$ are non-degenerate, up to the action of $\PSL(2,\R)$ on the space of possible Segre symbols. This equivalence class is well defined because of Lemma \ref{lem:PSL2ActionPencilBasis}. The Segre symbol of a $2$-pencil is invariant by the action of $GL(V)$, hence two $2$-pencils of quadrics with different Segre symbols are not isomorphic.

\begin{thm}
\label{thm:ADL}
The fiber of the fibration $\pi_{u_\mathfrak{h}}:\Omega_\rho\to \widetilde{S}$ is the base of a non-degenerate $2$-pencil of quadrics $\mathcal{P}$ whose Segre symbol is equal to $\lbrace i, -i\rbrace$ with multiplicity $n$ at each point, and a minimal partition $[n]$ of $n$ at each point. 

\medskip

For $\K=\R$ it is diffeomorphic to the Stiefel manifold:
$$O(k)/\left(O(k-2)\times O(1)\right).$$ 

For $\K=\C$ it is diffeomorphic to the Stiefel manifold: $$O(2k)/\left(O(2k-2)\times U(1)\right).$$
\end{thm} 

A similar fibration was already constructed by Alessandrini, Li and the author in \cite{ADL}. The fibers in both constructions are bases of pencils, but in \cite{ADL} the fiber of the fibration constructed is the base of a different $2$-pencil $\mathcal{P}_0$. However the base of $\mathcal{P}_0$ is diffeomorphic to the same Stiefel manifold.

\medskip

The Segre symbol of the pencil $\mathcal{P}_0$ also has two conjugate points with multiplicity $n$, but decorated with a maximal partition $[1,1,\cdots,1]$ of $n$. A difference for instance is that the stabilizer in $\PSL_{2n}(\K)$ is larger for $\mathcal{P}_0$ than for $\mathcal{P}$.

\medskip

The rest of the section is devoted to the proof of Theorem \ref{thm:ADL}.

\medskip

\begin{proof}

As in example \ref{exmp:PSL2N}, we identify $\C^n$ with the space $V_{2n}(\C)$ of homogeneous polynomials in $X$ and $Y$ of degree $2n-1$. We consider the change of variable $X=W+Z$ and $Y=i(W-Z)$. If $\K=\R$, we recall that the conjugation $\tau$ after the change of coordinates is given by of the $\C$-anti-linear map defined for $0\leq \ell<2k$ by:
 $$\tau: W^\ell Z^{2k-1-\ell}\mapsto Z^\ell W^{2k-1-\ell}.$$ 
 This change of basis diagonalizes the action of $\PSO(2)$. Let $R(\theta)$ be the rotation angle $\theta$ in the hyperbolic plane, or the class of the rotation of angle $\frac{\theta}{2}$ on $\R^2$. The action of $R_\theta$ satisfies for $0\leq \ell<2k$:
  $$\iota\left(R\left(\theta\right)\right)\cdot W^\ell Z^{2k-1+\ell}= e^{i\theta(2k-1-2\ell)}W^\ell Z^{2k-1+\ell}.$$
 
Consider the following basis of $V_{2n}(\C)$:
 $$(e_\ell)_{0\leq \ell<2k}=\left(W^\ell Z^{2k-1+\ell} \binom{2k-1}{\ell}^{\frac{1}{2}}\right)_{0\leq \ell<2k}$$  

 The pencil of tangent vectors $\mathrm{d}u_\mathfrak{h}(\mathbb{H}^2)\in T_{q_0}\X$ corresponds to a pencils of (Hermitian) quadrics generated by $q_1$ and $q_2$. The corresponding matrices in this basis are:

$$Q_1=\begin{pmatrix}
0& \lambda_0 & & & && \\
\lambda_0 & 0 & \lambda_1& & &&  \\
& \lambda_1 & 0&\lambda_2 & &&  \\
 &    &\lambda_2&\ddots &\ddots  && \\
  &    &&\ddots&\ddots &\lambda_2  & \\
 &  & && \lambda_2& 0& \lambda_1 \\
 &  & && & \lambda_1& 0& \lambda_0\\
 &  & && & &\lambda_0&  0\end{pmatrix}$$$$Q_2=i\begin{pmatrix}
0& -\lambda_0 & & & && \\
\lambda_0 & 0 & -\lambda_1& & &&  \\
& \lambda_1 & 0&-\lambda_2 & &&  \\
 &    &\lambda_2&\ddots &\ddots  && \\
  &    &&\ddots&\ddots &-\lambda_2  & \\
 &  & && \lambda_2& 0& -\lambda_1 \\
 &  & && & \lambda_1& 0& -\lambda_0\\
 &  & && & &\lambda_0&  0\end{pmatrix}$$
 
 where $\lambda_a=\lambda_{2k-1-a}=\sqrt{(a+1)(2k-1-a)}$ for $0\leq a\leq 2k-1$.

\smallskip

The matrix associated to $A_{q_1,q_2}$ is conjugate to:

$$\begin{pmatrix}
iT& 0 \\
0 & -iT \end{pmatrix}
$$
 
 where $T$ is the following $k\times k$ matrix :
 $$\begin{pmatrix}
1& 1 &  \cdots &1 &1 \\
 & 1 & \cdots &1 & 1 \\
&  &\ddots & \vdots & \vdots \\
 &   &  &1 &1 \\
 &   &  & &1 \end{pmatrix}.$$
 
The matrix $T-\I_n$ is nilpotent of order $k$. This implies that the characteristic values of $A_{q_1,q_2}$ are $i$ and $-1$ and the partition associated with the Jordan block decomposition is the minimal partition $[k]$ for both singular values. This determines the Segre symbol of $\mathcal{P}$.

\medskip

Now let's study the topology of the base of $\mathcal{P}$ in the case $\K=\R$. in order to do so we will deform this pencil to get a pencil that can be decomposed into smaller blocks. Let us define for $t\in \mathbb{R}$ the pencil of quadrics $\mathcal{P}_t$ generated by the quadratic forms $q_{1,t},q_{2,t}$ whose matrices $Q_{1,t},Q_{2,t}$ in the basis $(e_\ell)_{0\leq \ell<2k}$ are:

$$Q_{1,t}=\begin{pmatrix}
0& \lambda_0 & & & && \\
\lambda_0 & 0 & t\lambda_1& & &&  \\
& t\lambda_1 & 0&\lambda_2 & &&  \\
 &    &\lambda_2&\ddots &\ddots  && \\
  &    &&\ddots&\ddots &\lambda_2  & \\
 &  & && \lambda_2& 0& t\lambda_1 \\
 &  & && & t\lambda_1& 0& \lambda_0\\
 &  & && & &\lambda_0&  0\end{pmatrix}$$$$Q_{2,t}=i\begin{pmatrix}
0& -\lambda_0 & & & && \\
\lambda_0 & 0 & -t\lambda_1& & &&  \\
& t\lambda_1 & 0&-\lambda_2 & &&  \\
 &    &\lambda_2&\ddots &\ddots  && \\
  &    &&\ddots&\ddots &-\lambda_2  & \\
 &  & && \lambda_2& 0& -t\lambda_1 \\
 &  & && & t\lambda_1& 0& -\lambda_0\\
 &  & && & &\lambda_0&  0\end{pmatrix}.$$
 
In other words the coefficients $\lambda_k$ for $k$ odd are replaced by $t\lambda_k$.
The determinant of $Q_{1,t}$ and $Q_{2,t}$ is independent of $t$, hence the pencils of quadrics $\mathcal{P}_t$ are all non-degenerate. Moreover one has 
$$\iota\left(R\left(\theta\right)\right)\cdot Q_{1,t}=\cos(\theta)Q_{1,t}+\sin(\theta)Q_{2,t}$$

which implies that all the non-zero quadrics in the pencil $\mathcal{P}_t$ are non-degenerate. 

\medskip

Let $M=\lbrace([v],t)\in \mathbb{P}(V_n(\R))\times \R|q_{1,t}(v,v)=q_{2,t}(v,v)=0\rbrace$ and $p:\mathbb{P}(V_n(\K))\times \R \to \R$ be the projection onto the second factor. Because of Lemma \ref{lem:SubmertionPencils}, $M$ is a sub-manifold of  $\mathbb{P}(V_n(\R))\times \R$. The map $p$ is a proper submersion, hence by the Ehresmann principle it is a fibration. In particular the bases $b(\mathcal{P})$ and $b(\mathcal{P}_0)$ of these pencils are diffeomorphic.

\bigskip

The pencil of quadrics $\mathcal{P}_0$ is exactly the pencil that defines the Stiefel manifold in \cite{ADL}, so we get the desired result.
 
\end{proof}

\newpage
\bibliographystyle{hamsalpha}
\bibliography{biblio}
\end{document}